%% file: kuenneth.tex
\documentclass[reqno,oneside]{amsart} 

\usepackage{amsthm,
amssymb,amsmath,enumerate,stackrel,
stmaryrd, 
leftidx,
ifthen}
\usepackage{xcolor} 
\usepackage{xr-hyper}
\usepackage[all,cmtip]{xy}
\usepackage[utf8]{inputenc} 
\chardef\cprime"7E 
\usepackage{cases}

\usepackage[pagebackref,hyperindex,linktocpage=true]{hyperref}
\hypersetup{
    colorlinks,
    linkcolor={red!50!black},
    citecolor={red!50!black}, 
    urlcolor={red!50!black}, 
    filecolor={red!50!black} 
}

\usepackage{theoremref}
\definecolor{labelkey}{rgb}{1,0,0}

\makeatletter
\@addtoreset{equation}{section}
\makeatother

\input{constructible_citations} 
\def\consref#1{\cite[\myconsref{#1}]{HemoRicharzScholbach:Constructible}}
\def\consrefsect#1{\cite[\S\myconsref{sect--#1}]{HemoRicharzScholbach:Constructible}}

\def\conseqref#1{\cite[(\myconsref{#1})]{HemoRicharzScholbach:Constructible}}

\numberwithin{equation}{section}

\theoremstyle{definition}
\newtheorem{Defi}{Definition}[section] \newcommand{\defi}{\begin{Defi}} \newcommand{\xdefi}{\end{Defi}} 
\newtheorem{DefiLemm}[Defi]{Definition and Lemma} \newcommand{\defilemm}{\begin{DefiLemm}} \newcommand{\xdefilemm}{\end{DefiLemm}} 
\newtheorem{Bsp}[Defi]{Example} \newcommand{\exam}{\begin{Bsp}} \newcommand{\xexam}{\end{Bsp}} 
\newtheorem{Syno}[Defi]{Synopsis} \newcommand{\syno}{\begin{Syno}} \newcommand{\xsyno}{\end{Syno}} 
\newtheorem{Bem}[Defi]{Remark} \newcommand{\rema}{\begin{Bem}} \newcommand{\xrema}{\end{Bem}} 
\newtheorem{Notation}[Defi]{Notation} \newcommand{\nota}{\begin{Notation}} \newcommand{\xnota}{\end{Notation}} 
\newtheorem{Warning}[Defi]{Warning} \newcommand{\warn}{\begin{Warning}} \newcommand{\xwarn}{\end{Warning}} 
\newtheorem{Situation}[Defi]{Situation} \newcommand{\situ}{\begin{Situation}} \newcommand{\xsitu}{\end{Situation}} 

\theoremstyle{plain}
\newtheorem{Theo}[Defi]{Theorem} \newcommand{\theo}{\begin{Theo}} \newcommand{\xtheo}{\end{Theo}} 
\newtheorem{Satz}[Defi]{Proposition} \newcommand{\prop}{\begin{Satz}} \newcommand{\xprop}{\end{Satz}} 
\newtheorem{Lemm}[Defi]{Lemma} \newcommand{\lemm}{\begin{Lemm}} \newcommand{\xlemm}{\end{Lemm}} 
\newtheorem{Coro}[Defi]{Corollary} \newcommand{\coro}{\begin{Coro}} \newcommand{\xcoro}{\end{Coro}}
\newtheorem{Ques}[Defi]{Question} \newcommand{\ques}{\begin{Ques}} \newcommand{\xques}{\end{Ques}}
\newtheorem{Conj}[Defi]{Conjecture} \newcommand{\conj}{\begin{Conj}} \newcommand{\xconj}{\end{Conj}}

\newcommand{\refsect}[1]{Section \ref{sect--#1}}
\newcommand{\refit}[1]{(\ref{item--#1})}
\newcommand{\refeq}[1]{(\ref{eqn--#1})}

\newcommand{\eqn}{\begin{equation}} \newcommand{\xeqn}{\end{equation}}
\newcommand{\eqnarr}{\begin{eqnarray*}} \newcommand{\xeqnarr}{\end{eqnarray*}}
\newcommand{\eqnarra}{\begin{eqnarray}} \newcommand{\xeqnarra}{\end{eqnarray}}

\newcommand{\pf}{\begin{proof}} \newcommand{\xpf}{\end{proof}}

\usepackage[top=1.5cm, bottom=1.5cm, left=2cm, right=2cm, heightrounded, marginparwidth=3.5cm, marginparsep=.2cm]{geometry}

\usepackage{color}

\newcommand{\nc}{\newcommand}
\nc{\StP}[1]{\cite[Tag~\href{http://stacks.math.columbia.edu/tag/#1}{#1}]{StacksProject}}
\nc{\StPd}[2]{\cite[Tags~\href{http://stacks.math.columbia.edu/tag/#1}{#1}, \href{http://stacks.math.columbia.edu/tag/#2}{#2}]{StacksProject}} 

\nc{\on}{\operatorname}
\nc{\aff}{{\on{aff}}}
\nc{\modi}{{\on{mod}}} 
\nc{\even}{{\on{even}}}
\nc{\odd}{{\on{odd}}}
\nc{\naive}{{\on{naive}}}
\nc{\hofib}{\on{hofib}}
\nc{\Bun}{\on{Bun}}
\nc{\ad}{{\on{ad}}}
\nc{\lft}{{\on{lft}}}

\nc{\Weil}{{\on{Weil}}} 
\nc{\FWeil}{{\on{FWeil}}} 
\nc{\cons}{{\on{cons}}} 
\nc{\tot}{{\on{Tot}}} 

\nc{\str}{\on{-}}
\nc{\perf}{{\on{perf}}}
\nc{\Rel}{{\on{Pos}}}
\nc{\lan}{\langle}
\nc{\ran}{\rangle}

\nc{\bbA}{{\mathbb A}} 
\nc{\bbB}{{\mathbb B}}
\nc{\bbC}{{\mathbb C}}
\nc{\bbD}{{\mathbb D}}
\nc{\bbE}{{\mathbb E}}
\nc{\bbF}{{\mathbb F}}
\nc{\bbG}{{\mathbf G}}
\nc{\bbH}{{\mathbb H}}
\nc{\bbI}{{\mathbb I}}
\nc{\bbJ}{{\mathbb J}}
\nc{\bbK}{{\mathbb K}}
\nc{\bbL}{{\mathbb L}}
\nc{\bbM}{{\mathbb M}}
\nc{\bbN}{{\N}} 
\nc{\bbO}{{\mathbb O}}
\nc{\bbP}{{\mathbb P}} 
\nc{\bbQ}{{\mathbb Q}} 
\nc{\bbR}{{\mathbb R}}
\nc{\bbS}{{\mathbb S}}
\nc{\bbT}{{\mathbb T}}
\nc{\bbU}{{\mathbb U}}
\nc{\bbV}{{\mathbb V}}
\nc{\bbW}{{\mathbb W}}
\nc{\bbX}{{\mathbb X}}
\nc{\bbY}{{\mathbb Y}}
\nc{\bbZ}{{\mathbb Z}}

\nc{\calA}{{\mathcal A}}
\nc{\calB}{{\mathcal B}}
\nc{\calC}{{\mathcal C}}
\nc{\calD}{{\mathcal D}}
\nc{\calE}{{\mathcal E}}
\nc{\calF}{{\mathcal F}}
\nc{\calG}{{\mathcal G}}
\nc{\calH}{{\mathcal H}}
\nc{\calI}{{\mathcal I}}
\nc{\calJ}{{\mathcal J}}
\nc{\calK}{{\mathcal K}}
\nc{\calL}{{\mathcal L}}
\nc{\calM}{{\mathcal M}}
\nc{\calN}{{\mathcal N}}
\nc{\calO}{{\mathcal O}}
\nc{\calP}{{\mathcal P}}
\nc{\calQ}{{\mathcal Q}}
\nc{\calR}{{\mathcal R}}
\nc{\calS}{{\mathcal S}}
\nc{\calT}{{\mathcal T}}
\nc{\calU}{{\mathcal U}}
\nc{\calV}{{\mathcal V}}
\nc{\calW}{{\mathcal W}}
\nc{\calX}{{\mathcal X}}
\nc{\calY}{{\mathcal Y}}
\nc{\calZ}{{\mathcal Z}}

\nc{\Sht}{{\on{Sht}}}
\nc{\Frob}{{\on{Frob}}}
\nc{\Hecke}{{\on{Hecke}}}
\nc{\inv}{{\on{inv}}}
\nc{\Conv}{{\on{Conv}}}
\nc{\triv}{{\on{triv}}}
\nc{\Isom}{{\on{Isom}}}

\nc{\scrB}{{\mathscr{B}}}
\nc{\scrA}{{\mathscr{A}}}
\nc{\bbf}{{\mathbf{f}}}
\nc{\bba}{{\mathbf{a}}}
\nc{\rig}{{\mathrm rig}}

\nc{\Indft}{{\pi_1\!\on{-}\!\mathrm{Indft}}}
\nc{\IndPerf}{{\pi_1\!\on{-}\!\mathrm{IndPerf}}}
\nc{\AS}{{\on{AS}}}

\nc{\al}{\alpha}
\nc{\be}{\beta}
\nc{\ga}{\gamma}
\nc{\la}{\lambda}
\nc{\qcqs}{{\on{qcqs}}}
\nc{\Bmu}{{\boldsymbol \mu}}

\nc{\pot}[1]{ [\hspace{-0,5mm}[ {#1} ]\hspace{-0,5mm}] }
\nc{\rpot}[1]{ (\hspace{-0,7mm}( {#1} )\hspace{-0,7mm}) }

\nc{\defined}{\hspace{0.1cm}\stackrel{\text{\tiny \rm def}}{=}\hspace{0.1cm}}
\nc{\co}{\colon}
\nc{\specto}{{\leadsto}}

\newcommand{\category}[1]{\mathrm{#1}}
\newcommand{\FinSet}{\category{FinSet}} 
\newcommand{\Cat}{\category{Cat}} 
\newcommand{\Ex}{\category{Ex}} 
\newcommand{\CatEx}{\Cat^{\Ex}} 
\newcommand{\Fun}{\category{Fun}} 
\newcommand{\PreStk}{\category{PreStk}} 
\newcommand{\Ani}{\category{Ani}} 
\newcommand{\cts}{\mathrm{cts}}  
\newcommand{\Sch}{\category{Sch}} 
\renewcommand{\Pr}{\category{Pr}}
\newcommand{\PrL}{\Pr^\category{L}} 
\newcommand{\PrSt}{\Pr^{\category{St}}} 
\newcommand{\PrStL}{\PrSt_\Lambda} 
\newcommand{\Ind}{\category{Ind}} 
\newcommand{\Mod}{\category{Mod}} 
\newcommand{\Sp}{\category{Sp}} 
\newcommand{\Perf}{\category{Perf}} 
\newcommand{\Fix}{\mathrm {Fix}} 
\newcommand{\Idem}{\mathrm {Idem}} 

\newcommand{\Dcons}{\mathrm{D}_\cons}
\newcommand{\Dindcons}{\mathrm{D}_\indcons}
\newcommand{\Dlis}{\mathrm{D}_\lis}
\newcommand{\Dindlis}{\mathrm{D}_\indlis}

\newcommand{\Cht}{\mathrm{Cht}}
\newcommand{\Schfp}{\mathrm{Sch}^{\mathrm{fp}}}

\def\GL{\mathrm {GL}} 
\def\fp{\mathrm{fp}} 

\newcommand{\colim}{\operatornamewithlimits{colim}} 
\def\id{{\rm id}} 
\def\ev{{\operatorname {ev}}} 
\def\opp{{\rm op}} 
\def\To#1#2{\mathop{\count0=#1 \loop\ifnum\count0>0 \smash-\mkern-7mu \advance\count0 -1 \repeat \mathord\rightarrow}\limits^{#2}} 
\def\Hom{\mathop{\rm Hom}\nolimits} 
\def\Frac{\mathop{\rm Frac}} 

\def\Ind{\category{Ind}} 
\def\RHom{\mathop{\rm RHom}\nolimits} 
\def\Sht{\mathop{\rm Sht}\nolimits} 
\def\CAlg{\mathop{\rm CAlg}\nolimits} 
\def\End{\mathop{\rm End}\nolimits} 
\def\Rep{\category{Rep}} 
\def\et{\mathrm{\acute et}} 
\def\proet{\mathrm{pro\acute et}} 

\definecolor{hellgrau}{RGB}{200,200,200} 
\definecolor{dunkelgrau}{RGB}{160,160,160} 
\definecolor{hellblau}{RGB}{194, 215, 249} %
\definecolor{dunkelblau}{RGB}{68, 128, 226} %

\def\Z{{\bbZ}} 
\def\Fp{{{\mathbb F}_p}} %
\def\Fq{{{\mathbb F}_q}} %
\def\N{{\mathbb N}} 
\def\Q{{\bbQ}} 
\def\Qp{\Q_p} 
\def\Ql{{\Q_\ell}} 

\def\H{{\rm H}} 
\def\DM{\category{DM}} 
\def\Vect{\category{Vect}} 
\def\fd{\category{fd}} 
\def\ii{$\infty$}

\def\Gal{{\rm Gal}} 

\def\disc{{\rm disc}} 
\def\lis{{\rm lis}} 
\def\indlis{{\rm indlis}} 
\def\indcons{{\rm indcons}} 

\def\RG{\R \Gamma} 

\def\Spec{\mathop{\rm Spec}} 
\def\B{{\rm B}} 
\newcommand{\D}{\category{D}} 

\def\R{{\rm R}} 
\def\sbuildrel#1\over#2{\mathrel{\smash{\mathop{\kern0pt #2}\limits^{#1}}}}

\let\x\times
\let\ol\overline
\renewcommand{\t}{\otimes}

\newcommand{\bx}{\boxtimes}

\renewcommand{\r}{\rightarrow}
\newcommand{\lr}{\longrightarrow}
\newcommand{\hr}{\hookrightarrow}

\def\matrix#1{\null\,\vcenter{\normalbaselines
    \ialign{\hfil$##$\hfil&&\quad\hfil$##$\hfil\crcr
      \mathstrut\crcr\noalign{\kern-\baselineskip}
      #1\crcr\mathstrut\crcr\noalign{\kern-\baselineskip}}}\,}

\newdimen\harrowsize
\def\mapright#1{\smash{\mathop{\hbox to\harrowsize{\rightarrowfill}}\limits^{#1}}}

{\catcode`@=11
\gdef\cal{\fam\tw@}
\global\let\over\@@over
\global\let\atop\@@atop
\global\let\above\@@above
\global\let\overwithdelims\@@overwithdelims
\global\let\atopwithdelims\@@atopwithdelims
\global\let\abovewithdelims\@@abovewithdelims
\gdef\eqalign#1{\null\,\vcenter{\openup\jot\m@th
\ialign{\strut\hfil$\displaystyle{##}$&$\displaystyle{{}##}$\hfil
      \crcr#1\crcr}}\,}
\newskip\xcentering \global\xcentering=0pt plus 1000pt minus 1000pt
\gdef\eqalignno#1{\displ@y \tabskip\xcentering
  \halign to\displaywidth{\hfil$\@lign\displaystyle{##}$\tabskip\z@skip
    &$\@lign\displaystyle{{}##}$\hfil\tabskip\xcentering
    &\llap{$\@lign##$}\tabskip\z@skip\crcr
    #1\crcr}}
\global\def\cases#1{\left\{\,\vcenter{\normalbaselines\m@th
    \ialign{$##\hfil$&\quad##\hfil\crcr#1\crcr}}\right.}
\gdef\eqlabel#1{\refstepcounter{equation}\label{eqn--#1}\eqno\hbox{\@eqnnum}}
}

\begin{document}

\author{Tamir Hemo, Timo Richarz and Jakob Scholbach}

\title{A categorical K\"unneth formula for constructible Weil sheaves}

\thanks{*The second named author T.R.~is funded by the European Research Council (ERC) under Horizon Europe (grant agreement nº 101040935), by the Deutsche Forschungsgemeinschaft (DFG, German Research Foundation) TRR 326 \textit{Geometry and Arithmetic of Uniformized Structures}, project number 444845124 and the LOEWE professorship in Algebra.
J.S. was supported by Deutsche Forschungsgemeinschaft (DFG), EXC 2044–390685587, Mathematik Münster: Dynamik–Geometrie–Struktur. }

\address{Caltech, Department of Mathematics, 91125 Pasadena, CA, USA}
\email{themo@caltech.edu}

\address{Technical University of Darmstadt, Department of Mathematics, 64289 Darmstadt, Germany}
\email{richarz@mathematik.tu-darmstadt.de}

\address{Università degli Studi di Padova, Dipartimento di Matematica, 35139 Padova, Italy} 

\begin{abstract}
We prove a Künneth-type equivalence of derived categories of lisse and constructible Weil sheaves on schemes in characteristic $p>0$ for various coefficients, including finite discrete rings, algebraic field extensions $E\supset \bbQ_\ell$, $\ell\neq p$ and their rings of integers $\calO_E$. We also consider a variant for ind-construtible sheaves which applies to the cohomology of moduli stacks of shtukas over global function fields.



\end{abstract}

\maketitle

\tableofcontents

\section{Introduction}

\label{sect--categorical.kuenneth.sect}

The classical Künneth formula expresses the (co-)homology of a product of two spaces $X_1$ and $X_2$ in terms of the tensor product of the (co-)homology of the individual factors.
For two topological spaces, for example, one has under suitable finiteness hypothesis an isomorphism
\begin{equation}
\label{coho.kuenneth.singular.eq}
\bigoplus_{i+j = n} \H^{i}(X_1, \Q) \t_\Q \H^{j}(X_2, \Q) \cong \H^{n}(X_1 \x X_2, \Q)
\end{equation}
on singular cohomology with rational coefficients.
Such cohomology groups are naturally morphism groups in the derived categories of sheaves on these spaces.
So one may ask whether the K\"unneth formula can be extended to a categorical level, that is, whether it is possible to relate the derived categories of sheaves on $X_1$ and $X_2$ to those on their product $X_1 \x X_2$.
Statements in this direction are referred to as \textit{categorical K\"unneth formulas} and are known in different contexts: for example, for the respective derived categories of topological sheaves, for D-modules on varieties in characteristic $0$ and for quasi-coherent sheaves, see \cite[Section A.2]{GKRV.Toy.Model}.

In addition to \eqref{coho.kuenneth.singular.eq} above, categorical Künneth formulas require decomposing a sheaf on $X_1 \x X_2$ into exterior products $M_1 \boxtimes M_2$, with $M_1$, $M_2$ being sheaves on $X_1$, $X_2$, respectively.
For varieties in characteristic $p>0$, an analogous decomposition for constructible (pro-)étale sheaves fails in general, and so does a categorical K\"unneth formula in this context, see \thref{AS.exam.intro} below.
The main result of the manuscript at hand (see \thref{derived_Drinfeld}) shows how to rectify the failure by adding equivariance data under partial Frobenius morphisms, that is, one arrives at a categorical K\"unneth formula for constructible Weil sheaves. 
Our work relies on the analogous result \cite[Theorem 2.1]{Drinfeld:Langlands} for étale fundamental groups known as Drinfeld's lemma, see \refsect{drinfelds.lemma} for details and references.

\subsection{Definitions and results}
Weil sheaves are defined in \cite[Definition 1.1.10]{Deligne:Weil2}.
We start by explaining a site-theoretic approach which slightly differs from \cite{Geisser:Weil, Lichtenbaum:Weil}.

Let $X$ be a scheme over a finite field $\bbF_q$, where $q$ is a $p$-power.
Fix an algebraic closure $\bbF/\bbF_q$, and denote by $X_\bbF$ the base change. 
The partial ($q$-)Frobenius $\phi_X:=\Frob_X\x \id_\bbF$ defines an endomorphism of $X_\bbF$. 

\defi 
The \textit{Weil-pro\'{e}tale site} $X^\Weil_\proet$ is the following site: 
Objects are pairs $(U,\varphi)$ consisting of $U\in (X_{\bbF})_{\proet}$, the proétale site of $X_\bbF$ \cite{BhattScholze:ProEtale}, equipped with an endomorphism $\varphi\colon U\rightarrow U$ of $\bbF$-schemes covering $\phi_X$. 
Morphisms are given by equivariant maps. 
A family $\{(U_i,\varphi_i) \rightarrow (U,\varphi)\}$ of morphisms is a cover if the family $\{U_i \rightarrow U\}$ is a cover in $(X_{\bbF})_{\proet}$.
\xdefi

The Weil-pro\'etale site sits in the sequence of sites
\begin{equation}
\label{Weil.sites.maps.eq}
(X_\bbF)_\proet\r X^\Weil_\proet \r X_\proet
\end{equation}
given by the functors $U \mapsfrom (U, \varphi)$ and $(U_\bbF,\phi_U)\mapsfrom U$ in the opposite direction.
The maps \eqref{Weil.sites.maps.eq} commute over $*_\proet$, the proétale site of the point. 
Thus, for any condensed ring $\Lambda$ viewed as a sheaf of rings on $*_\proet$, we get pullback functors on derived categories of proétale $\Lambda$-sheaves
\[
\D(X,\Lambda)\r \D\big(X^\Weil,\Lambda\big) \r \D(X_\bbF,\Lambda).
\]
In analogy with the definition of lisse and constructible sheaves (as recalled in \thref{recall.definitions}), we introduce the categories of lisse and constructible Weil sheaves $\D_\lis\big(X^\Weil,\Lambda\big)\subset \D_\cons\big(X^\Weil,\Lambda\big)$ as the full subcategories of $\D\big(X^\Weil,\Lambda\big)$ that are dualizable, resp.~that are Zariski locally on $X$ dualizable along a constructible stratification. 
These categories are equivalent to the corresponding categories of sheaves on the prestack $X_\bbF/\phi_X$, that is, equivalent to the homotopy fixed points of the induced $\phi_X^*$-action:

\prop[\thref{presentation.Weil.prop}, \thref{prop.cons.weil.as.fixed.points}]
\thlabel{Weil.sheaves.limit}
The pullback of sheaves along $(X_\bbF)_\proet\r X^\Weil_\proet$ induces an equivalence of $\Lambda_*$-linear symmetric monoidal stable \ii-categories
\[
\D_\bullet\big(X^\Weil,\Lambda\big) \stackrel\cong\lr \D_\bullet(X_{\bbF},\Lambda)^{\phi_X^*=\id},
\]
for $\bullet\in \{\varnothing,\lis,\cons\}$.
\xprop

Thus, objects in $\D_\bullet\big(X^\Weil,\Lambda\big)$ are pairs $(M, \al)$ with $M\in \D_\bullet(X_\bbF,\Lambda)$ and $\al\co M \cong \phi_X^*M$.
On the abelian level, we recover the classical approach \cite[Definition 1.1.10]{Deligne:Weil2}.
If $\Lambda$ is a finite discrete ring, then every Weil descent datum on constructible $\Lambda$-sheaves is effective so that $\D_\cons\big(X^\Weil,\Lambda\big)\cong \D_\cons(X,\Lambda)$, see \thref{discrete_Geisser}.
However, the categories are not equivalent if $\Lambda=\bbZ,\bbZ_\ell, \bbQ_\ell$, say.  
This relates to the difference between continuous representations of Galois groups such as $\hat \bbZ$ versus Weil groups such as $\bbZ$.

For several $\bbF_q$-schemes $X_1,\ldots, X_n$, a similar process is carried out for their product $X:=X_1\x_{\bbF_q}\ldots \x_{\bbF_q} X_n$ equipped with the partial Frobenii $\phi_{X_i}\co X_\bbF\r X_\bbF$, see \refsect{Weil.sheaves.product}. 
Generalizing \thref{Weil.sheaves.limit}, there is an equivalence of $\Lambda_*$-linear symmetric monoidal stable \ii-categories
\eqn
\label{fixed.point.intro.eq}
\D_\bullet\big(X_1^\Weil\x\ldots\x X_n^\Weil,\Lambda\big) \stackrel\cong\lr \D_\bullet(X_{\bbF},\Lambda)^{\phi_{X_1}^*=\id,\ldots,\phi_{X_n}^*=\id}
\xeqn
for $\bullet\in \{\varnothing,\lis,\cons\}$. 
The category on the left is defined using the Weil-pro\'etale site $\big(X_1^\Weil\x\ldots\x X_n^\Weil\big)_\proet$ consisting of objects $(U,\varphi_1,\ldots,\varphi_n)$ with $U\in (X_\bbF)_\proet$ and pairwise commuting endomorphisms $\varphi_i\co U\r U$ covering the partial Frobenii $\phi_{X_i}\co X_\bbF\r X_\bbF$ for all $i=1,\ldots, n$.
The category on the right is the category of simultaneous homotopy fixed points, see \refsect{fixed.points.infty}.
For constructible Weil sheaves, \eqref{fixed.point.intro.eq} relies on decompositions of partial Frobenius invariant cycles in $X_\bbF$, see \thref{partial.Frobenius2}.

The following result is referred to as the categorical K\"unneth formula for constructible Weil sheaves (or, derived Drinfeld's lemma): 

\theo
[\thref{derived_Drinfeld.text}, \thref{limits.Drinfeld.rema}]
\thlabel{derived_Drinfeld}
Let $\bbF_q$ be a finite field of characteristic $p>0$.
Let $X_1, \dots, X_n$ be finite type $\Fq$-schemes.
Let $\Lambda$ be either a finite discrete ring of prime-to-$p$-torsion, or an algebraic field extension $E \supset \Ql$, $\ell \ne p$, or its ring of integers $\calO_E$.

Then the external tensor product of sheaves $(M_1,\ldots,M_n)\mapsto M_1\bx\ldots\bx M_n$ induces an equivalence
\eqn
\label{derived_Drinfeld_map.intro}
\D_\cons\big(X_1^\Weil,\Lambda\big)\t_{\Perf_{\Lambda_*}} \ldots\t_{\Perf_{\Lambda_*}}  \D_\cons\big(X_n^\Weil, \Lambda\big) \stackrel\cong\lr \D_\cons\big(X_1^\Weil\x\ldots\x X_n^\Weil,\Lambda\big),
\xeqn
and likewise for the categories of lisse Weil sheaves if, in the case $\Lambda=E$, one assumes the schemes $X_1,\ldots, X_n$ to be geometrically unibranch (for example, normal).
\xtheo

This statement can also be recast as the symmetric monoidality of the functor sending a Weil prestack $X^\Weil$, which is defined on $R$-points by $X^\Weil(R):= \colim (X(R) \stackrel[\id]{\phi_X}\rightrightarrows X(R))$, to its \ii-category of constructible sheaves (\thref{derived_Drinfeld.prestacks}).

The tensor product of \ii-categories (see \refsect{recollections}) is formed using the natural $\Lambda_*$-linear structures on the categories. 
We have an analogous equivalence for the categories of lisse Weil sheaves with coefficients $\Lambda$ in finite discrete $p$-torsion rings like $\bbZ/p^m$, $m\geq 1$, see \thref{derived_Drinfeld.text}.
As the following example shows, the use of Weil sheaves is necessary for the essential surjectivity to hold.
This behavior is mentioned in the first arXiv version of \cite[Equation (0.8)]{GKRV.Toy.Model} which is one of the main motivations for our work.

\exam[{compare \cite[Exposé X, \S1, Remarques 1.10]{SGA1}}]
\thlabel{AS.exam.intro}
Let $X_{1,\bbF}=X_{2,\bbF}=\bbA^1_{\bbF}$ be the affine line so that $X_\bbF=\bbA^2_\bbF$ with coordinates denoted by $x_1$ and $x_2$. 
Then 
\[
U:=\{t^p-t = x_1\cdot x_2\} \lr \bbA^2_{\bbF}
\]
defines a finite \'etale cover with Galois group $\bbZ/p$.
Let $M\in \D_\lis(\bbA^2_\bbF,\Lambda)$ be the sheaf in degree $0$ associated with some non-trivial character $\bbZ/p\r \Lambda^\x_*$.
For $\lambda,\mu\in\bbF$ not differing by a scalar in $\bbF_p^\x$, the fibers $U|_{\{x_1=\lambda\}}$, $U|_{\{x_1=\mu\}}$ are not isomorphic over $\bbA_\bbF^1$ by Artin-Schreier theory.
Hence, $M\not\simeq\phi_{X_i}^*M$ and one can show that $M\not \simeq M_1\bx M_2$ for any $M_i\in \D(\bbA^1_\bbF,\Lambda)$.
\xexam

If $\Lambda$ as above is $p$-torsion free, then the full faithfulness of \eqref{derived_Drinfeld_map.intro} is a direct consequence of the Künneth formula for $X_{i,\bbF}$, $i=1,\ldots,n$. 
For $\Lambda=\bbZ/p^m$, we use Artin--Schreier theory instead.
It would be interesting to see whether the lisse $p$-torsion case can be extended to constructible sheaves. 
In both cases, the essential surjectivity relies on a variant of Drinfeld's lemma for Weil group representations, see \thref{Drinfelds_lemma},
together with a characterization of partial-Frobenius stable algebraic cycles (\thref{partial.Frobenius2}) as well as a decomposition argument for representations of a product of abstract groups (\thref{factorization_lemm}).

With a view towards \cite{Lafforgue:Chtoucas}, we consider Weil sheaves whose underlying sheaf is ind-constructible, but where the action of the partial Frobenii do not necessarily preserve the constructible pieces.
For finite type $\bbF_q$-schemes $X_1,\ldots, X_n$ and $\Lambda$ as in \thref{derived_Drinfeld}, we consider the category of simultaneous homotopy fixed points
\[
\D_{\bullet}\big(X_1^\Weil\x\ldots\x X_n^\Weil,\Lambda\big) \defined \D_\bullet(X_{\bbF},\Lambda)^{\phi_{X_1}^*=\id, \ldots, \phi_{X_n}^*=\id}
\]
 for $\bullet\in \{\indlis,\indcons\}$.
 Then the external tensor product induces a fully faithful functor 
\eqn
\label{ind.cons.product.eq}
\D_{\bullet}\big(X_1^\Weil,\Lambda\big)\t_{\Mod_{\Lambda_*}}\ldots\t_{\Mod_{\Lambda_*}} \D_{\bullet} \big(X_n^\Weil,\Lambda\big) \lr \D_{\bullet} \big(X_1^\Weil\x \ldots\x X_n^\Weil,\Lambda\big).
\xeqn
Unlike the case of lisse or constructible sheaves, the functor is not essentially surjective as one can add freely actions by the partial Frobenii, see \thref{ind.cons.equivalence.example}.
However, we can identify a large class of objects in the essential image of \eqref{ind.cons.product.eq}. 
When combined with the smoothness results of Xue \cite[Theorem 4.2.3]{Xue:Smoothness}, we obtain, for example, that the compactly supported cohomology of moduli stacks of shtukas over global function fields lies in the essential image of \eqref{ind.cons.product.eq}, see \refsect{shtuka.cohomology} for details. 

\rema
Another motivation for this work is our (T.R.~and J.S.) ongoing project aiming for a motivic refinement of \cite{Lafforgue:Chtoucas}. 
In this project, we will need a motivic variant of Drinfeld's lemma. 
Since triangulated categories of motives such as $\DM(X,\bbQ)$ carry t-structures only conditionally, we need a Drinfeld lemma to be a statement about triangulated categories. 
In conjunction with the conjecture relating Weil-étale motivic cohomology to Weil-étale cohomology \cite{Kahn:Equivalences,Geisser:Weil,Lichtenbaum:Weil}, our results suggest to look for a Drinfeld lemma for constructible Weil motives. 
\xrema

\subsection*{Acknowledgements} 
We thank Clark Barwick, Jean-Fran\c{c}ois Dat, Christopher Deninger, Rune Haugseng, Claudius Heyer, Peter Schneider, Burt Totaro, Torsten Wedhorn, Alexander Yom Din, and Xinwen Zhu for helpful conversations and email exchanges.

\section{Recollections on \ii-categories}
\label{sect--recollections}

Throughout this section, $\Lambda$ denotes a unital, commutative ring. 
We briefly collect some notation pertaining to \ii-categories from \cite{Lurie:HA,Lurie:Higher}.
As in \cite[Section 5.5.3]{Lurie:Higher}, $\PrL$ denotes the \ii-category of presentable \ii-categories with colimit-preserving functors. 
It contains the subcategory $\PrSt \subset \PrL$ consisting of stable \ii-categories.

\subsection{Monoidal aspects}
\label{sect--monoidal.aspects}

The category $\PrL$ carries the Lurie tensor product \cite[Section 4.8.1]{Lurie:HA}.
This tensor product induces one on the full subcategory $\PrSt \subset \PrL$ consisting of stable \ii-categories \cite[Proposition~4.8.2.18]{Lurie:HA}.
For our commutative ring $\Lambda$, the \ii-category $\Mod_\Lambda$ of chain complexes of $\Lambda$-modules, up to quasi-isomorphism, is a commutative monoid in $\PrSt$ with respect to this tensor product.
This structure includes, in particular, the existence of a functor
$$\Mod_\Lambda \x \Mod_\Lambda \r \Mod_\Lambda$$
which, after passing to the homotopy categories is the classical \emph{derived} tensor product on the unbounded derived category of $\Lambda$-modules.

We define $\PrStL$ to be the category of modules, in $\PrSt$, over $\Mod_\Lambda$.
Noting that modules over $\Mod_\Lambda$ are in particular modules over $\Sp$, the \ii-category of spectra, $\PrStL$ can be described as the \ii-category consisting of \emph{stable} presentable \ii-categories together with a $\Lambda$-linear structure, such that functors are continuous and $\Lambda$-linear.
Therefore $\PrStL$ carries a symmetric monoidal structure, whose unit is $\Mod_\Lambda$. We will also denote by $\PrSt_{\omega}$ the category of compactly generated presentable with functors that send compact objects to compact objects (equivaletly, those whose right adjoint is continuous).

In order to express monoidal properties of \ii-categories consisting, say, of bounded complexes, 
recall from \cite[Corollary~4.8.1.4 joint with Lemma~5.3.2.11]{Lurie:HA} or \cite[Proposition~4.4]{BenZviFrancisNadler:Integral} the symmetric monoidal structure on the \ii-category $\CatEx_\infty(\Idem)$ of idempotent complete stable \ii-categories and exact functors: it is characterized by
$$D_1 \t D_2 \defined \big(\Ind(D_1) \t \Ind(D_2)\big)^\omega,\eqlabel{Cat.idem.monoidal}$$
that is, the compact objects in the Lurie tensor product of the Ind-completions. 
With respect to these monoidal structures, the Ind-completion functor (taking values in compactly generated presentable \ii-categories with the Lurie tensor product) and the functor forgetting the compact generatedness
\eqn\CatEx_\infty(\Idem) \stackrel [\cong]{\Ind} \lr \PrSt_{\omega} \lr \PrSt
\label{eqn--CatExIdem.etc}
\xeqn
are both symmetric monoidal \cite[Lemmas~5.3.2.9, 5.3.2.11]{Lurie:HA}.

The subcategory of compact objects in $\Mod_\Lambda$ is given by perfect complexes of $\Lambda$-modules \cite[Proposition~7.2.4.2.]{Lurie:HA}.
It is denoted $\Perf_\Lambda$.
Under the equivalence in \refeq{CatExIdem.etc}, the category $\Perf_\Lambda \in \CatEx_\infty(\Idem)$ corresponds to $\Mod_\Lambda$. Moreover, $\Perf_\Lambda$ is a commutative monoid in $\CatEx_\infty(\Idem)$, so that we can consider its category of modules, denoted as
$\CatEx_{\infty, \Lambda}(\Idem)$.
This category inherits a symmetric monoidal structure denoted by $D_1 \t_{\Perf_\Lambda} D_2$.

Any stable \ii-category $D$ is canonically enriched over the category of spectra $\Sp$.
We write $\Hom_D(-, -)$ for the mapping spectrum.
Any category in $\PrStL$ is canonically enriched over $\Mod_\Lambda$, so that we refer to $\Hom_D(-, -) \in \Mod_\Lambda$ as the mapping complex.
For example, for $M, N\in \Mod_\Lambda$, then $\Hom_{\Mod_\Lambda}(M, N)$ is commonly also denoted by $\RHom(M, N)$.
Its $n$-th cohomology is the Hom-group $\Hom(M, N[n])$ in the classical derived category.

\subsection{Fixed points of \ii-categories}
\label{sect--fixed.points.infty}
A basic structure in Drinfeld's lemma is the equivariance datum for the partial Frobenii.
In this section, we assemble some abstract results where such \ii-categorical constructions are carried out.

\defi
\thlabel{fixed.point.definition}
Let $\phi\co D \r D$ be an endofunctor in $\CatEx_{\infty}(\Idem)$.
The category of \emph{$\phi$-fixed points} is
$$D^{\phi = \id} \defined \Fix(D, \phi) \defined \lim \left (D \stackrel[\id_D]{\phi}\rightrightarrows D \right).$$
\xdefi


Recall that for a symmetric monoidal \ii-category $D$, a commutative monoid object $\Lambda \in \CAlg(D)$, the forgetful functors $\CAlg(D) \r D$ and $\Mod_\Lambda(D) \r D$ preserve limits \cite[Corollary~3.2.2.5, Corollary~4.2.3.3]{Lurie:HA}.
In particular, if $D$ is in addition $\Lambda$-linear, that is, an object in $\CatEx_{\infty, \Lambda}(\Idem)$, and $\phi$ is also $\Lambda$-linear, then $\Fix(D, \phi)$ admits a natural $\Lambda$-linear structure as well.

Because of these facts, we will usually not specify where the limit above is formed.
Note that all functors 
\eqn\CatEx_\infty(\Idem) \stackrel [\cong]{\Ind} \lr \PrSt_{\omega} \stackrel {\text{(}*\text{)}} \lr \PrSt \lr \PrL \lr \widehat{\Cat_\infty}
\label{eqn--CatExIdem.etc2}
\xeqn
except for the forgetful functor marked ($*$) preserve limits, see \cite[Corollary~4.2.3.3]{Lurie:HA} and \cite[Proposition~5.5.3.13]{Lurie:Higher} for the rightmost two functors.
To give a concrete example of that failure in our situation, note that $\Fix(D, \id_D) = \Fun(B\Z, D)$, that is, objects are pairs $(M, \alpha)$ consisting of some $M \in D$ and some automorphism $\alpha\co M \cong M$.
Now consider $D = \Vect^\fd_\Lambda$, the (abelian) category of finite-dimensional vector spaces over a field $\Lambda$. 
The natural functor
$$\Ind \left (\lim \big(\Vect^\fd_\Lambda \rightrightarrows \Vect^\fd_\Lambda\big) \right) \r 
\lim \left (\Ind \big(\Vect^\fd_\Lambda\big) \rightrightarrows \Ind\big( \Vect^\fd_\Lambda\big) \right ) = 
\lim \big(\Vect_\Lambda \rightrightarrows \Vect_\Lambda\big)$$
is fully faithful, but \emph{not} essentially surjective: given an automorphism $\alpha$ of an infinite-dimensional vector space $M$, there need not be a filtration $M = \bigcup M_i$ by finite-dimensional subspaces $M_i$ that is compatible with $\alpha$.

Fixed point categories inherit t-structures as follows:

\lemm
\thlabel{Fix.t-structure}
Let $\phi\co D \r D$ be a functor in $\CatEx_\infty(\Idem)$.
Suppose $D$ carries a t-structure such that $\phi$ is t-exact.
Then $\Fix(D, \phi)$ carries a unique t-structure such that the evaluation functor is t-exact. 
There is a natural equivalence
$$\Fix (D^\heartsuit, \phi)\stackrel\cong\lr \Fix (D, \phi)^\heartsuit.$$
\xlemm

\pf
Let us abbreviate $\widetilde D := \Fix(D, \phi)$.
For $\bullet$ being either ``$\le 0$'' or ``$\ge 0$'', we put $\widetilde D^\bullet := \Fix(D^\bullet, \phi)$, which is a (non-stable) \ii-category. 
This is clearly the only choice for a t-structure making $\ev$ a t-exact functor.
It satisfies the claim about the hearts of the t-structure by definition.

We need to show that it is a t-structure. 
Being a limit of full subcategories, the categories $\widetilde D^\bullet$ are full subcategories of $\widetilde D$.
Since $\phi$, being t-exact, commutes with $\tau_D^{\le 0}$ and $\tau_D^{\ge 0}$, these two functors also yield truncation functors for $\widetilde D$.
For $M \in \widetilde D^{\le 0}$, $N \in \widetilde D^{\ge 1}$ (we use cohomological conventions), we have
$$\Hom_{\widetilde D}\big(M, N\big) \;=\; \lim \big (\Hom_D(M, N) \rightrightarrows \Hom_D(M, N) \big ),$$
where on the right hand side $M$, $N$ denote the underlying objects in $D$.
Since $M \in D^{\le 0}$, $N \in D^{\ge 1}$, we have $\H^i\Hom_D(M, N)=0$ for $i=-1,0$.
Thus, $\H^0\Hom_{\widetilde D}(M, N)=0$ as well.
\xpf

\thref{fixed.point.definition} can be generalized as follows: 
Let $\varphi\co B\Z^n \r \CatEx_\infty(\Idem)$ be a diagram. 
For example, for $n=1$, this amounts to giving $D=\varphi(*) \in \CatEx_\infty(\Idem)$ and an equivalence $\phi=\varphi(1)\co D \r D$. 
For $n = 2$, such a datum corresponds to giving $D$, equivalences $\phi_1, \phi_2\co D \stackrel \cong \r D$ together with an equivalence $\phi_1 \circ \phi_2 \stackrel \cong \r \phi_2 \circ \phi_1$.
So we define the \ii-category of \emph{simultaneous fixed points} as
$$\Fix(D, \phi_1, \dots, \phi_n) \defined\lim \varphi \in \CatEx_\infty(\Idem).$$

\rema
\thlabel{simultaneous.Fix}
The statement of \thref{Fix.t-structure} carries over verbatim assuming that $D$ has a t-structure and all $\phi_i$ are t-exact, noting that $B\Z^n = (S^1)^n$ is a finite simplicial set.
\xrema

\lemm
\thlabel{ind.to.fixed.fully.faithful.big.category}
Let $\varphi\co B\Z^n \r \CatEx_\infty(\Idem)$ be a diagram. 
Denote $D = \varphi(*)$ and $\phi_i=\varphi(e_i)$ for the $i$-th standard vector $e_i\in \bbZ^n$.  
The functor 
$$\Fix\big(D, \phi_1,\ldots,\phi_n\big) \r \Fix \big(\Ind (D), \phi_1,\ldots,\phi_n\big)$$
induced from the inclusion $D \subset \Ind (D)$ is fully faithful and takes values in compact objects. 
In particular, it yields a fully faithful functor
$$\Ind\big(\Fix (D, \phi_1, \ldots, \phi_n)\big) \r \Fix \big(\Ind (D), \phi_1, \ldots, \phi_n\big).$$
\xlemm

\pf
Let $M \in \Fix(D, \phi_1, \dots, \phi_n)$ and denote its underlying object in $D$ by the same symbol.
For every $N \in \Fix\big(\Ind (D), \phi_1, \dots, \phi_n)$, we have a limit diagram of mapping complexes
\[
\Hom_{\Fix(\Ind (D))}(M, N) \cong \Fix \big(\Hom_{\Ind (D)}(M,N), \phi_1, \dots, \phi_n \big).
\]
Since filtered colimits commute with finite limits in the \ii-category of anima (a.k.a.~spaces) \cite[Proposition 5.3.3.3.]{Lurie:Higher}, we see that $M$ is compact in $\Fix(\Ind(D))$ because $M$ is so in $\Ind (D)$.
\xpf

\lemm
\thlabel{fix.tensor.product.big.categories}
Let $\varphi_i\co B\Z \r \CatEx_{\infty, \Lambda}(\Idem)$, $i = 1, \dots, n$ be given. 
Denote $D_i=\varphi_i(*)$, $\phi_i=\varphi_i(1)$ and $\widetilde D_i=\Ind(D_i)$.
Then there is a canonical equivalence
$$
\Fix\big(\widetilde D_1, \phi_1\big)\t_{\Mod_\Lambda}\ldots\t_{\Mod_\Lambda} \Fix\big(\widetilde D_n, \phi_n\big) \stackrel \cong \r \Fix \left (\widetilde D_1\t_{\Mod_\Lambda}\ldots\t_{\Mod_\Lambda}\widetilde D_n, \phi_1, \dots, \phi_n \right ).$$
\xlemm

\pf
The categories $\Fix(\widetilde D_i, \phi_i)$ are compactly generated: the forgetful functor $U : \Fix(\widetilde D_i, \phi_i) \r \widetilde D_i = \Ind(D_i)$ preserves colimits, so its left adjoint $L$ preserves compact objects. Moreover, $U$ is conservative, so that the objects $L(d_i)$, for $d_i \in D_i$, form a family of compact generators. 
Then, we use that any compactly generated category in $\PrStL$ is dualizable \cite[Remark D.7.7.6 (1)]{Lurie:SAG} so that tensoring with it preserves limits.
\xpf

\section{Lisse and constructible sheaves}

In order to state and prove the categorical Künneth formula for Weil sheaves, we use the framework for lisse and constructible sheaves provided by \cite{HemoRicharzScholbach:Constructible}.
For the convenience of the reader, we collect here some basics of the formalism.

Throughout, $\Lambda$ denotes a condensed ring, for example any T1-topological ring such as discrete rings, algebraic extensions $E / \Ql$ or their ring of integers $\calO_E$.
In the synopsis below, we refer to the latter choices of $\Lambda$ as the \emph{standard coefficient rings}.
We write $\Lambda_*$ for the underlying ring.
Let $\D(X, \Lambda)$ be the derived category of sheaves of $\Lambda$-modules on the proétale site $X_\proet$. 

\defi[{\cite[\myconsref{lisse.constructible.defi}, \myconsref{ind.lis.cons.defi}]{HemoRicharzScholbach:Constructible}}]
\thlabel{recall.definitions}
For every scheme $X$ and every condensed ring $\Lambda$, there are the full subcategories

$$\Dlis(X, \Lambda) \subset \Dcons(X, \Lambda) \subset \D(X, \Lambda).\eqlabel{categories.roster}$$
By definition, the left hand category of \emph{lisse sheaves} consists of the dualizable objects in the right-most category.
An object (henceforth referred to as a sheaf) $M$ in the right hand category is \emph{constructible}, if on any affine $U \subset X$ there is a finite stratification into constructible locally closed subschemes $U_i \subset U$ such that $M|_{U_i}$ is lisse, that is, dualizable.
Finally, an \emph{ind-lisse} (respectively,~\emph{ind-constructible}) sheaf is a filtered colimit, in the category $\D(X, \Lambda)$, of lisse (respectively,~constructible) sheaves.
The corresponding full subcategories of $\D(X, \Lambda)$ are denoted by
$$\Dindlis(X, \Lambda) \subset \Dindcons(X, \Lambda) \subset \D(X, \Lambda).$$
\xdefi

For the standard coefficient rings $\Lambda$ above and quasi-compact quasi-separated (qcqs) schemes $X$, that definition of lisse and constructible sheaves agrees with the classical ones, see \consrefsect{comparison.results} for details.

The categories enjoy the following properties:
\syno
\thlabel{sheaves}
\begin{enumerate}[(i)]
  \item
  \label{item--Mod.Lambda}
  Via the natural functor $\Mod_{\Lambda_*}\r \D(X,\Lambda)$, $M\mapsto \underline{M}\otimes_{\underline{\Lambda_*}}\Lambda_X$ (see around \conseqref{constant.objects.functor}, the category $\D(X,\Lambda)$ is an object in $\PrSt_{\Lambda_*}$. 
  The functor restricts to a functor $\Perf_{\Lambda_*}\r \D_\lis(X,\Lambda)$, and the categories $\Dlis(X,\Lambda)\subset \Dcons(X,\Lambda)$ are objects in $\CatEx_{\infty, \Lambda}(\Idem)$.
  In particular, all categories listed in \refeq{categories.roster} are stable idempotent complete $\Lambda_*$-linear \ii-categories.

  \item
  \label{item--preservation.constructibility}
  The extension-by-zero functor along any constructible locally closed immersion and quasi-compact étale morphisms preserves constructibility, see \cite[\myconsref{t.exactness.lemm}, \myconsref{preservation.constructibility}]{HemoRicharzScholbach:Constructible}.

  \item
  \label{item--descent}
  The functors $X \mapsto \Dcons(X, \Lambda)$ and $X\mapsto \Dlis(X, \Lambda)$ satisfy proétale hyperdescent \consref{lisse.constructible.hyperdescent.coro}.
  (According to \cite[Theorem~2.2]{HansenScholze:RelativePerversity}, it also satisfies v-descent, but we will not need this in this paper.)
  The functor $X \mapsto \Dindcons(X, \Lambda)$, resp.~$X\mapsto \Dindlis(X, \Lambda)$ satisfies hyperdescent for quasi-compact étale, resp.~finite étale covers, see \consref{descent.ind.cons}.

  \item
  \label{item--manipulations.Lambda}
  If $\Lambda = \colim \Lambda_i$ is a filtered colimit of condensed rings and $X$ is qcqs, then the natural functors
\[
\colim \D_\lis(X,\Lambda_i)\stackrel \cong \lr \D_\lis(X,\Lambda), \;\; \colim \D_\cons(X,\Lambda_i)\stackrel \cong \lr \D_\cons(X,\Lambda)
\]
  are equivalences \consref{colimit.coefficients.lemm}.

  \item
  \label{item--constructible.bounded}
  If $X$ is qcqs, then any constructible sheaf is bounded with respect to the t-structure on $\D(X, \Lambda)$ \consref{constructible.bounded.coro}.

  \item
  \label{item--t-structure}
  For $X$ locally Noetherian (and much more generally), the t-structure on $\D(X, \Lambda)$ restricts to one on $\Dlis(X, \Lambda)$ and $\Dcons(X, \Lambda)$ provided that $\Lambda$ is t-admissible in the sense of \consref{t.admissible.defi}.
  Here, t-admissibility is a combination of an algebraic and a topological condition: first, $\Lambda_*$ needs to be regular coherent (for example, any regular Noetherian ring of finite Krull dimension, but $\Z/\ell^2$ is excluded). 
  The topological condition on the condensed structure of $\Lambda$ is satisfied for all the standard coefficient rings listed above, see \consref{t.structure.condensed}.

  \item
  \label{item--lisse.locally.constant}
  For $X$ locally Noetherian (and again more generally), a sheaf is lisse if and only if it is proétale locally the constant sheaf associated to a perfect complex of $\Lambda_*$-modules, see \consref{locally.constant.prop}.

  \item
  \label{item--compact.objects}
  Let $X$ be a qcqs scheme.
  If the $\Lambda$-cohomological dimension is uniformly bounded for all proétale affines $U=\lim_i U_i$ over $X$, then $\Ind(\Dcons(X, \Lambda)) = \Dindcons(X, \Lambda)$ and likewise for ind-lisse sheaves.
  If $X$ is of finite type over $\Fq$ or a separably closed field, this condition holds for any of the above standard rings.
  For discrete $p$-torsion rings, algebraic extensions $E/\Qp$ and their ring of integers $\calO_E$, this holds for arbitrary qcqs schemes in characteristic $p$.
  See \cite[\myconsref{Artin.vanishing}, \myconsref{compact.objects}]{HemoRicharzScholbach:Constructible}.
\end{enumerate}
\xsyno

For schemes $X_1,\ldots, X_n$ over a fixed base scheme $S$ (for example, the spectrum of a field) and a condensed ring $\Lambda$, we denote the external product in the usual way:
\begin{align*}
\bx\co \D(X_{1},\Lambda)\x\ldots \x \D(X_{n},\Lambda) &\lr \D\left(X_1\x_S\ldots \x_S X_n,\Lambda\right), \\
(M_1, \ldots, M_n) &\longmapsto M_1 \bx\ldots\bx M_n := p_1^* (M_1)\t_{\Lambda_X}\ldots \t_{\Lambda_X} p_n^* (M_n).
\end{align*}
Here $p_i\co X := X_1\x_S\ldots\x_S X_n\to X_i$ are the projections.
This functor induces the functor
\eqn
\label{tensor.product}
\bx\co \D(X_{1},\Lambda)\t_{\Mod_{\Lambda_*}}\ldots\t_{\Mod_{\Lambda_*}} \D(X_{n},\Lambda) \r \D\left(X_1\x_S\ldots \x_S X_n,\Lambda\right),
\xeqn
in $\PrSt_{\Lambda_*}$.
Here we regard $\D(X_{i},\Lambda)$ as objects in $\PrSt_{\Lambda_*}$, like in \refit{Mod.Lambda} in the synopsis above.
The external tensor product of constructible sheaves is again constructible, and hence induces a functor 
\eqn
\label{tensor.product.constructible}
\bx\co \D_\cons(X_1,\Lambda)\t_{\Perf_{\Lambda_*}}\ldots\t_{\Perf_{\Lambda_*}} \D_\cons(X_{n},\Lambda) \r \D_\cons\left(X_1\x_S\ldots \x_S X_n,\Lambda\right),
\xeqn
in $\CatEx_{\infty,\Lambda_*}(\Idem)$ and likewise for the categories of ind-constructible, resp.~(ind-)lisse sheaves.

\section{Weil sheaves}\label{sect--section.weil.sheaves}

In this section, we introduce the categories 
$$\D_\lis\big(X^\Weil, \Lambda\big) \;\subset\; \D_\cons\big(X^\Weil, \Lambda\big) \;\subset\; \D\big(X^\Weil, \Lambda\big)$$
consisting of lisse, resp.~constructible, resp.~all Weil sheaves. 
These are the categories featuring in the categorical Künneth formula (\thref{derived_Drinfeld}).

Throughout this section, $X$ is a scheme over a finite field $\Fq$ of characteristic $p>0$.
Unless the contrary is mentioned, we impose no conditions on $X$.
Moreover, $\Lambda$ is a condensed ring. 
We fix an algebraic closure $\bbF$ of $\Fq$, and denote by $X_\bbF:=X\x_{\Fq}\Spec \bbF$ the base change.
Denote by $\phi_X$ (resp.~$\phi_\bbF$) the endomorphism of $X_\bbF$ that is the $q$-Frobenius on $X$ (resp.~$\Spec\bbF$) and the identity on the other factor. 

Let 
$$\D_\lis(X_\bbF,\Lambda)\subset \D_\cons(X_\bbF,\Lambda) \subset \D(X_\bbF, \Lambda)$$
be the categories of lisse, resp.~constructible, resp.~all proétale sheaves of $\Lambda$-modules on $X_\bbF$ (\thref{recall.definitions}).
These categories are objects in $\CatEx_{\infty,\Lambda_*}(\Idem)$, that is, $\Lambda_*$-linear stable idempotent complete symmetric monoidal \ii-categories where $\Lambda_*=\Gamma(*,\Lambda)$ is the underlying ring.

\subsection{The Weil-pro\'{e}tale site}
\label{sect--Def.Weil} 
The Weil-\'etale topology for schemes over finite field is introduced in \cite{Lichtenbaum:Weil}, see also \cite{Geisser:Weil}. 
Our approach for the pro\'etale topology is slightly different: 

\defi 
The \textit{Weil-pro\'{e}tale site} of $X$, denoted by $X^\Weil_\proet$, is the following site: 
Objects in $X^\Weil_\proet$ are pairs $(U,\varphi)$ consisting of $U\in (X_{\bbF})_{\proet}$ equipped with an endomorphism $\varphi\colon U\rightarrow U$ of $\bbF$-schemes such that the map $U\rightarrow X_{\bbF}$ intertwines $\varphi$ and $\phi_X$. Morphisms in $X^\Weil_\proet$ are given by equivariant maps, and a family $\{(U_i,\varphi_i) \rightarrow (U,\varphi)\}$ of morphisms is a cover if the family $\{U_i \rightarrow U\}$ is a cover in $(X_{\bbF})_{\proet}$.
\xdefi

Note that $X^\Weil_\proet$ admits small limits formed componentwise as $\lim (U_i,\varphi_i)=(\lim U_i,\lim\varphi_i)$.
In particular, there are limit-preserving maps of sites
\begin{equation}
\label{map.Weil.sites.eq}
(X_\bbF)_\proet\r X^\Weil_\proet \r X_\proet
\end{equation}
given by the functors (in the opposite direction) $U \mapsfrom (U, \varphi)$ and $(U_\bbF,\phi_U)\mapsfrom U$.
We denote by $\D(X^\Weil,\Lambda)$ the unbounded derived category of sheaves of $\Lambda_X$-modules on $X^\Weil_\proet$. 
The maps of sites \eqref{map.Weil.sites.eq} induce functors 
\eqn
\label{map.Weil.category.eq}
\D(X,\Lambda)\r \D(X^\Weil,\Lambda)  \r \D(X_{\bbF},\Lambda),
\xeqn
whose composition is the usual pullback functor along $X_\bbF\r X$.

\rema
The functor $\D(X,\Lambda) \rightarrow \D(X^\Weil,\Lambda)$ is not an equivalence in general. 
This relates to the difference between continuous representations Galois versus Weil groups.
See, however, \thref{discrete_Geisser} for filtered colimits of finite discrete rings $\Lambda$.
\xrema

We have the following basic functoriality: 
Let $j \colon U \rightarrow X$ be a weakly \'{e}tale morphism and consider the corresponding object $(U_{\bbF},\phi_U)$ of $X^\Weil_{\proet}$. Then the slice site $(X^\Weil_{\proet})_{/(U_{\bbF},\phi_U)}$ is equivalent to $U^\Weil_{\proet}$. This gives a functor $(X_\proet)^\opp \rightarrow  \PrStL$, $U\mapsto \D(U^\Weil,\Lambda)$ which is a hypercomplete sheaf of $\Lambda_*$-linear presentable stable categories.

Also, we obtain an adjunction
\[
j_! \colon \D(U^\Weil,\Lambda) \rightleftarrows \D(X^\Weil,\Lambda) {\colon j^*} 
\]
that is compatible with the $((j_\bbF)_!, (j_\bbF)^*)$-adunction under \eqref{map.Weil.category.eq}. The category $\D(X^\Weil,\Lambda)$ is equivalent to the category of $\phi_X$-equivariant sheaves on $X_{\bbF}$, as we will now explain. 

For each $i\geq 0$, consider the object $(X_i,\Phi_i)\in X^\Weil_{\proet}$ with $X_i=\bbZ^{i+1}\times X_{\bbF}$ the countably disjoint union of $X_\bbF$, the map $X_i\r X_\bbF$ given by projection and the endomorphism $\Phi_i\colon X_i\r X_i$ given by $(\underline n,x)\mapsto (\underline n-(1,\ldots,1),\phi_X(x))$ on sections.
The inclusion $\bbZ^i\r \bbZ^{i+1}$, $\underline n\mapsto (0,\underline n)$ induces a map of schemes $X_{i-1}\r X_{i}$ where $X_{-1}:=X_\bbF$.
By pullback, we get a limit-preserving map of sites
\begin{equation}
\label{comparison.Weil.site.eq}
(X_{i-1})_\proet\r \left(X_\proet^\Weil\right)_{/(X_i,\Phi_i)}
\end{equation}

\lemm
\thlabel{lemma.trivialization.equivariant.Weil.site}
For each $i\geq 0$, the map \eqref{comparison.Weil.site.eq} induces an equivalence on the associated $1$-topoi. 
\xlemm
\pf
As universal homeomorphisms induce equivalences on pro\'etale $1$-topoi \cite[Lemma~5.4.2]{BhattScholze:ProEtale}, we may assume that $X$ is perfect. 
In this case, the sites \eqref{comparison.Weil.site.eq} are equivalent because $\phi_X$ is an isomorphism. 
Explicitly, an inverse is given by sending an object $U\in (X_{i-1})_\proet$ to the object $V=\bigsqcup_{\underline n\in \bbZ^{i+1}} V_{\underline n}$, $V_{\underline n}\r \{\underline n\}\x X_\bbF$ defined by 
\[
V_{\underline n}=U_{(n_2-n_1,\ldots,n_{i+1}-n_1)}\x_{X_\bbF, \phi_X^{n_1}} X_\bbF, 
\]
and with endomorphism $\varphi \colon V \rightarrow V$ defined by the maps $V_{\underline n}=V_{\underline n -(1,\ldots,1)}\x_{X_\bbF, \phi_X} X_\bbF\r V_{\underline n -(1,\ldots,1)}$.
\xpf

Weil sheaves admit the following presentation as the $\phi_X^*$-fixed points of $\D(X_{\bbF},\Lambda)$, see \thref{fixed.point.definition}:

\prop
\thlabel{presentation.Weil.prop}
The last functor in \eqref{map.Weil.category.eq} induces an equivalence
\eqn
\label{eqn--Weil.sheaves}
\D(X^\Weil,\Lambda)\ \cong \lim\left (\D(X_{\bbF},\Lambda) \stackrel[\id]{\phi_X^*}\rightrightarrows \D(X_{\bbF},\Lambda) \right).
\xeqn
\xprop

\rema 
\thlabel{remark.weil.phi-X.phi-F}
Objects in \refeq{Weil.sheaves} are pairs $(M,\al)$ where $M \in \D(X_\bbF, \Lambda)$ and $\al$ is an isomorphism $M\cong \phi_X^*M$.
Note that the composition $\phi_X\circ \phi_{\bbF}$ is the absolute $q$-Frobenius of $X_\bbF$. 
In particular, it induces the identity on pro\'etale topoi, see \cite[Lemma 5.4.2]{BhattScholze:ProEtale}. 
Therefore, replacing $\phi_X^*$ by $\phi_{\bbF}^*$ in \refeq{Weil.sheaves} 
yields an equivalent category.
\xrema

\pf[Proof of \thref{presentation.Weil.prop}]
The structural morphism $(X_0,\Phi_0) \rightarrow (X_\bbF,\phi_X)$ is a cover in $X_{\proet}^\Weil$. 
Its \v{C}ech nerve has objects $(X_i,\Phi_i)\in X_{\proet}^{\Weil}$, $i\geq 0$ as above.
By descent, there is an equivalence
\eqn
\label{descent.topoi.twisted.bar}
\D\big(X^\Weil,\Lambda\big) \stackrel \cong \lr \tot \left(\D\big( (X_{\proet}^{\Weil})_{/(X_{\bullet},\Phi_\bullet)},\Lambda\big)\right).
\xeqn
Under \thref{lemma.trivialization.equivariant.Weil.site}, the cosimplicial $1$-topos associated with $(X_{\proet}^{\Weil})_{/(X_{\bullet},\Phi_\bullet)}$ is equivalent to the cosimplicial $1$-topos associated with the action of $\phi_X^*$ on $(X_\bbF)_\proet$.
The equivalence \eqref{descent.topoi.twisted.bar} then becomes
\begin{equation*}
\D\big(X^\Weil,\Lambda\big) \stackrel \cong \lr \lim_{B\bbZ} \D(X_{\bbF},\Lambda),
\end{equation*}
for the diagram $B\bbZ\to \PrStL$ corresponding to the endomorphism $\phi^*_X$ of $ \D\big(X_{\bbF},\Lambda\big)$.
That is, $\D(X^\Weil,\Lambda)$ is equivalent to the homotopy fixed points of $\D(X_{\bbF},\Lambda)$ with respect to the action of $\phi_X^{*}$, which is our claim.
\xpf

\subsection{Weil sheaves on products}
\label{sect--Weil.sheaves.product}
The discussion of the previous section generalizes to products of schemes as follows. 
Let $X_1,\ldots, X_n$ be schemes over $\Fq$, and denote by $X:=X_1\x_{\Fq}\ldots\x_{\Fq}X_n$ their product. 
For every $1\leq i\leq n$, we have a morphism $\phi_{X_{i}}\colon X_{i,\bbF} \rightarrow X_{i,\bbF}$ as in the previous section. 
We use the notation $\phi_{X_{i}}$ to also denote the corresponding map on $X_\bbF=X_{1,\bbF}\times_{\bbF}\ldots \times_{\bbF} X_{n,\bbF}$ which is $\phi_{X_{i}}$ on the $i$-th factor and the identity on the other factors. 

We define the site \((X_1^\Weil\x\ldots\x X_n^\Weil)_{\proet}\) whose underlying category consists of tuples $(U,\varphi_1,\dots,\varphi_n)$ with $U\in (X_\bbF)_{\proet}$ and pairwise commuting endomomorphisms $\varphi_i \colon U\rightarrow U$ such that the following diagram commutes
\[
\xymatrix{
U \ar[r]^{\varphi_i}\ar[d] &  U\ar[d]\\
X_\bbF \ar[r]^{\phi_{X_i}} & X_\bbF,
}
\]
 for all $1\leq i\leq n$. 
 As before, we denote by $\D\big(X_1^\Weil\x\ldots\x X_n^\Weil,\Lambda\big)$ the corresponding derived category of $\Lambda$-sheaves.
 
Using a similar reasoning as in the previous section, we can identify this category of sheaves with the homotopy fixed points
\eqn
\label{eqn--derived.cat.weil.sheaves.products.limit}
\D\big(X_1^\Weil\x\ldots\x X_n^\Weil,\Lambda\big) \stackrel\cong\lr \Fix \big(\D(X_{\bbF},\Lambda), \phi_{X_1}^*, \dots, \phi_{X_n}^*\big)
\xeqn
of the commuting family of the functors $\phi_{X_i}^*$, see \thref{simultaneous.Fix}. 
Explicitly, for $n=2$, this is the homotopy limit of the diagram 
\[
\xymatrix{
\D\big(X_{\bbF},\Lambda\big) \ar@<1ex>[r]^{\phi_{X_1}^*} \ar@<-1ex>[r]_{\id} \ar@<-1ex>[d]_{\id} \ar@<1ex>[d]^{\phi_{X_2}^*} & 
\D\big(X_{\bbF},\Lambda\big) \ar@<-1ex>[d]_{\id} \ar@<1ex>[d]^{\phi_{X_2}^*}
\\
\D\big(X_{\bbF},\Lambda\big) \ar@<1ex>[r]^{\phi_{X_1}^*} \ar@<-1ex>[r]_{\id} & 
\D\big(X_{\bbF},\Lambda\big).
}
\]

Roughly speaking, objects in the category $\D(X_1^\Weil\x\ldots\x X_n^\Weil,\Lambda)$ are given by tuples $(M,\al_1,\ldots, \al_n)$ with $M \in \D(X_{\bbF}, \Lambda)$ and with pairwise commuting equivalences $\al_i\co M\cong \phi_{X_i}^*M$. 
That is, equipped with a collection of equivalences $\phi_{X_j}^*(\al_i)\circ \al_j\simeq \phi_{X_i}^*(\al_j)\circ \al_i$ for all $i,j$ satisfying higher coherence conditions.

\subsection{Partial-Frobenius stability}
For schemes $X_1, \dots, X_n$ over $\Fq$, we denote by $X := X_1 \x_{\Fq} \dots \x_{\Fq} X_n$ their product together with the partial Frobenii $\Frob_{X_i}\co X\r X$, $1\leq i\leq n$.
To give a reasonable definition of lisse and constructible Weil sheaves, we need to understand the relation between partial-Frobenius invariant constructible subsets in $X$ and constructible subsets in the single factors $X_i$:

\defi
A subset $Z \subset X$ is called \emph{partial-Frobenius invariant} if $\Frob_{X_i}(Z) = Z$ for all $1\leq i\leq n$.
\xdefi

The composition $\Frob_{X_1} \circ \dots \circ \Frob_{X_n}$ is the absolute $q$-Frobenius on $X$ and thus induces the identity on the topological space underlying $X$.
Therefore, in order to check that $Z \subset X$ is partial-Frobenius invariant, it suffices that, for any fixed $i$, the subset $Z$ is $\Frob_{X_j}$-invariant for all $j \ne i$.
This remark, which also applies to $X_\bbF = X_1 \x_{\Fq} \dots \x_{\Fq} X_n \x_{\Fq} \Spec \bbF$, will be used below without further comment.

We first investigate the case of two factors with one being a separably closed field. 
This eventually rests on Drinfeld's descent result \cite[Proposition 1.1]{Drinfeld:FSheaves} for coherent sheaves:

\lemm
\thlabel{partial.Frobenius1}
Let $X$ be a qcqs $\Fq$-scheme, and let $k/\Fq$ be a separably closed field.
Denote by $p\co X_k\to X$ the projection.
Then $Z\mapsto p^{-1}(Z)$ induces a bijection
$$\{\text{constructible subsets in $X$}\,\} \leftrightarrow \{\text{partial-Frobenius invariant, constructible subsets in $X_k$}\,\}.$$ 
\xlemm

\pf
The injectivity is clear because $p$ is surjective. 
It remains to check the surjectivity.
Without loss of generality we may assume that $k$ is algebraically closed, and replace $\Frob_X$ by $\Frob_k$ which is an automorphism. 
Given that $Z\mapsto p^{-1}(Z)$ is compatible with passing to complements, unions and localizations on $X$,
we are reduced to proving the bijection for constructible closed subsets $Z$ and for $X$ affine over $\Fq$.
By Noetherian approximation (\thref{Noetherian.approximation.lemm}), we reduce further to the case where $X$ is of finite type over $\Fq$ and still affine.
Now we choose a locally closed embedding $X\to \bbP^n_{\Fq}$ into projective space. 
A closed subset $Z'\subset X_k$ is $\phi_k$-invariant if and only if its closure inside $\bbP^n_{k}$ is so. 
Hence, it is enough to consider the case where $X=\bbP^n_{\Fq}$ is the projective space. 
Let $Z'$ be a closed $\Frob_k$-invariant subset of $X_{k}$. 
When viewed as a reduced subscheme, the isomorphism $\phi_k$ restricts to an isomorphism of $Z'$. 
In particular, $\calO_{Z'}$ is a coherent $\calO_{X_k}$-module equipped with an isomorphism $\calO_{Z'}\cong \phi_k^*\calO_{Z'}$.
Hence, Drinfeld's descent result \cite[Proposition 1.1]{Drinfeld:FSheaves} (see also \cite[Section 4.2]{Kedlaya:Shtukas} for a recent exposition) yields $Z'=Z_k$ for a unique closed subscheme $Z\subset X$.
\xpf

The following proposition generalizes the results \cite[Lemma~9.2.1]{Lau:Shtukas} and \cite[Lemme~8.12]{Lafforgue:Chtoucas} in the case of curves.

\prop
\thlabel{partial.Frobenius2}
Let $X_1,\ldots, X_n$ be qcqs $\bbF_q$-schemes, and denote $X=X_1\x_{\bbF_q}\ldots\x_{\bbF_q}X_n$.
Then any partial-Frobenius invariant constructible closed subset $Z \subset X$ is a finite set-theoretic union of subsets of the form $Z_1 \x_{\bbF_q} \dots \x_{\bbF_q} Z_n$, for appropriate constructible closed subschemes $Z_i \subset X_i$.

In particular, any partial-Frobenius invariant constructible open subscheme $U \subset X$ is a finite union of constructible open subschemes of the form $U_1 \x_{\Fq} \dots \x_{\Fq} U_n$, for appropriate constructible open subschemes $U_i \subset X_i$.
\xprop

\pf
By induction, we may assume $n=2$.
By Noetherian approximation (\thref{Noetherian.approximation.lemm}), we reduce to the case where both $X_1,X_2$ are of finite type over $\bbF_q$. 
In the following, all products are formed over $\bbF_q$, and locally closed subschemes are equipped with their reduced subscheme structure. 
Let $Z\subset X_1\x X_2$ be a partial-Frobenius invariant closed subscheme. 
The complement $U = X_1 \x X_2 \setminus Z$ is also partial-Frobenius invariant.

In the proof, we can replace $X_1$ (and likewise $X_2$) by a stratification in the following sense: 
Suppose $X_1 = A' \sqcup A''$ is a set-theoretic stratification into a closed subset $A'$ with open complement $A''$.
Once we know $Z \cap A' \x X_2 = \bigcup_j Z'_{1j} \x Z'_{2j}$ and $Z \cap A'' \x X_2 = \bigcup_j Z''_{1j} \x Z''_{2j}$ for appropriate closed subschemes $Z'_{1j} \subset A'$, $Z''_{1j} \subset A''$ and $Z'_{2j}, Z''_{2j} \subset X_2$, we have the set-theoretic equality
$$Z = \bigcup_j Z'_{1i} \x Z'_{2j} \cup \bigcup_j \ol{Z''_{1j}} \x Z''_{2j},$$
where $\ol{Z''_{1j}}\subset X_1$ denotes the scheme-theoretic closure.
Here we note that taking scheme-theoretic closures commutes with products because the projections $X_1 \x X_2 \r X_i$ are flat, and that the topological space underlying the scheme-theoretic closure agrees with the topological closure because all schemes involved are of finite type.

The proof is now by Noetherian induction on $X_2$, the case $X_2 = \emptyset$ being clear (or, if the reader prefers the case where $X_2$ is zero dimensional reduces to \thref{partial.Frobenius1}).
In the induction step, we may assume, using the above stratification argument, that both $X_i$ are irreducible with generic point $\eta_i$.
We let $\ol \eta_i$ be a geometric generic point over $\eta_i$, and denote by $p_i\co X_1 \x X_2 \r X_i$ the two projections.
Both $p_i$ are faithfully flat of finite type and in particular open, so that $p_i(U)$ is open in $X_i$.
We have a set-theoretic equality
$$Z = \big ( (X_1 \setminus p_1(U)) \x X_2 \big ) \cup \big ( X_1 \x (X_2 \setminus p_2(U)) \big ) \cup \big ( Z \cap p_1(U) \x p_2(U) \big ).$$
Once we know $Z \cap p_1(U) \x p_2(U) = \bigcup_j Z_{1j} \x Z_{2j}$ for appropriate closed $Z_{ij} \subset p_i(U)$, we are done.
We can therefore replace $X_i$ by $p_i(U)$ and assume that both $p_i\co U \r X_i$ are surjective.

The base change $U\x_{X_2} \ol \eta_2$ is a $\phi_{\ol \eta_2}$-invariant subset of $X_1 \x \ol \eta_2$. 
By Lemma \ref{partial.Frobenius1}, it is thus of the form $U_1\x \ol \eta_2$ for some open subset $U_1\subset X_1$.
There is an inclusion (of open subschemes of $X_1 \x \eta_2$): $U \x_{X_2} \eta_2 \subset U_1 \x \eta_2$.
It becomes a set-theoretic equality, and therefore an isomorphism of schemes, after base change along $\ol \eta_2 \r \eta_2$. 
By faithfully flat descent, 
this implies that the two mentioned subsets of $X_1 \x \eta_2$ agree. 
We claim $U_1 = X_1$.
Since the projection $U \r X_2$ is surjective, in particular its image contains $\eta_2$, so that $U_1$ is a non-empty subset, and therefore open dense in the irreducible scheme $X_1$.
Let $x_1\in X_1$ be a point.
Since the projection $U \r X_1$ is surjective, $U\cap (\{x_1\} \x X_2)$ is a non-empty open subscheme of $\{x_1\}\x X_2$.
So it contains a point lying over $(x_1,\eta_2)$.
We conclude $X_1 \x \eta_2 \subset U$.

We claim that there is a non-empty open subset $A_2 \subset X_2$ such that 
$$X_1 \x A_2 \subset U \text{ or, equivalently, } X_1 \x (X_2 \setminus A_2) \supset X_1 \x X_2 \setminus U.$$
The underlying topological space of $V = X_1 \x X_2 \setminus U$ is Noetherian and thus has finitely many irreducible components $V_j$.
The closure of the projection $\overline{p_2(V_j)} \subset X_2$ does not contain $\eta_2$, since $X_1 \x \eta_2 \subset U$.
Thus, $A_2 := \bigcap_j X_1 \setminus \overline{p_2(V_j)}$ satisfies our requirements.

Now we continue by Noetherian induction applied to the stratification $X_2 = A_2 \sqcup (X_2 \setminus A_2)$:
We have $Z \cap X_1 \x A_2 = \emptyset$, so that we may replace $X_2$ by the proper closed subscheme $X_2 \setminus A_2$.
Hence, the proposition follows by Noetherian induction.
\xpf

The following lemma on Noetherian approximation of partial Frobenius invariant subsets is needed for the reduction to finite type schemes:

\lemm
\thlabel{Noetherian.approximation.lemm}
Let $X_1,\ldots, X_n$ be qcqs $\bbF_q$-schemes, and denote $X=X_1\x_{\bbF_q}\ldots\x_{\Fq}X_n$.
Let $X_i=\lim_j X_{ij}$ be a cofiltered limit of finite type $\bbF_q$-schemes with affine transition maps, and write $X=\lim_j X_j$, $X_j:=X_{1j}\x_{\bbF_q}\ldots\x_\Fq X_{nj}$ (see \StP{01ZA} for the existence of such presentations). 
Let $Z\subset X$ be a constructible closed subset.
Then the intersection
\[
Z'=\bigcap_{i=1}^n \bigcap_{m\in \bbZ} \Frob_{X_i}^m(Z)
\]
is partial Frobenius invariant, constructible closed and there exists an index $j$ and a partial Frobenius invariant closed subset $Z'_{j}\subset X_{j}$ such that $Z'=Z'_{j}\x_{X_{j}}X$ as sets. 
\xlemm

We note that each $\Frob_{X_i}$ induces a homeomorphism on the underlying topological space of $X$ so that $Z'$ is well-defined. 
This lemma applies, in particular, to partial Frobenius invariant constructible closed subsets $Z\subset X$ in which case we have $Z=Z'$. 

\pf
As $Z$ is constructible, there exists an index $j$ and a constructible closed subscheme $Z_j\subset X_j$ such that $Z=Z_j\x_{X_j}X$ as sets.
We put $Z'_j=\cap_{i=1}^n \cap_{m\in \bbZ} \Frob_{X_{ij}}^m(Z_j)$.
As $X_j$ is of finite type over $\bbF_q$, the subset $Z'_j$ is still constructible closed. 
As partial Frobenii induce bijections on the underlying topological spaces, one checks that $\Frob_{X_{ij}}^m(Z_j)\x_{X_j}X=\Frob_{X_i}^m(Z)$ as sets for all $m\in \bbZ$.
Thus, $Z'=Z_j'\x_{X_j}X$ which, also, is constructible closed because $X\r X_j$ is affine. 
\xpf

\subsection{Lisse and constructible Weil sheaves}
In this subsection, we define the subcategories of lisse and constructible Weil sheaves and establish a presentation similar to \refeq{Weil.sheaves}.
Let $X_1, \dots, X_n$ be schemes over $\Fq$, and denote $X := X_1 \x_{\Fq} \dots \x_{\Fq} X_n$. 
Let $\Lambda$ be a condensed ring.

\defi 
\thlabel{definition-Weil-on-products}
Let $M\in \D\big(X_1^\Weil\x\ldots\x X_n^\Weil,\Lambda\big)$.
\begin{enumerate}
    \item The Weil sheaf $M$ is called {\it lisse} if it is dualizable.
   (Here dualizability refers to the symmetric monoidal structure on $\D\big(X_1^\Weil\x\ldots\x X_n^\Weil,\Lambda\big)$, given by the derived tensor product of $\Lambda$-sheaves on the Weil-pro\'etale topos.)    
    
    \item 
    \label{item--definition-Weil-on-products.2}
    The Weil sheaf $M$ is called {\it constructible} if for any open affine $U_i\subset X_i$ there exists a finite subdivision into constructible locally closed subschemes $U_{ij}\subseteq  U_i$ such that each restriction $M|_{U_{1j}^\Weil\x \ldots \x U_{nj}^\Weil }\in \D\big(U_{1j}^\Weil\times\ldots \times U_{nj}^\Weil ,\Lambda\big)$ is lisse.
\end{enumerate}
\xdefi

The full subcategories of $\D\big(X_1^\Weil\x\ldots\x X_n^\Weil,\Lambda\big)$ consisting of lisse, resp.~constructible Weil sheaves are denoted by
\[
\D_\lis\big(X_1^\Weil\x\ldots\x X_n^\Weil,\Lambda\big) \subset \D_\cons\big(X_1^\Weil\x\ldots\x X_n^\Weil,\Lambda\big).
\]
Both categories are idempotent complete stable $\Gamma(X,\Lambda)$-linear symmetric monoidal \ii-categories.

From the presentation \refeq{derived.cat.weil.sheaves.products.limit}, we get that a Weil sheaf $M$ is lisse if and only if the underlying object $M_{\bbF}\in \D(X_\bbF,\Lambda)$ is lisse. 
So \refeq{derived.cat.weil.sheaves.products.limit} restricts to an equivalence
\eqn
\label{eq.lisse.weil.product.limit}
\D_\lis\big(X_1^\Weil\x\ldots\x X_n^\Weil,\Lambda\big) \cong \Fix \left(\D_\lis\big(X_\bbF,\Lambda\big), \phi_{X_1}^*, \dots, \phi_{X_n}^*\right).
\xeqn
The same is true for constructible Weil sheaves by the following proposition: 

\prop\thlabel{prop.cons.weil.as.fixed.points}
A Weil sheaf $M\in \D\big(X_1^\Weil\x\ldots\x X_n^\Weil,\Lambda\big)$ is constructible if and only if the underlying sheaf $M_{\bbF}\in \D(X_\bbF,\Lambda)$ is constructible. 
Consequently, \refeq{derived.cat.weil.sheaves.products.limit} restricts to an equivalence
\eqn
\label{eq.construtible.weil.product.limit}
\D_\cons\big(X_1^\Weil\x\ldots\x X_n^\Weil,\Lambda\big) \cong \Fix \left(\D_\cons\big(X_\bbF,\Lambda\big), \phi_{X_1}^*, \dots, \phi_{X_n}^*\right).
\xeqn
\xprop
\pf 
Clearly, if $M$ is constructible, so is $M_\bbF$ by \thref{definition-Weil-on-products}.
Let $M\in \D\big(X_1^\Weil\x\ldots\x X_n^\Weil,\Lambda\big)$ such that $M_\bbF$ is constructible. 
We may assume that all $X_{i}$ are affine.  
We claim that there is a finite subdivision $X_\bbF=\sqcup X_\al$ into constructible locally closed subsets such that $M_\bbF|_{X_\al}$ is lisse and such that each $X_\al$ is partial Frobenius invariant. 

Assuming the claim we finish the argument as follows.
By \thref{partial.Frobenius2}, any open stratum $U=X_{j_0}\subset X_\bbF$ is a finite union of subsets of the form $U_{1,\bbF}\x_\bbF\ldots\x_\bbF U_{n,\bbF}$ and the restriction of $M$ to each of them is lisse. 
In particular, the complement $X_\bbF\backslash U$ is defined over $\bbF_q$ and arises as a finite union of schemes of the form $X'=X_1'\x_\Fq\ldots\x_\Fq X_n'$ for suitable qcqs schemes $X_i'$ over $\Fq$.
Intersecting each $X_\bbF'$ with the remaining strata $\sqcup_{j\neq j_0}X_j$, we conclude by induction on the number of strata. 

It remains to prove the claim.
We start with any finite subdivision $X_\bbF=\sqcup X'_j$ into constructible locally closed subsets such that $M_\bbF|_{X'_j}$ is lisse.
Pick an open stratum $X'_{j_0}$, and set
\eqn
\label{partial.Frobenius.orbit.eq}
X_{j_0}=\bigcup_{i=1}^n\bigcup_{m\in \bbZ}\phi_{X_i}^m(X_{j_0}').
\xeqn
This is a constructible open subset of $X_\bbF$ by \thref{Noetherian.approximation.lemm} applied to its closed complement.
Furthermore, $M_\bbF|_{X_{j_0}}$ is lisse by its partial Frobenius equivariance, noting that $\phi_{X_i}^*$ induces equivalences on pro\'etale topoi to treat the negative powers in \eqref{partial.Frobenius.orbit.eq}.
As before, $X_\bbF\backslash X_{j_0}$ is defined over $\Fq$.
So replacing $X_j'$, $j\neq j_0$ by $X_j'\cap (X_\bbF\backslash X_{j_0})$, the claim follows by induction on the number of strata. 
\xpf

In the case of a single factor $X=X_1$, the preceding discussion implies
\eqn\label{eqn--lisse.weil.sheaves.equalizer}
\D_\bullet\big(X^\Weil,\Lambda\big)\ \cong \lim\left (\D_\bullet(X_{\bbF},\Lambda) \stackrel[\id]{\phi_X^*}\rightrightarrows \D_\bullet(X_{\bbF},\Lambda) \right),
\xeqn
for $\bullet \in \{\varnothing, \lis, \cons\}$.

\subsection{Relation with the Weil groupoid} 
In this subsection, we relate lisse Weil sheaves with representations of the Weil groupoid.
Throughout, we work with \'etale fundamental groups as opposed to their pro\'etale variants in order to have Drinfeld's lemma available, see \refsect{drinfelds.lemma}.  
The two concepts differ in general, but agree for geometrically unibranch (for example, normal) Noetherian schemes, see \cite[Lemma 7.4.10]{BhattScholze:ProEtale}.

For a Noetherian scheme $X$, let $\pi_1(X)$ be the \emph{étale fundamental groupoid} of $X$ as defined in \cite[Exposé V, \S7 and \S9]{SGA1}.
Its objects are geometric points of $X$, and its morphisms are isomorphisms of fiber functors on the finite étale site of $X$.
This is an essentially small category. 
The automorphism group in $\pi_1(X)$ at a geometric point $x \to X$ is profinite. 
It is denoted $\pi_1(X,x)$ and called the \emph{étale fundamental group} of $(X,x)$.
If $X$ is connected, then the natural map $B\pi_1(X,x)\to \pi_1(X)$ is an equivalence for any $x\to X$.
If $X$ is the disjoint sum of schemes $X_i$, $i\in I$, then $\pi_1(X)$ is the disjoint sum of the $\pi_1(X_i)$, $i\in I$.
In this case, if $x\to X$ factors through $X_i$, then $\pi_1(X,x)=\pi_1(X_i,x)$. 

\defi
\thlabel{FWeil}
Let $X_1,\dots, X_n$ be Noetherian schemes over $\bbF_q$, and write $X = X_1\times_{\bbF_q} \ldots \times_{\bbF_q} X_n$. 
The \textit{Frobenius-Weil groupoid} is the stacky quotient
\eqn 
\label{Frobenius_Weil_Intro}
\FWeil(X)= \pi_1(X_\bbF)/\lan\phi_{X_1}^\bbZ,\ldots, \phi_{X_n}^\bbZ\ran,
\xeqn 
where we use that the partial Frobenii $\phi_{X_i}$ induce automorphisms on the finite étale site of $X_{\bbF}$. 
\xdefi

For $n=1$, we denote $\FWeil(X)=\Weil(X)$.
Even if $X$ is connected, its base change $X_\bbF$ might be disconnected in which case the action of $\phi_X$ permutes some connected components. 
Therefore, fixing a geometric point of $X_\bbF$ is inconvenient, and the reason for us to work with fundmental groupoids as opposed to fundmental groups. 
The automorphism groups in $\Weil(X)$ carry the structure of locally profinite groups:
indeed, if $X$ is connected, then $\Weil(X)$ is, for any choice of a geometric point $x\to X_\bbF$, equivalent to the classifying space of the Weil group $\Weil(X,x)$ from \cite[Définition 1.1.10]{Deligne:Weil2}.
Recall that this group sits in an exact sequence of topological groups
\eqn
1\to \pi_1(X_\bbF,x)\to \Weil(X,x) \to \Weil(\bbF/\Fq)\simeq \bbZ,
\label{eqn--Weil.Z}
\xeqn
where $\pi_1(X_\bbF,x)$ carries its profinite topology and $\bbZ$ the discrete topology.
The topology on the morphism groups in $\Weil(X)$ obtained in this way is independent from the choice of $x\to X_\bbF$.
The image of $\Weil(X,x)\to \bbZ$ is the subgroup $m\bbZ$ where $m$ is the degree of the largest finite subfield in $\Gamma(X,\calO_X)$.
In particular, we have $m=1$ if $X_\bbF$ is connected.
Let us add that if $x\to X_\bbF$ is fixed under $\phi_X$, then the action of $\phi_X$ on $\pi_1(X_\bbF,x)$ corresponds by virtue of the formula $\phi_X^*=(\phi_\bbF^*)^{-1}$  to the action of the geometric Frobenius, that is, the inverse of the $q$-Frobenius in $\Weil(\bbF/\Fq)$.

Likewise, for every $n\geq 1$, the stabilizers of the Frobenius-Weil groupoid are related to the partial Frobenius-Weil groups introduced in \cite[Proposition 6.1]{Drinfeld:FSheaves} and \cite[Remarque 8.18]{Lafforgue:Chtoucas}. 
In particular, there is an exact sequence
\[
1 \rightarrow \pi_1(X_\bbF,x) \rightarrow \FWeil(X,x) \rightarrow  \bbZ^n,
\]
for each geometric point $x\r X_\bbF$.
This gives $\FWeil(X)$ the structure of a locally profinite groupoid. 

Let $\Lambda$ be either of the following coherent topological rings: a coherent discrete ring, an algebraic field extension $E\supset \bbQ_\ell$ for some prime $\ell$, or its ring of integers $\calO_E\supset \bbZ_\ell$. 
For a topological groupoid $W$, we will denote by $\Rep_{\Lambda}(W)$ the category of continuous representations of $W$ with values in finitely presented $\Lambda$-modules and by $\Rep^\fp_{\Lambda}(W) \subset \Rep_{\Lambda}(W)$ its full subcategory of representations on finite projective $\Lambda$-modules. 
Here finitely presented $\Lambda$-modules $M$ carry the quotient topology induced from the choice of any surjection $\Lambda^n\to M$, $n\geq 0$ and the product topology on $\Lambda^n$. 

\lemm
\thlabel{elemntary.rep.lemm}
In the situation above, the category $\Rep_\Lambda(W)$ is $\Lambda_*$-linear and abelian.
In particular, its full subcategory $\Rep^\fp_{\Lambda}(W)$ is $\Lambda_*$-linear and additive. 
\xlemm
\pf
Let $W_\disc$ be the discrete groupoid underlying $W$, and denote by $\Rep_\Lambda(W_\disc)$ the category of $W_\disc$-representations on finitely presented $\Lambda$-modules. 
Evidently, this category is $\Lambda_*$-linear. It is abelian since $\Lambda$ is coherent (\thref{sheaves}~\refit{t-structure}, see also \consref{coherent.ring.Mod.fp.abelian}. 
We claim that $\Rep_\Lambda(W)\subset \Rep_\Lambda(W_\disc)$ is a $\Lambda_*$-linear full abelian subcategory.  
If $\Lambda$ is discrete (and coherent), then every finitely presented $\Lambda$-module carries the discrete topology and the claim is immediate, see also \StP{0A2H}.
For $\Lambda=E,\calO_E$, one checks that every map of finitely presented $\Lambda$-modules is continuous, every surjective map is a topological quotient and every injective map is a closed embedding.
For the latter, we use that every finitely presented $\Lambda$-module can be written as a countable filtered colimit of compact Hausdorff spaces along injections, and that every injection of compact Hausdorff spaces is a closed embedding.  
This implies the claim. 
\xpf

We apply this for $W$ being either of the locally profinite groupoids $\pi_1(X)$, $\pi_1(X_\bbF)$ or $\FWeil(X)$. 
Note that restricting representations along $\pi_1(X_{\bbF})\r \FWeil(X)$ induces an equivalence of $\Lambda_*$-linear abelian categories
\eqn\label{eqn-Weil.groupoid.rep.fixed.points}
\Rep_{\Lambda}\big(\FWeil(X)\big) \cong \Fix \left(\Rep_{\Lambda}\big(\pi_1(X_{\bbF})\big), \phi_{X_1}, \dots, \phi_{X_n}\right),
\xeqn
and similarly for the $\Lambda_*$-linear additive category $\Rep^{\fp}_{\Lambda}\big(\FWeil(X)\big)$. 

\defi
\thlabel{bounded.complexes.nota}
For an integer $n\geq 0$, we write $\D_\lis^{\{-n,n\}}(X,\Lambda)$ for the full subcategory of $\D_\lis(X,\Lambda)$ of objects $M$ such that $M$ and its dual $M^\vee$ lie in degrees $[-n,n]$ with respect to the t-structure on $\D(X, \Lambda)$.
\xdefi

\lemm \thlabel{FWeil.heart.and.almost.compact.objects}
In the situation above, there is a natural functor 
    \eqn\label{eqn-weil.sheaves.map}
    \Rep_\Lambda\big(\FWeil(X)\big) \r \D\big(X_1^\Weil\times \ldots \times X_n^\Weil,\Lambda\big)^{\heartsuit},
    \xeqn
that is fully faithful. 
Moreover, the following properties hold if $\Lambda$ is either finite discrete or $\Lambda=\calO_E$ for $E\supset \bbQ_\ell$ finite: 
\begin{enumerate}
    \item \label{item--weil.rep.essential.image}
    An object $M$ lies in the essential image of \eqref{eqn-weil.sheaves.map} if and only if its underlying sheaf $M_\bbF$ is locally on $(X_\bbF)_\proet$ isomorphic to $\underline N\t_{\underline{\Lambda_*}}\Lambda_{X_\bbF}$ for some finitely presented $\Lambda_*$-module $N$. 

    \item \label{item--weil.rep.amplitude} 
    The functor \eqref{eqn-weil.sheaves.map} restricts to an equivalence of $\Lambda_*$-linear additive categories
    \[
    \Rep^\fp_\Lambda\big(\FWeil(X)\big) \stackrel\cong \lr  \D^{\{0,0\}}_\lis\big(X_1^\Weil\x \ldots \x X_n^\Weil ,\Lambda\big).
    \]
    \item \label{item--weil.rep.t.structure} 
    
    If $\Lambda_*$ is regular (so that $\Lambda$ is t-admissible, cf.~\thref{sheaves}~\refit{t-structure}), then \eqref{eqn-weil.sheaves.map} restricts to an equivalence of $\Lambda_*$-linear abelian categories
    \[
    \Rep_\Lambda\big(\FWeil(X)\big) \stackrel\cong \lr \D_\lis\big(X_1^\Weil\x \ldots \x X_n^\Weil ,\Lambda\big)^\heartsuit.
    \] 
\end{enumerate}
    If all $X_i$, $i=1,\ldots, n$ are geometrically unibranch, then \refit{weil.rep.essential.image}, \refit{weil.rep.amplitude} and \refit{weil.rep.t.structure} hold for general coherent topological rings $\Lambda$ as above.
\xlemm 
\pf  
There is a canonical equivalence of topological groupoids $\pi_1(X_\bbF)\cong\widehat{\pi_1^\proet(X_\bbF)}$ with the profinite completion of the pro\'etale fundamental groupoid, see \cite[Lemma 7.4.3]{BhattScholze:ProEtale}.
It follows from \cite[Lemmas 7.4.5, 7.4.7]{BhattScholze:ProEtale} that restricting representations along $\pi_1^\proet(X_\bbF)\r \pi_1(X_\bbF)$ induces full embeddings
\eqn\label{eqn-functor.pi_1.rep.to.sheaves}
\Rep_{\Lambda}\big(\pi_1(X_\bbF)\big) \hr  \Rep_{\Lambda}\big(\pi_1^\proet(X_\bbF)\big)\hr \D(X_{\bbF},\Lambda)^\heartsuit,
\xeqn
that are compatible with the action of $\phi_{X_i}$ for all $i=1,\ldots,n$. 
So we obtain the fully faithful functor \eqref{eqn-weil.sheaves.map} by passing to fixed points, see \eqref{eqn-Weil.groupoid.rep.fixed.points}, \eqref{eq.lisse.weil.product.limit} and \thref{Fix.t-structure} (see also \thref{simultaneous.Fix}).

Part \refit{weil.rep.essential.image} describes the essential image of $\Rep_{\Lambda}\big(\pi_1^\proet(X_\bbF)\big)\hr \D(X_{\bbF},\Lambda)^\heartsuit$.
So if $\Lambda$ is finite discrete or profinite, then the first functor in \eqref{eqn-functor.pi_1.rep.to.sheaves} is an equivalence, and we are done.
Part \refit{weil.rep.amplitude} is immediate from \refit{weil.rep.essential.image}, noting that an object in the essential image of \eqref{eqn-functor.pi_1.rep.to.sheaves} is lisse if and only if its underlying module is finite projective.
Likewise, part \refit{weil.rep.t.structure} is immediate from \refit{weil.rep.essential.image}, using \thref{sheaves}~\refit{lisse.locally.constant}.
Here we need to exclude rings like $\Lambda=\bbZ/\ell^2$ in order to have a t-structure on lisse sheaves.

Finally, if all $X_i$ are geometrically unibranch, so is $X_\bbF$ which follows from the characterization \StP{0BQ4}. 
In this case, we get $\pi_1(X_\bbF)\cong\pi_1^\proet(X_\bbF)$ by \cite[Lemma 7.4.10]{BhattScholze:ProEtale}.
This finishes the proof. 
\xpf

\subsection{Weil-\'{e}tale versus \'{e}tale sheaves}
We end this section with the following description of Weil sheaves with (ind-)finite coefficients.
Note that such a simplification in terms of ordinary sheaves is not possible for $\Lambda = \bbZ, \bbZ_\ell, \bbQ_\ell$, say.

\prop
\thlabel{discrete_Geisser}
Let $X$ be a qcqs $\Fq$-scheme.
Let $\Lambda$ be a finite discrete ring or a filtered colimit of such rings. 
Then the natural functors
$$\D_\lis(X, \Lambda) \r \D_\lis\big(X^\Weil, \Lambda\big), \;\; \D_\cons(X, \Lambda) \r \D_\cons\big(X^\Weil, \Lambda\big),$$
are equivalences.
\xprop
\pf 
Throughout, we repeatedly use that filtered colimits commute with finite limits in $\Cat_\infty$. 
Using compatibility of $\Dcons$ with filtered colimits in $\Lambda$ (\thref{sheaves}~\refit{manipulations.Lambda}), we may assume that $\Lambda$ is finite discrete. 
By the comparison result with the classical bounded derived category of constructible sheaves \consref{discrete.comparison.prop}, we can identify the categories $\D_\bullet(X,\Lambda)$, resp.~$\D_\bullet(X_\bbF,\Lambda)$ for $\bullet\in \{\lis,\cons\}$ with full subcategories of the derived category of \'etale $\Lambda$-sheaves $\D(X_\et,\Lambda)$, resp.~$\D(X_{\bbF,\et},\Lambda)$.
Write $X = \lim X_{i}$ as a cofiltered limit of finite type $\bbF_q$-schemes $X_{i}$ with affine transition maps \StP{01ZA}. 
Using the continuity of \'etale sites \StP{03Q4}, there are natural equivalences
\eqn
\label{Noetherian.approx.discrete.coefficients}
\colim\D_\bullet(X_i,\Lambda) \stackrel\cong \lr \D_\bullet(X,\Lambda), \quad \colim\D_\bullet(X_i^\Weil,\Lambda) \stackrel\cong \lr \D_\bullet(X^\Weil,\Lambda)
\xeqn
for $\bullet\in \{\lis,\cons\}$.
Hence, we can assume that $X$ is of finite type over $\bbF_q$.

To show full faithfulness, we claim more generally that the natural map
\[
\D(X_\et,\Lambda) \r \lim\left (\D(X_{\bbF,\et},\Lambda) \stackrel[\id]{\phi_X^*}\rightrightarrows \D(X_{\bbF,\et},\Lambda) \right)=: \D\big(X_{\et}^\Weil,\Lambda\big)
\]
is fully faithful. 
As $\Lambda$ is torsion, this is immediate from \cite[Corollary~5.2]{Geisser:Weil} applied to the inner homomorphisms between sheaves. 
Let us add that this induces fully faithful functors 
\eqn\label{eqn-Geisser.et.to.weil.all.sheaves}
\D^+(X_\et,\Lambda) \r \D^+\big(X_{\et}^\Weil,\Lambda\big) \r \D\big(X^\Weil,\Lambda\big)
\xeqn
on bounded below objects, see \cite[Proposition~5.2.6 (1)]{BhattScholze:ProEtale}. 

It remains to prove essential surjectivity. 
Using a stratification as in \thref{definition-Weil-on-products}, it is enough to consider the lisse case. 
Pick $M\in \D_\lis(X^\Weil,\Lambda)$. 
It is enough to show that $M$ lies is in the essential image of \eqref{eqn-Geisser.et.to.weil.all.sheaves}, noting that the functor detects dualizability. 
As $M$ is bounded, this will follow from showing that for every $j\in \bbZ$, the cohomology sheaf $\H^j(M) \in \D(X^\Weil,\Lambda)^\heartsuit$ is in the essential image of \eqref{eqn-Geisser.et.to.weil.all.sheaves}.

Fix $j\in \bbZ$. 
As $M$ is lisse, the underlying sheaf $\H^j(M)_{\bbF}\in \D(X_\bbF,\Lambda)^\heartsuit$ is pro\'etale-locally constant (\thref{sheaves}~\refit{lisse.locally.constant}) and valued in finitely presented $\Lambda$-modules. 
By \thref{FWeil.heart.and.almost.compact.objects} \refit{weil.rep.essential.image}, it comes from a representation of $\Weil(X)$. 
Restriction of representations along $\Weil(X) \rightarrow \pi_1(X)$ fits into a commutative diagram 
\[
\xymatrix{
\Rep_{\Lambda}\big(\pi_1(X)\big)\ar[d]\ar[r]^{\cong} &
\Rep_{\Lambda}\big(\Weil(X)\big)\ar[d] \\
\D(X_\et,\Lambda)^\heartsuit \ar[r] & \D\big(X^\Weil,\Lambda\big)^\heartsuit,
}
\]
where the upper horizontal arrow is an equivalence since $\Lambda$ is finite. In particular, the object $\H^j(M)$ is in the essential image of the fully faithful functor \eqref{eqn-Geisser.et.to.weil.all.sheaves}. 
\xpf 
\section{The categorical K\"{u}nneth formula}
We continue with the notation of \refsect{section.weil.sheaves}.
In particular, $\bbF_q$ denotes a finite field of characteristic $p>0$.
Recall from \refsect{recollections} the tensor product of $\Lambda_*$-linear idempotent complete stable \ii-categories. 
The external tensor product of sheaves $(M_1,\ldots,M_n)\mapsto M_1\bx\ldots\bx M_n$ as in \eqref{tensor.product} induces a functor
\eqn
\label{Drinfelds_map}
\D_\bullet\big(X_1^\Weil,\Lambda\big)\t_{\Perf_{\Lambda_*}}\ldots\t_{\Perf_{\Lambda_*}} \D_\bullet\big(X_n^\Weil,\Lambda\big) \r \D_\bullet\big(X_1^\Weil\x\ldots\x X_n^\Weil,\Lambda\big),
\xeqn
for $\bullet\in \{\lis,\cons\}$. 
Throughout, we consider the following situation.
In \thref{limits.Drinfeld.rema} we explain the compatibility of \eqref{Drinfelds_map} with certain (co-)limits in the schemes $X_i$ and coefficients $\Lambda$, which allows 
to relax these assumptions on $X$ and $\Lambda$ somewhat.

\situ
\thlabel{l.standard.situ.approximation}
The schemes $X_1,\dots,X_n$ are of finite type over $\bbF_q$, and $\Lambda$ is the condensed ring associated with one of the following topological rings:
\begin{enumerate}[(a)]
    \item \label{item--q.invertible.finite.Lambda} a finite discrete ring of prime-to-$p$-torsion;
    \item \label{item--q.invertible.adic.Lambda} the ring of integers $\calO_E$ of an algebraic field extension $E \supset\Ql$ for $\ell\neq p$ (for example $\bar \bbZ_\ell$);
    \item \label{item--q.invertible.rational.Lambda} an algebraic field extension $E \supset\Ql$ for $\ell\neq p$ (for example $\bar \bbQ_\ell$);
    \item \label{item--p.torsion.Lambda} a finite discrete $p$-torsion ring that is flat over $\bbZ/p^m$ for some $m\geq 1$.
\end{enumerate}
\xsitu

\theo
\thlabel{derived_Drinfeld.text}
In \thref{l.standard.situ.approximation}, the functor \eqref{Drinfelds_map} is an equivalence in each of the following cases:
\begin{enumerate}
    \item \label{item--cons.Drinfeld} $\bullet=\cons$ and $\Lambda$ is as in \refit{q.invertible.finite.Lambda}, \refit{q.invertible.adic.Lambda} or \refit{q.invertible.rational.Lambda}; 
    \item \label{item--lisse.Drinfeld} $\bullet=\lis$ and $\Lambda$ is as in \refit{q.invertible.finite.Lambda}, \refit{q.invertible.adic.Lambda}, \refit{p.torsion.Lambda} or as in \refit{q.invertible.rational.Lambda} if all $X_i$, $i=1,\ldots, n$ are geometrically unibranch (for example, normal).  
\end{enumerate}
\xtheo

In the $p$-torsion free cases \refit{q.invertible.finite.Lambda}, \refit{q.invertible.adic.Lambda} and \refit{q.invertible.rational.Lambda}, the full faithfulness is a direct consequence of the Künneth formula applied to the $X_{i,\bbF}$. 
In the $p$-torsion case \refit{p.torsion.Lambda}, we use Artin--Schreier theory instead.
It would be interesting to see whether this part can be extended to constructible sheaves using the mod-$p$-Riemann--Hilbert correspondence as in, say, \cite{BhattLurie:Riemann}.
In all cases, the essential surjectivity relies on a variant of Drinfeld's lemma for Weil group representations. 

Before turning to the proof of \thref{derived_Drinfeld.text}, we record the following compatibility of the functor \eqref{Drinfelds_map} with (co-)limits. 
This can be used to reduce the case of an (infinite) algebraic extension $E \supset \Ql$ in cases \refit{q.invertible.adic.Lambda} and \refit{q.invertible.rational.Lambda} above to the case where $E \supset \Ql$ is finite.
In the sequel we will therefore assume $E$ is finite in these cases.
\thref{limits.Drinfeld.rema} can further be used to extend \thref{derived_Drinfeld.text} to qcqs $\bbF_q$-schemes $X_i$ and finite discrete rings like $\bbZ/m$ for any integer $m\geq 1$ in cases \refit{q.invertible.finite.Lambda} and \refit{p.torsion.Lambda}.

\rema[Compatibility of \eqref{Drinfelds_map} with certain (co-)limits]
\thlabel{limits.Drinfeld.rema}
Throughout, we repeatedly use that filtered colimits commute with finite limits in $\CatEx_{\infty, \Lambda_*}(\Idem)$: the forgetful functors $\CatEx_{\infty, \Lambda_*}(\Idem) \r \CatEx_{\infty}(\Idem) \r \Cat_\infty$ create these (co)limits \cite[Theorem~1.1.4.4]{Lurie:HA}, \cite[Corollary 4.4.5.21]{Lurie:Higher}, 
and the statement holds in any compactly generated \ii-category, such as $\Cat_\infty$ \cite[Example~3.6(3)]{BhattMathew:Arc}.
We will also throughout use that in all the stable \ii-categories encountered below the tensor product preserves colimits and in particular finite limits.

\begin{enumerate}
\item
\label{item--limits.Drinfeld.rema.1}
\textit{Filtered colimits in $\Lambda$.}
First off, extension of scalars along any map of condensed rings $\Lambda \r \Lambda'$ induces a commutative diagram in $\CatEx_{\infty, \Lambda_*}(\Idem)$:
\[
\xymatrix{
\D_\bullet\big(X_1^\Weil,\Lambda\big)\t_{\Perf_{\Lambda_*}}\ldots\t_{\Perf_{\Lambda_*}} \D_\bullet\big(X_n^\Weil,\Lambda\big) \ar[r]\ar[d] & 
\D_\bullet(X^\Weil_{1} \ldots \x X^\Weil_{n},\Lambda)\ar[d]\\
\D_\bullet\big(X_1^\Weil,\Lambda'\big)\t_{\Perf_{\Lambda'_*}}\ldots\t_{\Perf_{\Lambda'_*}} \D_\bullet\big(X_n^\Weil,\Lambda'\big) \ar[r] &
\D_\bullet(X^\Weil_{1} \ldots \x X^\Weil_{n},\Lambda')
}
\]
It follows from the compatibility of $\Dcons$ with filtered colimits in $\Lambda$ (\thref{sheaves}~\refit{manipulations.Lambda}) that both sides of \eqref{Drinfelds_map} are compatible with filtered colimits in $\Lambda$. 

\item 
\label{item--limits.Drinfeld.rema.2}
\textit{Finite products in $\Lambda$.} 
Let $\Lambda=\prod \Lambda_i$ be a finite product of condensed rings. 
For any scheme $X$, the natural map $\D_\bullet(X,\Lambda)\r \prod \D_\bullet(X,\Lambda_i)$ is an equivalence for $\bullet\in \{\varnothing, \lis, \cons\}$, and likewise for Weil sheaves if $X$ is defined over $\bbF_q$. 
As $\Lambda_*=\prod \Lambda_{i,*}$, we see that \eqref{Drinfelds_map} is compatible finite products in the coefficients.

\item
\label{item--limits.Drinfeld.rema.3}
\textit{Limits in $X_i$ for discrete $\Lambda$.}
Assume that $\Lambda$ is finite discrete, see \thref{l.standard.situ.approximation} \refit{q.invertible.finite.Lambda}, \refit{p.torsion.Lambda}. 
Let $X_1,\ldots, X_n$ be qcqs $\bbF_q$-schemes.
Write each $X_i$ as a cofiltered limit $X_i = \lim X_{ij}$ of finite type $\bbF_q$-schemes $X_{ij}$ with affine transition maps \StP{01ZA}. 
As $\Lambda$ is finite discrete, we can use the continuity of \'etale sites as in \eqref{Noetherian.approx.discrete.coefficients} to show that the natural map
\[
\colim_{j} \D_\bullet\big(X^\Weil_{1j}\x \ldots \x X^\Weil_{nj},\Lambda\big) \stackrel\cong\lr \D_\bullet\big(X^\Weil_{1} \ldots \x X^\Weil_{n},\Lambda\big),
\]
is an equivalence for $\bullet \in \{\lis, \cons\}$. 
Thus, \eqref{Drinfelds_map} is compatible with cofiltered limits of finite type $\bbF_q$-schemes with affine transition maps. 
\end{enumerate}
\xrema


\subsection{A formulation in terms of prestacks}
Before turning to the proof, we point out a formulation of the results of the previous subsection in terms of symmetric monoidality of a certain sheaf theory. This formulation
makes the connection with constructions in the geometric approaches to the Langlands program  \cite{GKRV.Toy.Model, zhu2021coherent, LafforgueZhu:Decomposition} more manifest. Readers not familiar with prestacks
and formulations of sheaf theories on them can safely skip this section. The categories of constructible, resp.~lisse $\Lambda$-sheaves assemble into a lax symmetric monoidal functor
\eqn
\label{prestacks.sheaves.schemes}
\D_{\bullet, \Lambda}  \colon (\Sch_\bbF)^{\mathrm{op}} \r \CatEx_{\infty, \Lambda}(\Idem) \ (\bullet = \lis \text{ or } \cons). 
\xeqn
Namely, as a functor it sends a scheme $X$ to the category of constructible, resp.~lisse $\Lambda$-sheaves on $X$, and a morphism $f\colon X \r Y$ to the functor $f^*\colon \D_\bullet(Y, \Lambda) \r \D_\bullet(X, \Lambda)$.
These are objects, resp.~maps in the \ii-category $\CatEx_{\infty, \Lambda}(\Idem) := \Mod_{\Perf_\Lambda}(\CatEx_\infty(\Idem))$, cf.~\refsect{monoidal.aspects} for notation. 
The lax monoidal structure is given by the external tensor product of sheaves:
 \begin{align*}
  \boxtimes \colon \D_\bullet(X_\proet, \Lambda) \t_{\Perf_\Lambda} \D_\bullet(Y_\proet, \Lambda) &\rightarrow\D_\bullet((X\times_\bbF Y)_\proet, \Lambda).
\end{align*} 
That is, we consider the category of schemes as symmetric monoidal with respect to 
the fiber product over $\bbF$, and the external tensor product is natural on $X$ and $Y$ in the appropriate sense, see \cite[Section 3.1]{GaitsgoryLurie:Weil}, \cite[Section III.2]{GaitsgoryRozenblyum:StudyI} for details and precise statements. 
This functor $\boxtimes$ often fails to be an equivalence, so $\D_{\bullet, \Lambda}$ is not symmetric monoidal.
The assertion of \thref{derived_Drinfeld.text} is that this issue is resolved by replacing sheaves with Weil sheaves. 
In order to formulate \thref{derived_Drinfeld.text} as the monoidality of a certain functor, we need to replace the category of schemes by a category of objects that model Weil sheaves. 
We will represent these by taking the appropriate formal quotient by the partial Frobenius automorphism. Such formal quotients can be taken in the category of prestacks. 

We denote by $\PreStk_\bbF$ the category of (accessible) functors from the category $\CAlg_\bbF$ of commutative algebras
over $\bbF$ to the \ii-category $\Ani$ of Anima. 
The functor of taking points embeds the category of schemes fully faithfully into $\PreStk_\bbF$. 
We denote by
\eqn
\label{prestacks.sheaves.all}
  \D_{\bullet, \Lambda} \colon (\PreStk_\bbF)^{\mathrm{op}} \r \CatEx_{\infty, \Lambda}(\Idem)
\xeqn
the functor obtained by right Kan extension \cite[\textsection 4.3.2]{Lurie:Higher} along the inclusion $(\Schfp_\bbF)^{\mathrm{op}}\subset (\PreStk_\bbF)^{\mathrm{op}} $. Concretely, \cite[Proposition 6.2.1.9, Proposition 6.2.3.1]{Lurie:SAG}, given a prestack $Y$ which can be written as a colimit of schemes $Y_\alpha$ over some indexing category $A$ we have a canonical equivalence
\eqn
\D_\bullet(Y, \Lambda) \cong \lim_{\alpha} \D_\bullet(Y_\alpha, \Lambda).
\xeqn
This limit is formed in $\CatEx_{\infty, \Lambda}(\Idem)$; recall from around \refeq{CatExIdem.etc2} that the $\Ind$-completion functor to $\CatEx_{\infty, \Lambda}(\Idem) \r \PrStL$ does \emph{not} preserve (even finite) limits.

With this general sheaf theory in place, we can restrict our attention to the class of prestacks that is relevant to the derived Drinfeld lemma. 
\defi
Let $X$ be a scheme over $\Fq$. The \textit{Weil prestack} is defined as
$$X^\Weil := \colim \left (X \x_{\Fq} \bbF  \stackrel[\id]{\phi_X} \rightrightarrows X \x_{\Fq} \bbF \right ) \in \PreStk_\bbF,$$
i.e., it is the prestack sending $R\in \CAlg_\bbF$ to the colimit
\eqn
\label{Weil.prestack}
X^\Weil(R) = \colim \left (X(R) \stackrel[\id]{\phi_X}\rightrightarrows X(R) \right).
\xeqn
\xdefi
We denote by $\Schfp_\Weil$ the smallest full monoidal subcategory of $\PreStk_\bbF$ containing the Weil prestacks of finite type schemes $X / \Fq$. Equivalently, this is the full subcategory consisting of finite products of the form $X_1^\Weil \times \cdots \times X_n^\Weil$.
\lemm
Let $X_1,\ldots, X_n$ be schemes over $\Fq$. There is a canonical equivalence
\eqn
\label{Weil.prestack.colim}
\D_\bullet(X_1^\Weil \times_\bbF \cdots \times_\bbF X_n^\Weil)  \stackrel\cong\lr \Fix \big(\D_\bullet(X_{\bbF},\Lambda), \phi_{X_1}^*, \dots, \phi_{X_n}^*\big) 
\xeqn
\xlemm
\pf
Let $\Phi\colon \B \bbZ^n \r \PreStk_\bbF$ be the functor corresponding to the commuting automorphisms $\phi_{X_i}$. Then the claim follows immediately from the identification of $X_1^\Weil \times_\bbF \cdots \times_\bbF X_n^\Weil $ with the colimit of $\Phi$ (as an object in $\PreStk_\bbF$). 
\xpf

 \theo
 \thlabel{derived_Drinfeld.prestacks}
Suppose $\bullet$ and $\Lambda$ are as in \thref{derived_Drinfeld.text}.
Then the restriction of $\D_{\bullet, \Lambda}$ to Weil prestacks, i.e., the following composite
 \eqn
 \label{prestack.Drinfeld}
 \D_{\bullet, \Lambda}  \colon (\Schfp_\Weil)^{\mathrm{op}} \subset \PreStk_\bbF \r \CatEx_{\infty, \Lambda}(\Idem), 
\xeqn
is symmetric monoidal. 
\xtheo

\pf
As was noted above, the functor in \eqref{prestacks.sheaves.schemes} is lax symmetric monoidal.
By \cite[Proposition~2.7]{Torii:Perfect}, the Kan extension in \eqref{prestacks.sheaves.all} is still lax symmetric monoidal.
To check its restriction to the (symmetric monoidal) subcategory $\Schfp_\Weil$ is symmetric monoidal it suffices to show that the lax monoidal maps are in fact isomorphisms.
This is precisely the content of \thref{derived_Drinfeld.text}.
\xpf

\subsection{Full faithfulness}
\label{sect--sec.fully.faithfulness.Weil}
In this section, we prove that the functor \eqref{Drinfelds_map} is fully faithful under the conditions of \thref{derived_Drinfeld.text}.
We first consider the $p$-torsion free cases:

\prop
\thlabel{fully.faithul.q.invertible}
Let $X_1,\dots,X_n$ and $\Lambda$ be as in \thref{l.standard.situ.approximation} \refit{q.invertible.finite.Lambda}, \refit{q.invertible.adic.Lambda} or \refit{q.invertible.rational.Lambda}. 
Then the functor \eqref{Drinfelds_map} is fully faithful for $\bullet\in \{\lis,\cons\}$.
\xprop

\pf
For constructible sheaves on $X_{i,\bbF}$ (as opposed to $X_i^\Weil$), this interpretation of the Künneth formula appears already in \cite[Section A.2]{GKRV.Toy.Model}. 
Throughout, we drop $\Lambda$ from the notation.
It is enough to verify that for all $M_i, N_i \in \D_{\cons}(X_i^\Weil)$ the natural map
\eqn 
\label{tensor_equivalence}
\bigotimes_{i=1}^n\Hom_{\D(X_i^\Weil)}(M_i,N_i) \r \Hom_{\D(X_1^\Weil \x\ldots\x X_n^\Weil)}(M_1\bx\ldots\bx M_n,N_1\bx\ldots\bx N_n)
\xeqn 
is an equivalence. 
As \eqref{tensor_equivalence} is functorial in the objects and compatible with shifts, it suffices, by \thref{definition-Weil-on-products}, to consider the case where $M_i$, $i=1,\ldots,n$ is the extension by zero of a lisse Weil $\Lambda$-sheaf on some locally closed subscheme $Z_{i} \subset X_{i}$.
Using the adjunction 
$$(\iota_i)_! : \D_\cons(Z_i^\Weil) \rightleftarrows \D_\cons(X_i^\Weil) : (\iota_i)^!,$$
and the dualizability of lisse sheaves, we reduce to the case $M_i=\Lambda_{X_i}$, $i=1,\ldots,n$.
That is, \eqref{tensor_equivalence} becomes a map of cohomology complexes. 
By \thref{presentation.Weil.prop}, we have
$$\RG\big(X_i^\Weil, N_i\big) =  \mathrm{Fib} \left (\RG(X_{i,\bbF}, N_i) \stackrel{\phi_{X_i}^*-\id} \lr\RG(X_{i,\bbF}, N_i) \right).\eqlabel{mapping.bla}$$
A similar computation holds for the mapping complexes in $\D(X_1^\Weil\x\ldots \x X_n^\Weil)$, see \refeq{derived.cat.weil.sheaves.products.limit}.
Such finite limits commute with the tensor product in $\Mod_\Lambda$. 
Thus, \refeq{mapping.bla} reduces to the Künneth formula
\[
\RG\big(X_{1,\bbF},N_1\big)\t\ldots\t \RG\big(X_{n,\bbF},N_n\big) \overset{\cong}{\lr} \RG\big(X_{1,\bbF}\x_{\bbF}\ldots \x_{\bbF} X_{n,\bbF},N_1\bx\ldots \bx N_n\big),
\]
where we use that the $X_i$ are of finite type and the coprimality assumptions on $\Lambda$, see \StP{0F1P}. 
\xpf

Next, we consider the $p$-torsion case:

\prop
\thlabel{fully.faithfull.q.torsion}
Let $X_1,\dots,X_n$ and $\Lambda$ be as in \thref{l.standard.situ.approximation} \refit{p.torsion.Lambda}. 
Then the functor \eqref{Drinfelds_map} is fully faithful for $\bullet=\lis$.
\xprop

\pf
As in the proof of \thref{fully.faithul.q.invertible}, we need to show that the map
\eqn
\label{bla}
\bigotimes_{i=1}^n\RG\big(X_i^\Weil,N_i\big) \r \RG\big(X_1^\Weil \x\ldots\x X_n^\Weil, N_1\bx\ldots\bx N_n\big)
\xeqn
is an equivalence for any $N_i \in \D_\lis(X_i^\Weil)$.
Using Zariski descent for both sides, we may assume that each $X_i $ is affine.
As $\Lambda$ is finite discrete (see also the discussion around \eqref{Noetherian.approx.discrete.coefficients}), the invariance of the étale site under perfection reduces us to the case where each $X_i$ is perfect.
The proof now proceeds by several reduction steps: 1) reduce to $N_i = \Lambda_{X_i}$; 2) reduce to $\Lambda = \Z/p$; 3) reduce to $q = p$ being a prime. 
The last step 4) is then an easy computation.

\textbf{Step 1):} \textit{We may assume $N_i = \Lambda_{X_i}$.} 
In order to show \eqref{bla} is a quasi-isomorphism, it suffices to show this after applying $\tau^{\le r}$ for arbitrary $r$.
The complexes $N_i$ are bounded (\thref{sheaves}~\refit{constructible.bounded}).
By shifting them appropriately, we may assume $r = 0$.
Note that $\RG(X_i^\Weil,N_i)\cong \RG(X_i,N_i)$, see \thref{discrete_Geisser}.
By right exactness of the tensor product, we have $\tau^{\le 0} \left ( \bigotimes_i \RG(X_i, N_i) \right ) = \bigotimes_i \tau^{\le 0} \RG(X_i,N_i)$.
By the comparison with the classical notion of constructible sheaves (for discrete coefficients, see \consref{discrete.comparison.prop} and the discussion preceding it), there is an étale covering $U_i \r X_i$ such that $N_i |_{U_i}$ is perfect-constant.
Let $U_{i,\bullet}$ be the \v Cech nerve of this covering.
By étale descent, we have
$$\RG(X_i, N_i) = \lim_{[j]\in \Delta} \RG(U_{i,j},N_i).$$
For each $r \in \Z$, there is some $j_r$ such that 
$$\tau^{\le r} \lim_{[j] \in \Delta} \RG(U_{i,j}, N_i) = 
\lim_{[j] \in \Delta, j \le j_r} \tau^{\le r} \RG(U_{i,j}, N_i).$$
This can be seen from the spectral sequence
(note that it is concentrated in degrees $j \ge 0$ and degrees $j' \ge r$ for some $r$, since the complexes $N_i$ are bounded from below)
$$\H^{j'}(U_{i,j}, N_i) \Rightarrow \H^{j'+j} \lim_{j \in \Delta} \RG(U_{i,j}, N_i) =\H^{j'+j}(X_i, N_i).$$
As the tensor product in \eqref{bla} commutes with \emph{finite} limits, we may thus assume that each $N_i$ is perfect-constant.
Another dévissage reduces us to the case $N_i = \Lambda_{X_i}$, the constant sheaf itself.

\textbf{Step 2):} \textit{We may assume $\Lambda = \Z/p$.} 
By assumption, $\Lambda$ is flat over $\Z / p^m$ for some $m\geq 1$.
We immediately reduce to $\Lambda = \Z / p^m$.
For any perfect affine scheme $X=\Spec R$ in characteristic $p> 0$, we claim that $\RG(X,\bbZ/p^m)\t_{\bbZ/p^m}\bbZ/p^r\cong \RG(X,\bbZ/p^r)$.
Assuming the claim, we finish the reduction step by tensoring \eqref{bla} with the short exact sequence of $\bbZ/p^m$-modules $0 \r \Z/p^{m-1} \r \Z/p^m \r \Z/p \r 0$, using that finite limits commutes with tensor products. 
It remains to prove the claim.
The Artin--Schreier--Witt exact sequence of sheaves on $X_\et$ yields
$$\RG(X,\Z/p^m) = [W_m (R) \stackrel{F - \id} \r W_m (R)].$$
Now we use that $W_m (R) \t_{\Z/p^m} \Z/p^r \stackrel \cong \r W_r(R)$ compatibly with $F$, which holds since $R$ is perfect.
This shows the claim, and we have accomplished Step 2).

\textbf{Step 3):} \textit{We may assume $q$ is prime.} 
Recall that $q = p^r$ is a prime power. 
In order to reduce to the case $r=1$, let $X'_i := X_i$, but now regarded as a scheme over $\Fp$.
We have $X'_{i, \bbF} = \bigsqcup_{i = 1}^r X_{i, \bbF}$.
The Galois group $\Gal(\Fq / \Fp)$ is generated by the $p$-Frobenius, which acts by permuting the components in this disjoint union. 
Thus, we have $\D((X'_i)^\Weil) = \D(X_i^\Weil)$.
The same reasoning also applies to several factors $X_i^\Weil$, so we may assume our ground field to be $\Fp$.

\textbf{Step 4):} 
Set $R := \bigotimes_{i, \Fp} R_i$, $R_\bbF := R\t_{\Fp}\bbF$.
We write $\phi_i$ for the $p$-Frobenius on $R_i$ and also for any map on a tensor product involving $R_i$, by taking the identity on the remaining tensor factors.
By Artin--Schreier theory, we have 
$$
\eqalign{
\RG(X_i^\Weil,\bbZ/p) \stackrel{\text{\ref{discrete_Geisser}}} = \RG(X_i,\bbZ/p) &= [R_i \stackrel {\phi_i - \id} \r R_i], \cr 
\RG(X_{1, \bbF} \x_\bbF \dots \x_\bbF X_{n, \bbF}, \bbZ/p) & = [R_\bbF \stackrel{\phi - \id} \r R_\bbF],
}
$$
where $\phi$ is the absolute $p$-Frobenius of $R_\bbF$.
Thus, the right hand side in \eqref{bla} is the homotopy orbits of the action of $\Z^{n+1}$ on $R_\bbF$, whose basis vectors act as $\phi_1, \dots, \phi_n$ and $\phi$. 
Note that $\phi$ is the composite $\phi_\bbF \circ \phi_1 \circ \dots \circ \phi_n$,
where $\phi_\bbF$ is the Frobenius on $\bbF$.
Thus, the previously mentioned $\Z^{n+1}$-action on $R_\bbF$ is equivalent to the one where the basis vectors act as $\phi_1, \dots, \phi_n$ and $\phi_\bbF$.
We conclude our claim by using that $[R_\bbF \stackrel{\id - \phi_\bbF} \r R_\bbF]$ is quasi-isomorphic to $R[0]$.
\xpf


\subsection{Drinfeld's lemma}\label{sect--drinfelds.lemma}
The essential surjectivity in \thref{derived_Drinfeld.text} is based on the following variant of Drinfeld's lemma \cite[Theorem 2.1]{Drinfeld:Langlands} (see also \cite[IV.2, Theorem 4]{LaurentLafforgue:Chtoucas}, \cite[Theorem 8.1.4]{Lau:Shtukas}, \cite[Lemme 8.11]{Lafforgue:Chtoucas}, and \cite[Theorem 4.2.12]{Kedlaya:Shtukas}, \cite[Lemma 6.3]{Heinloth:Langlands}, \cite[Theorem 16.2.4]{ScholzeWeinstein:Berkeley} for expositions).
Its formulation is close to \cite[Theorem 8.1.4]{Lau:Shtukas}, and in this form is a slight extension of \cite[Lemme 8.2]{Lafforgue:Chtoucas} for $\bbZ_\ell$-coefficients and \cite[Lemma 3.3.2]{Xue:Finiteness} for $\bbQ_\ell$-coefficients.
We will drop the coefficient ring $\Lambda$ from the notation whenever convenient.

Let $X_1,\dots, X_n$ be Noetherian schemes over $\bbF_q$, and denote $X = X_1\times_{\bbF_q} \ldots \times_{\bbF_q} X_n$. 
Recall the Frobenius--Weil groupoid $\FWeil(X)$, see \thref{FWeil}. 
The projections $X_\bbF\r X_{i,\bbF}$ onto the single factors induce a continuous map of locally profinite groupoids
\eqn
\label{Weil.groupoids.map}
\mu\co \FWeil(X)\to \Weil(X_1)\x\ldots\x \Weil(X_n).
\xeqn

\theo[Version of Drinfelds's lemma]
\thlabel{Drinfelds_lemma}
Let $\Lambda$ be as in \thref{l.standard.situ.approximation}.
Restriction along the map \eqref{Weil.groupoids.map} induces an equivalence
\eqn
\label{eqn--Drinfelds_lemma_eq}
\Rep_\Lambda\big( \Weil(X_1)\x\ldots\x\Weil(X_n)\big) \stackrel \cong \r \Rep_\Lambda\big(\FWeil(X)\big),
\xeqn
between the abelian categories of continuous representations on finitely presented $\Lambda$-modules.
\xtheo
\pf
For all objects $x \in \FWeil(X)$, that is, all geometric points $x\to X_\bbF$, passing to the automorphism groups induces a commutative diagram of locally profinite groups
$$\xymatrix{
1 \ar[r] & \pi_1(X_\bbF,x) \ar[r] \ar[d] & \FWeil(X,x) \ar[r] \ar[d]^{\mu_x} & \bbZ^n \ar@{=}[d] \\
1 \ar[r] & \prod_{i=1}^n\pi_1(X_{i,\bbF},x) \ar[r] & \prod_{i=1}^n\Weil(X_{i},x) \ar[r] & \bbZ^n.
}$$
The left vertical arrow is surjective \StPd{0BN6}{0385}.
 Thus $\mu_x$ is surjective as well and hence \refeq{Drinfelds_lemma_eq} is fully faithful. 
For essential surjectivity, it remains to show that any continuous representation $\FWeil(X,x)\to \GL(M)$ on a finitely presented $\Lambda$-module $M$ factors through $\mu_x$.
The key input is Drinfeld's lemma: it implies that $\mu_x$ induces an isomorphism on profinite completions. 
Therefore, it is enough to apply \thref{bla.bla.blup.blup.lemma} below with $H:=\FWeil(X,x)\to \Weil(X_1)\x\ldots\x\Weil(X_n)=:G$ and $K:=\pi_1(X_\bbF,x)$.
This completes the proof of \refeq{Drinfelds_lemma_eq}.
\xpf

The following lemma formalizes a few arguments from \cite[\S 3.2.3]{Xue:Finiteness}, and we reproduce the proof for the convenience of the reader:

\lemm[Drinfeld, Xue]
\thlabel{bla.bla.blup.blup.lemma}
Let $\Lambda$ be as in \thref{l.standard.situ.approximation}.
Let $\mu\co H\to G$ be a continuous surjection of locally profinite groups that induces an isomorphism on profinite completions. 
Assume that there exists a compact open normal subgroup $K\subset H$ containing $\ker \mu$ such that $H/K$ is finitely generated and injects into its profinite completion. 
Then $\mu$ induces an equivalence
\[
\Rep_\Lambda(G) \cong \Rep_\Lambda(H)
\]
between their categories of continuous representations on finitely presented $\Lambda$-modules. 
\xlemm
\pf
The case where $\Lambda$ is finite discrete is obvious, and hence so is the case $\Lambda = \calO_E$ for some finite field extension $E\supset \Ql$. 
The case $\Lambda = E$ is reduced to $\Lambda = \Ql$.
As $\mu$ is surjective, it remains to show that every continuous representation $\rho\co H\to \GL(M)$ on a finite-dimensional $\bbQ_\ell$-vector space factors through $G$, that is, $\ker \mu\subset \ker \rho$.
One shows the following properties:
\begin{enumerate}
\item The group $\ker\mu$ is the intersection over all open subgroups in $K$ which are normal in $H$.
\item The group $\ker\rho\cap K$ is a closed normal subgroup in $H$ such that $K/\ker\rho\cap K\cong \rho(K)$ is topologically finitely generated. 
\end{enumerate}
These properties imply $\ker\mu\subset \ker\rho\cap K$ as follows:
For a finite group $L$, let $U_L:= \cap \ker(K\to L)$ where the intersection is over all continuous morphisms $K\to L$ that are trivial on $\ker\rho\cap K$.
Because of the topologically finitely generatedness in (2), this is a finite intersection so that $U_L$ is open in $K$.
Also, it is normal in $H$, and hence $\ker\mu \subset U_L$ by (1). 
On the other hand, it is evident that $\ker\rho\cap K=\cap_L U_L$ because $K$ is profinite. 

For the proof of (1) observe that $\ker\mu$ agrees with the kernel of $H\to H^\wedge\cong G^\wedge$ by our assumption on the profinite completions. 
Using $\ker \mu \subset K$ and the injection $H/K\to (H/K)^\wedge$ implies (1).

For (2) it is evident that $\ker\rho\cap K$ is a closed normal subgroup in $H$. 
Since $K$ is compact, its image $\rho(K)$ is a closed subgroup of the $\ell$-adic Lie group $\GL(M)$, hence an $\ell$-adic Lie group itself.
The final assertion follows from \cite[théorème 2]{Serre:GroupesAnal}.
\xpf

For the overall goal of proving essential surjectivity in \thref{derived_Drinfeld.text}, we need to investigate how representations of product groups factorize into external tensor products of representations. 
In view of \thref{elemntary.rep.lemm} and its proof, it is enough to consider representations of abstract groups, disregarding the topology. 
This is done in the next section.

\subsection{Factorizing representations}\label{Factorizing.representations}
In this subsection, let $\Lambda$ be a Dedekind domain \StP{034X}. 
Thus, 
any submodule $N$ of a finite projective $\Lambda$-module $M$ is again finite projective.

Given any group $W$, we write $\Rep^\fp_\Lambda(W)$ for the category of $W$-representations on finite projective $\Lambda$-modules.
As in \cite[Sections 73.8, 75]{CurtisReiner:RepresentationsOfFiniteGroups}, we say that such a $W$-representation $M$ is \emph{fp-simple} if any subrepresentation $0 \ne N \subset M$ has maximal rank.
By induction on the rank, every non-zero representation in $\Rep^\fp_\Lambda(W)$ admits a non-zero fp-simple subrepresentation. 
The proof of the following lemma is left to the reader. It parallels \cite[Theorem~75.6]{CurtisReiner:RepresentationsOfFiniteGroups}.

\lemm
\thlabel{Lambda.modules.reps}
A representation $M \in \Rep_\Lambda^\fp(W)$ is fp-simple if and only if $M \t_\Lambda \Frac(\Lambda)$ is fp-simple (hence, simple).
\xlemm


The following proposition will serve in the proof of \thref{derived_Drinfeld.text} using \thref{Drinfelds_lemma}, where we will need to decompose representations of a product of Weil groups into decompositions of the individual Weil groups.

\prop
\thlabel{factorization_lemm}
Let $W = W_1 \x W_2$ be a product of two groups. 
Let $M\in \Rep_\Lambda^\fp(W)$ be fp-simple.
Fix a $W_1$-subrepresentation $M_1 \subset M$ that is fp-simple.
Consider the $W_2$-representation $M_2 := \Hom_{W_1}(M_1, M)$ and the associated evaluation map 
$$\ev\co M_1 \bx M_2 \r M.$$
\begin{enumerate}
\item 
\label{item--factorization_lemm.1}
If $\Lambda$ is an algebraically closed field, then $\ev$ is an isomorphism and $M_2$ is simple.

\item 
\label{item--factorization_lemm.2}
If $\Lambda$ is a perfect field, then $\ev$ is a split surjection and $M_2$ is semi-simple.

\item 
\label{item--factorization_lemm.3}
If $\Lambda$ is a Dedekind domain of Krull dimension $1$ with perfect fraction field, then there is a short exact sequence
$$0 \r M \oplus \ker \ev \r M_1 \boxtimes M_2 \r T \r 0,
\eqlabel{extension.summand}$$
where $T$ is $\Lambda$-torsion.
\end{enumerate}
\xprop

\pf
Note that $\ev$ is a map in $\Rep_\Lambda^\fp(W)$. 
Its image has maximal rank by the fp-simplicity of $M$. 
Thus, if $\Lambda$ is a field, then it is surjective.

In case \refit{factorization_lemm.1}, we claim that $\ev$ is an isomorphism.
The following argument was explained to us by Jean-Fran\c{c}ois Dat:
For injectivity, observe that $M_1\bx M_2=M_1^{\oplus\dim M_2}$ as $W_1$-representations. 
Hence, if the kernel of $\ev$ is non-trivial, then it contains $M_1$ as an irreducible constituent. 
Therefore, it suffices to prove that $\Hom_{W_1}(M_1,\ev)$ is injective. 
Since $\Lambda$ is algebraically closed, we have $\End_{W_1}(M_1)=\Lambda$ by Schur's lemma.
Hence, the composition
\[
M_2=\Hom_{W_1}(\End_{W_1}(M_1),M_2)\cong\Hom_{W_1}(M_1,M_1\bx M_2) \to \Hom_{W_1}(M_1,M)=M_2
\]
is the identity. 
This shows that $\Hom_{W_1}(M_1,\ev)$ is an isomorphism.

In case \refit{factorization_lemm.2}, we claim that $M_1\bx M_2$ is semi-simple, and hence that $M$ appears as a direct summand.
Using \cite[Section 13.4 Corollaire]{Bourbaki:Algebre8} applied to the group algebras it is enough to show that $M_1$ and $M_2$ are absolutely semi-simple. 
Since $\Lambda$ is perfect, any finite-dimensional representation is semi-simple if and only if it is absolutely semi-simple, see \cite[Section 13.1]{Bourbaki:Algebre8}.
Hence, it remains to check that $M_{2,\bar\Lambda}=M_2\t_{\Lambda}\bar\Lambda$ is semi-simple where $\bar\Lambda/\Lambda$ is an algebraic closure.
The module $M_{2,\bar\Lambda}=\Hom_{W_1}(M_{1,\bar\Lambda},M_{\bar\Lambda})$ splits as a direct sum according to the simple constituents $\bar M_1\subset M_{1,\bar\Lambda}$ and $\bar M\subset M_{\bar\Lambda}$.
Finally, each $\bar M_2 =\Hom_{W_1}(\bar M_{1},\bar M)$ is either simple or vanishes: If there exists a non-zero $W_1$-equivariant map $\bar M_{1}\r \bar M$, then it must be injective by the simplicity of $\bar M_1$.
As $\bar\Lambda$ is algebraically closed, the proof of \refit{factorization_lemm.1} shows that $\bar M\cong \bar M_{1}\bx \bar M_2$ so that $\bar M_2$ must be simple because $\bar M$ is so. 
This shows that $M_2$ is absolutely semi-simple as well. 

In case \refit{factorization_lemm.3}, abbreviate $\Lambda' := \Frac \Lambda$, $M' := M \t_\Lambda \Lambda'$ and so on.
We will repeatedly use that $(\str)\t_\Lambda\Lambda'$ preserves and detects fp-simplicity of representations, see \thref{Lambda.modules.reps}. 
By \refit{factorization_lemm.2}, the evaluation map $\ev' := \ev \t \Lambda'$ admits a $\Lambda'$-linear section $\tilde i\co M' \r (M_1 \bx M_2)'$. As $M'$ is finitely presented, there is some $0 \ne \lambda \in \Lambda $ such that $\lambda \tilde i$ arises by scalar extension of a map $i\co M \r M_1 \bx M_2$. 
By construction, the map $i \oplus \on{incl}\co M \oplus \ker(\ev) \r M_1 \bx M_2$ is an isomorphism after tensoring with $\Lambda'$.
So its cokernel is $\Lambda$-torsion, and it is injective as both modules at the left are projective (hence $\Lambda$-torsion free).
This finishes the proof of the proposition. 
\xpf

\subsection{Essential surjectivity}\label{sect--section.essential.surjectivity}
In this section, we prove the essential surjectivity asserted in \thref{derived_Drinfeld.text}. 
Throughout, we freely use the full faithfulness proven in \thref{fully.faithul.q.invertible} and \thref{fully.faithfull.q.torsion}.

Recall that $X_1,\ldots, X_n$ are finite type $\Fq$-schemes, and write $X:=X_1\x_{\Fq}\ldots\x_{\Fq} X_n$.
Let $\Lambda$ be either a finite discrete ring, a finite field extension $E\supset \Ql$ for $\ell\neq p$ or its ring of integers $\calO_E$.
Note that this covers all cases from \thref{l.standard.situ.approximation}.

First, we show that it suffices to prove containment in the essential image \'{e}tale locally:

\lemm 
\thlabel{lemm.etale.descent.on.weil.equiv.site}
Let $U_i\rightarrow X_i$ be quasi-compact \'{e}tale surjections for $i=1,\ldots,n$. 
Then the following properties hold:
\begin{enumerate}
    \item \label{item--etale.descent.cons} 
    An object $M\in \D\big(X_1^\Weil \x \ldots\x X_n^{\Weil},\Lambda\big)$ belongs to the full subcategory 
    \[
    \D_\cons\big(X_1^\Weil,\Lambda\big)\t_{\Perf_{\Lambda_*}} \ldots \t_{\Perf_{\Lambda_*}}\D_\cons\big(X_n^\Weil,\Lambda\big) 
    \]
    if and only if its restriction $M|_{U_1^\Weil\x \ldots \x U_n^\Weil}$ belongs to the full subcategory 
    \[
    \D_\cons\big(U_1^\Weil,\Lambda\big) \t_{\Perf_{\Lambda_*}} \ldots \t_{\Perf_{\Lambda_*}}\D_\cons\big(U_n^\Weil,\Lambda\big) \subset
    \D\big(U_1^\Weil \x \ldots \x  U_n^\Weil,\Lambda\big).
    \]
    \item \label{item--etale.descent.lis}
    Assume that all $U_i\rightarrow X_i$ are finite \'{e}tale. 
    Then \refit{etale.descent.cons} holds for the categories of lisse sheaves. 
\end{enumerate}
\xlemm
\pf 
The only if direction in part \refit{etale.descent.cons} is clear. 
Conversely, assume that $M|_{U_1^\Weil\times \ldots \times U_n^\Weil}$ lies in the essential image of the external tensor product. 
By \'{e}tale descent, we have an equivalence:
\[
\D\big(X_1^\Weil\times \ldots\times X_n^\Weil,\Lambda\big) \stackrel\cong\lr \tot \left(\D\big(U_{1,\bullet}^\Weil\x \ldots \x U_{n,\bullet}^\Weil,\Lambda\big)\right).
\]
In particular, we get an equivalence $|(j_\bullet)_{!}\circ  j^*_\bullet M| \xrightarrow{\sim} M$ where $j_\bullet := j_{1,\bullet}\x\ldots\x j_{n,\bullet}$ with $j_{i,\bullet}\co U_{i,\bullet}\r X_i$ for $i=1,\ldots,n$. 
For each $m\geq 0$, the object $j_m^* M$ lies in 
\[
\D_\cons \big(U_{1,m}^\Weil,\Lambda\big)\t_{\Perf_{\Lambda_*}} \ldots \otimes_{\Perf_{\Lambda_*}} \D_\cons\big(U_{n,m}^\Weil,\Lambda\big). 
\]
It follows from \thref{sheaves}~\refit{preservation.constructibility} that these subcategories are preserved under $(j_m)_{!}$. 
So we see 
\[
(j_m)_{!}j_m^*(M)\in \D_\cons\big(X_1^\Weil,\Lambda\big)\otimes_{\Perf_{\Lambda_*}} \ldots \otimes_{\Perf_{\Lambda_*}}\D_\cons\big(X_n^\Weil,\Lambda\big) 
\]
for all $m\geq 0$. 
For every $m\geq 0$, let $M_m$ denote the realization of the $m$-th skeleton of the simplicial object $(j_\bullet)_{!}\circ  j^*_\bullet M$ so that we have a natural equivalence $\colim M_m \xrightarrow{\cong} M$ in $\D(X_1^\Weil\times \ldots\times X_n^\Weil,\Lambda)$. 
We claim that $M$ is a retract of some $M_m$, and hence lies in $\D_\cons\big(X_1^\Weil,\Lambda\big)\otimes_{\Perf_{\Lambda_*}} \ldots \otimes_{\Perf_{\Lambda_*}}\D_\cons\big(X_n^\Weil,\Lambda\big)$ by idempotent completeness.
To prove the claim, note that the sheaf $M_{\bbF}\in \D_\cons(X_\bbF,\Lambda)$ underlying $M$ is compact in the category of ind-constructible sheaves $\D_\indcons(X_\bbF,\Lambda)$, see \thref{sheaves}~\refit{compact.objects}.
As taking partial Frobenius fixed points is a finite limit, so commutes with filtered colimits, we see that the natural map of mapping complexes
\[
 \colim\Hom_{\D(X_1^\Weil \x \ldots\x X_n^{\Weil},\Lambda)}(M, M_m)\stackrel \cong\lr \Hom_{\D(X_1^\Weil \x \ldots\x X_n^{\Weil},\Lambda)}(M,\colim M_m)
\]
is an equivalence. 
In particular, the inverse equivalence $M \xrightarrow{\cong}\colim M_m$ factors through some $M_m$, presenting $M$ as a retract of $M_m$.
This proves the claim, and hence \refit{etale.descent.cons}.


For \refit{etale.descent.lis}, note that if $U_i \rightarrow X_i$ are finite \'{e}tale, then the functors $(j_{m})_!$ preserve the lisse categories, see \thref{sheaves}~\refit{preservation.constructibility}. 
In particular, for every $m\geq 0$ the object $(j_m)_{!}j_m^*(M)$ is lisse and so is $M_m$. 
We conclude using compactness as before. 
\xpf

Using \thref{ind.to.fixed.fully.faithful.big.category} and \thref{sheaves}~\refit{compact.objects}, the fully faithful functor \eqref{Drinfelds_map} uniquely extends to a fully faithful functor
\eqn\label{eqn:drinfeld-map-finite-ring-ind-cat}
\Ind\big(\D_\bullet(X^\Weil_{1},\Lambda)\big)\t_{\Mod_{\Lambda_*}}\ldots\t_{\Mod_{\Lambda_*}} \Ind\big(\D_{\bullet}(X^\Weil_{n},\Lambda)\big) \rightarrow \D(X_1^\Weil\times \ldots \times X_n^\Weil,\Lambda)
\xeqn
for $\bullet\in \{\lis,\cons\}$. 
We use this in the following variant of \thref{lemm.etale.descent.on.weil.equiv.site}. 

\lemm \thlabel{etale.descent.ind.small.Weil.categories}
The statements \refit{etale.descent.cons} and \refit{etale.descent.lis} of \thref{lemm.etale.descent.on.weil.equiv.site} hold for the functor \eqref{eqn:drinfeld-map-finite-ring-ind-cat} with $\bullet\in \{\lis,\cons\}$. 
Namely, to check that an object lies in the essential image of \eqref{eqn:drinfeld-map-finite-ring-ind-cat}, one can pass to a quasi-compact \'{e}tale cover if $\bullet = \cons$, and to a finite \'{e}tale cover if $\bullet = \lis$.
\xlemm
\pf
This is immediate from the proof of \thref{lemm.etale.descent.on.weil.equiv.site}: 
Arguing as above and using \'etale descent for ind-constructible, resp.~ind-lisse sheaves (\thref{sheaves}~\refit{descent}), we see that $M\cong \colim M_m$ with
\[
M_m\in \Ind\big(\D_\bullet(X^\Weil_{1},\Lambda)\big)\t_{\Mod_{\Lambda_*}}\ldots\t_{\Mod_{\Lambda_*}} \Ind\big(\D_{\bullet}(X^\Weil_{n},\Lambda)\big)
\]
for all $m\geq0 $ and $\bullet=\cons$, resp.~$\bullet=\lis$.
As the essential image of \eqref{eqn:drinfeld-map-finite-ring-ind-cat} is closed under colimits, $M$ lies in the corresponding subcategory as well. 
\xpf

Now we have enough tools to prove the categorical K\"unneth formula alias derived Drinfeld's lemma:

\pf [Proof of \thref{derived_Drinfeld.text}]
In view of \thref{fully.faithul.q.invertible} and \thref{fully.faithfull.q.torsion}, it remains to show the essential surjectivity of the external tensor product functor on Weil sheaves \eqref{Drinfelds_map} under the assumptions in \thref{derived_Drinfeld.text}. 
Part \refit{cons.Drinfeld}, the case of constructible sheaves, is reduced to part \refit{lisse.Drinfeld}, the case of lisse sheaves, by taking a stratification as in \thref{definition-Weil-on-products} \refit{definition-Weil-on-products.2} and using the full faithfulness already proven. 
Here we note that by refining the stratification witnessing the constructibility if necessary, we can even assume all strata to be smooth, so in particular geometrically unibranch. 
Hence, it remains to prove part \refit{lisse.Drinfeld}, that is, the essential surjectivity of the fully faithful functor
\eqn
\label{Drinfelds_map.proof}
\bx\co \D_\lis\big(X_1^\Weil,\Lambda\big)\t_{\Perf_{\Lambda_*}}\ldots\t_{\Perf_{\Lambda_*}} \D_\lis\big(X_n^\Weil,\Lambda\big) \r \D_\lis\big(X_1^\Weil\x\ldots\x X_n^\Weil,\Lambda\big),
\xeqn
when either $\Lambda$ is finite discrete as in cases \refit{q.invertible.finite.Lambda}, \refit{p.torsion.Lambda} in \thref{derived_Drinfeld.text} \refit{lisse.Drinfeld}, or $\Lambda=\calO_E$ for a finite field extension $E\supset \Ql$, $\ell\neq p$ as in case \refit{q.invertible.adic.Lambda}, or $\Lambda=E$ and the $X_i$ are geometrically unibranch as in the remaining case \refit{q.invertible.rational.Lambda}. 
In fact, the latter two cases are easier to handle due to the presence of natural t-structures on the categories of lisse sheaves (\thref{sheaves}~\refit{t-structure}).  
So we will distinguish two cases below: 1) $\Lambda=\calO_E$, or $\Lambda=E$ and all $X_i$ geometrically unibranch; 2) $\Lambda$ is finite discrete. 

Now pick $M\in \D_\lis(X_1^\Weil\x\ldots\x X_n^\Weil,\Lambda)$.
By \thref{sheaves}~\refit{constructible.bounded}, $M$ is bounded in the standard t-structure on $\D(X_1^\Weil\x\ldots\x X_n^\Weil,\Lambda)$.
So $M$ is a successive extension of its cohomology sheaves $\H^j(M)$, $j\in \bbZ$.
As $M$ is lisse, \thref{FWeil.heart.and.almost.compact.objects} \refit{weil.rep.essential.image} shows in both cases 1) and 2) that each $\H^j(M)$ comes from a continuous representation on a finitely presented $\Lambda$-module in
\eqn
\label{}
\Rep_{\Lambda}\big(\FWeil(X)\big) \stackrel{\ref{Drinfelds_lemma}} \cong \Rep_\Lambda( W),
\xeqn
where we denote $W:= W_1\x\ldots \x W_n$ with $W_i :=\Weil(X_i)$. 

Throughout, we repeatedly use that the functor \eqref{Drinfelds_map.proof} is fully faithful, commutes with finite (co-)limits and shifts, and that its essential image is closed under retracts (as the source category is idempotent complete, by definition) and contains all perfect-constant sheaves.

\textbf{Case 1):} \textit{Assume $\Lambda=\calO_E$, or $\Lambda=E$ and all $X_i$ geometrically unibranch.}
In this case, we have a t-structure on lisse Weil sheaves so that each $\H^j(M)$ belongs to $\D_\lis(X_1^\Weil\x\ldots\x X_n^\Weil,\Lambda)^\heartsuit$.
By induction on the length of $M$, using the full faithfulness of \eqref{Drinfelds_map.proof}, we reduce to the case where $M$ is abelian, that is, a continuous $W$-representation on a finitely presented $\Lambda$-module. 
The external tensor product induces a commutative diagram
\[
\xymatrix{
\Rep_\Lambda(W_1) \times \ldots \times \Rep_\Lambda(W_n) \ar[r]^{\quad\quad\quad\quad\boxtimes}\ar[d]^{\cong} & \Rep_\Lambda(W)\ar[d]^{\cong}\\
\D_\lis(X_1^{\Weil},\Lambda)^\heartsuit \times \ldots \times \D_\lis(X_n^{\Weil},\Lambda)^\heartsuit \ar[r]^{\quad\;\boxtimes} & \D_\lis(X_1^{\Weil}\times \ldots \times X_n^{\Weil},\Lambda)^\heartsuit,
}
\]
where the vertical equivalences are given by \thref{FWeil.heart.and.almost.compact.objects}.
Note that $M$ splits into a direct sum $M_{\on{tor}}\oplus M_\fp$ where the finitely presented $\Lambda$-module underlying $M_{\on{tor}}$ is $\Lambda$-torsion and $M_\fp$ is projective. 
So we can treat either case separately. 
Using that the essential image of \eqref{Drinfelds_map.proof} is closed under extensions (by full faithfulness) and retracts, the finite projective case is reduced to the fp-simple case and, by \thref{factorization_lemm}, to the finite torsion case.
Note that the $W_i$-representations constructed in, say \refeq{extension.summand}, are obtained from $M_\fp$ by taking subquotients and tensor products, so are automatically continuous.
Next, as the $\Lambda$-module underlying $M_{\on{tor}}$ is finite torsion, the $\Lambda$-sheaf $M_{\on{tor}}$ is perfect-constant along some finite \'{e}tale cover. 
So we conclude by \thref{lemm.etale.descent.on.weil.equiv.site} \refit{etale.descent.lis}. 

\textbf{Case 2):} \textit{Assume $\Lambda$ is finite discrete as above.}
In a nutshell, the argument is similar to the last step in case 1), but a little more involved due to the absence of natural t-structures on the categories of lisse sheaves in general, see \thref{sheaves}~\refit{t-structure} and \consref{Zl.no.t-structure}.
More precisely, in the special case, where $\Lambda$ is a finite field, the argument of case 1) applies, but not so if $\Lambda=\bbZ/\ell^2$, say. 
So, instead, we extend \eqref{Drinfelds_map.proof} by passing to Ind-completions to a commutative diagram
\[
\xymatrix{
\D_\lis\big(X_1^\Weil,\Lambda\big)\t_{\Perf_{\Lambda_*}}\ldots\t_{\Perf_{\Lambda_*}} \D_\lis\big(X_n^\Weil,\Lambda\big) \ar[r]^{\quad\quad\quad\quad\boxtimes}\ar[d] & \D_\lis\big(X_1^\Weil\x\ldots\x X_n^\Weil,\Lambda\big)\ar[d]\\
\Ind\big(\D_\lis(X^\Weil_{1},\Lambda)\big)\t_{\Mod_{\Lambda_*}}\ldots\t_{\Mod_{\Lambda_*}} \Ind\big(\D_{\lis}(X^\Weil_{n},\Lambda)\big) \ar[r]^{\quad\quad\quad\quad\quad\Ind(\boxtimes)} & \Ind\big(\D_\lis(X_1^\Weil\times \ldots \times X_n^\Weil,\Lambda)\big),
}
\] 
of full subcategories of $\D(X_1^\Weil\times \ldots \times X_n^\Weil,\Lambda)$, see the discussion around \eqref{eqn:drinfeld-map-finite-ring-ind-cat}.
Note that the fully faithful embedding \eqref{eqn:drinfeld-map-finite-ring-ind-cat} factors through $\Ind(\bx)$.
Both vertical arrows are the inclusion of the subcategories of compact objects by idempotent completeness of the involved categories and \refeq{Cat.idem.monoidal}.
Thus, if $M$ lies in the essential image of $\Ind(\bx)$, then it is a retract of a finite colimit of objects in the essential image of $\bx$, so lies itself in this essential image. 
As $M$ is a successive extension of its cohomology sheaves $\H^j(M)$, it suffices to show
\[
\H^j(M)\in \Ind\big(\D_\lis(X^\Weil_{1},\Lambda)\big)\t_{\Mod_{\Lambda_*}}\ldots\t_{\Mod_{\Lambda_*}} \Ind\big(\D_{\lis}(X^\Weil_{n},\Lambda)\big),
\]
for all $j\in \bbZ$.
So fix $j$ and denote $N:=\H^j(M)$ viewed as a continuous $W$-representation on a finitely presented $\Lambda$-module.
As $\Lambda$ is finite, $N$ comes from a continuous representation of $\pi_1(X_1)\times \ldots \pi_1(X_n)$ on which some open subgroup acts trivially. 
Hence, there exist finite \'{e}tale surjections $U_i\rightarrow X_i$ such that the subgroup $\pi_1(U_1)\times \ldots \times \pi_1(U_n)$ acts trivially on $N$. 
In particular, $N|_{U_1^{\Weil}\times \ldots \times U_n^{\Weil}}$ is constant, and hence lies in the essential image of the functor
\[
\Mod_{R}\cong\Ind\big(\Perf_R\big) \r \Ind\big(\D_\lis(U^\Weil_1 \x \ldots \x U^\Weil_n,\Lambda)\big),
\]
where $R:=\Gamma(\pi_0(U_1)\x \ldots \x\pi_0(U_n),\Lambda)$. 
As the sets $\pi_0(U_i)$ are finite discrete, each $R_i:=\Gamma(\pi_0(U_i),\Lambda)$ is a finite free $\Lambda_*$-algebra, and we have $R\cong R_1\t_{\Lambda_*}\ldots\t_{\Lambda_*}R_n$.
Thus, the external tensor product induces a commutative diagram
\[
\xymatrix{
\Mod_{R_1} \t_{\Mod_{\Lambda_*}} \ldots \t_{\Mod_{\Lambda_*}} \Mod_{R_n} \ar[r]^{\quad\quad\quad\quad\cong}\ar[d] & \Mod_{R}\ar[d]\\
\Ind\big(\D_\lis(U^\Weil_{1},\Lambda)\big)\t_{\Mod_{\Lambda_*}}\ldots\t_{\Mod_{\Lambda_*}} \Ind\big(\D_{\lis}(U^\Weil_{n},\Lambda)\big) \ar[r]^{\quad\quad\quad\quad\quad\Ind(\boxtimes)} & \Ind\big(\D_\lis(U_1^\Weil\times \ldots \times U_n^\Weil,\Lambda)\big),
}
\] 
where the upper horizontal arrow is an equivalence.
So $N|_{U_1^{\Weil}\times \ldots \times U_n^{\Weil}}$ lies in the essential image of $\Ind(\boxtimes)$, and we conclude by \thref{etale.descent.ind.small.Weil.categories} applied to the finite \'etale covers $U_i\r X_i$ and $\bullet=\lis$.
\xpf


\section{Ind-constructible Weil sheaves}\label{sec.ind.constructib.whole.section}
In this section, we introduce the full subcategories
\[
\D_\indlis\big(X^\Weil,\Lambda\big)\subset \D_\indcons\big(X^\Weil,\Lambda\big)
\]
of $\D\big(X^\Weil,\Lambda\big)$ consisting of ind-objects of lisse, resp.~constructible sheaves equipped with partial Frobenius action. 
That is, the partial Frobenius only preserves the ind-system of objects, but not necessarily each member. 
We will define analogous categories for a product of schemes. 
Similarly to the lisse, resp.~constructible case, there is a fully faithful functor
\[
\D_\indcons\big(X_1^\Weil,\Lambda\big)\t_{\Mod_{\Lambda_*}}\ldots\t_{\Mod_{\Lambda_*}} \D_\indcons\big(X_n^\Weil,\Lambda\big) \r \D_\indcons \big(X_1^\Weil\x\ldots\x X_n^\Weil,\Lambda\big),
\]
which, however, will not be an equivalence in general, see \thref{ind.cons.equivalence.example}. 
Nevertheless, we can identify a class of objects that lie in the essential image and that include many cases of interest such as the shtuka cohomology studied in \cite{Lafforgue:Chtoucas, LafforgueZhu:Decomposition, Xue:Finiteness, Xue:Smoothness}.

\subsection{Ind-constructible Weil sheaves}\label{sec.ind.cons.Weil}
Let $\Fq$ be a finite field of characteristic $p>0$, and fix an algebraic closure $\bbF$.
Let $X_1, \ldots, X_n$ be schemes of finite type over $\Fq$.  
Let $\Lambda$ be a condensed ring associated with the one of the following topological rings: a discrete coherent torsion ring (for example, a discrete finite ring), an algebraic field extension $E\supset \bbQ_\ell$, or its ring of integers $\calO_E$. 
We write $X := X_1 \x_{\Fq} \dots \x_{\Fq} X_n$, and denote by $X_{i, \bbF} := X_i \x_{\Fq} \Spec \bbF$ and $X_\bbF:=X\x_{\Fq} \Spec \bbF$ the base change.
Recall that under these assumptions, by \thref{sheaves}~\refit{compact.objects}, we have a fully faithful embedding
\eqn
\label{full.embedding.ind.cons}
\Ind\big(\D_\cons(X_{\bbF},\Lambda)\big)\stackrel\cong\lr  \D_{\indcons}(X_{\bbF},\Lambda) \subset \D(X_{\bbF},\Lambda),
\xeqn
and likewise for (ind-)lisse sheaves.
 
\defi
\thlabel{ind.lisse.cons.defi}
An object $M\in \D(X_1^\Weil\x\ldots\x X_n^\Weil,\Lambda)$ is called \textit{ind-lisse}, resp.~\textit{ind-constructible} if the underlying sheaf $M_{\bbF}\in \D(X_{\bbF},\Lambda)$ is ind-lisse, resp.~ind-constructible in the sense of \thref{recall.definitions}. 
\xdefi

We denote by 
\[
\D_{\indlis}\big(X_1^\Weil\x\ldots\x X_n^\Weil,\Lambda\big)\subset \D_{\indcons}\big(X_1^\Weil\x\ldots\x X_n^\Weil,\Lambda\big) 
\]
the resulting full subcategories of $\D(X_1^\Weil\x\ldots\x X_n^\Weil,\Lambda)$ consisting of ind-lisse, resp.~ind-constructible objects. 
Both categories are naturally commutative algebra objects in $\PrSt_{\Lambda_*}$ (see the notation from \refsect{recollections}), that is, presentable stable $\Lambda_*$-linear symmetric monoidal \ii-categories where $\Lambda_*:=\Gamma(*,\Lambda)$ is the ring underlying $\Lambda$. 
It is immediate from \thref{ind.lisse.cons.defi} that the equivalence \eqref{eqn--derived.cat.weil.sheaves.products.limit} restricts to an equivalence
\begin{align*}
 \D_{\bullet}\big(X_1^\Weil\x\ldots\x X_n^\Weil,\Lambda\big) &\cong \Fix \left(\D_{\bullet}(X_{\bbF},\Lambda), \phi_{X_1}^*, \dots, \phi_{X_n}^*\right)
\end{align*}
for $\bullet \in \{\indlis, \indcons\}$. 

\rema
Note that that we have a fully faithful embedding of $D_\cons(X^\Weil)$ into $\D_{\indcons}(X^\Weil)$ whose image consists of compact objects. However, the latter category is not generated by this image. Indeed, even
in the case of a point, the ind-cons category consists of $\Lambda$-modules with an action of an endomorphism, whereas the image of the embedding consists of $\Lambda$-modules with an action of an automorphism. This automorphism
does not have to fix any finitely generated submodule, which would be the case for any objects generated by the image of the constructible Weil complexes.
\xrema

Our goal in this section is to obtain a categorical Künneth formula for the categories of ind-lisse, resp.~ind-constructible Weil sheaves.  
In order to state the result, we need the following terminology.
Under our assumptions on $\Lambda$, each cohomology sheaf $\H^j(M)$, $j\in \bbZ$ for $M\in \D_{\lis} (X_{\bbF},\Lambda\big)$ is naturally a continuous representation of the pro\'etale fundamental groupoid $\pi_1^\proet(X_\bbF)$ on a finitely presented $\Lambda$-module, see \thref{FWeil.heart.and.almost.compact.objects}.
Further, the projections $X_{\bbF}\r X_{i,\bbF}$ induce a full surjective map of topological groupoids
\eqn
\label{pro.etale.fundmental.groupoids}
\pi_1^\proet(X_\bbF) \r \pi_1^\proet(X_{1,\bbF})\x\ldots\x \pi_1^\proet(X_{n,\bbF}).
\xeqn

\defi
\thlabel{split.ind.cons.definition}
Let $M\in \D(X_\bbF,\Lambda)$.
\begin{enumerate}
    \item 
    The sheaf $M$ is called \textit{split lisse} if it is lisse and the action of $\pi^{\proet}_1(X_{\bbF})$ on $\H^j(M)$ factors through \eqref{pro.etale.fundmental.groupoids} for all $j\in \bbZ$.
    
    \item 
     The sheaf $M$ is called \textit{split constructible} if it is constructible and there exists a finite subdivision into locally closed subschemes $X_{i, \alpha}\subseteq X_i$ such that for each $X_\alpha = \prod_i X_{i, \alpha} \subseteq X$, each restriction $M|_{X_{\alpha}}$ is split lisse.
\end{enumerate}
\xdefi

\defi 
\thlabel{split.weil.definition}
An object $M\in \D(X_1^\Weil\x \ldots \x X_n^\Weil,\Lambda)$ is called \textit{ind-(split lisse)}, resp.~\textit{ind-(split constructible)} if the underlying object $M_\bbF\in \D(X_\bbF,\Lambda)$ is a colimit of split lisse, resp.~split constructible objects.
\xdefi

As the category $\D_\bullet(X_\bbF,\Lambda)$, $\bullet\in\{\indlis, \indcons\}$ is cocomplete, every ind-(split lisse) object is ind-lisse, and likewise, every ind-(split constructible) object is ind-constructible.

\theo
\thlabel{split.objects.essential.image}
Assume that $\Lambda$ is either a finite discrete ring of prime-to-$p$ torsion, an algebraic field extension $E\supset \bbQ_\ell$ for $\ell\neq p$, or its ring of integers $\calO_E$. 
Then the functor induced by the external tensor product
\eqn\label{eqn:ind-lisse-external-functor}
\D_{\bullet} (X_1^{\Weil},\Lambda)\t_{\Mod_{\Lambda_*}}\ldots\t_{\Mod_{\Lambda_*}} \D_{\bullet} (X_n^{\Weil},\Lambda) \rightarrow
\D_{\bullet} (X_1^{\Weil}\times\ldots \times X_n^{\Weil},\Lambda)
\xeqn
is fully faithful for $\bullet \in \{\indlis,\  \indcons \}$. 
For $\bullet = \indlis$, resp.~$\bullet = \indcons$ the essential image contains the ind-(split lisse), resp.~ind-(split constructible) objects.
\xtheo

\pf
For full faithfulness, it is enough to consider the case $\bullet=\indcons$. 
Using \thref{fix.tensor.product.big.categories}, it remains to show that the functor
\eqn\label{eqn:functor-ind-lisse-underlying-tensor}
\bigotimes_i \D_{\indcons} \big(X_{i,\bbF},\Lambda) \cong \Ind\left( \bigotimes_i \Dcons(X_{i, \bbF},\Lambda) \right ) \rightarrow
\D_\bullet (X_{\bbF},\Lambda).
\xeqn
is fully faithful. 
In view of \eqref{full.embedding.ind.cons}, this is immediate from the K\"unneth formula for constructible $\Lambda$-sheaves as explained in \refsect{sec.fully.faithfulness.Weil}.

To identify objects in the essential image, we note that the fully faithful functors \eqref{eqn:ind-lisse-external-functor} and \eqref{eqn:functor-ind-lisse-underlying-tensor} induce a Cartesian diagram (see \thref{fix.tensor.product.big.categories}):
\eqn\label{eqn:cartesian-square-ind}
\xymatrix{
\bigotimes_i \D_{\bullet} \big(X_i^{\Weil},\Lambda) \ar[r]\ar[d] &
 \D_{\bullet}(X_1^{\Weil}\times\ldots \times X_n^{\Weil},\Lambda)\ar[d]\\
 \bigotimes_i \D_{\bullet}  (X_{i,\bbF},\Lambda) \ar[r] &
\D_{\bullet} (X_{\bbF},\Lambda),
}
\xeqn
for $\bullet \in \{\indlis,\  \indcons \}$. 
Thus, it is enough to show that the object $M_{\bbF}$ underlying an ind-split object $M$ lies in the image of the lower horizontal arrow.
Since this essential image is closed under colimits, it remains to show it contains the split lisse objects for $\bullet = \indlis$, resp.~the split constructible objects for $\bullet = \indcons$.

By the full faithfulness of \eqref{eqn:functor-ind-lisse-underlying-tensor}, the split constructible case reduces to the split lisse case, see also the proof of \thref{derived_Drinfeld.text} in \refsect{section.essential.surjectivity}.
So assume $\bullet=\indlis$ and let $M_\bbF\in \D(X_{\bbF},\Lambda)$ be split lisse. 
As each cohomology sheaf $\H^j(M_\bbF)$, $j\in\bbZ$ is at least ind-lisse (see also \consref{ind.cons.t.structure}), an induction on the cohomological length of $M_\bbF$ reduces us to show that $\H^j(M_\bbF)$ lies in the essential image.
By definition, being split lisse implies that the action of $\pi^{\proet}_1(X_\bbF)$ on $\H^j(M_\bbF)$ factors through $\pi^{\proet}_1(X_{1,\bbF}) \times \ldots \times \pi^{\proet}_1(X_{n,\bbF})$. 
Then the arguments of \refsect{section.essential.surjectivity} show that $\H^j(M_\bbF)$ lies is in the essential image of the lower horizontal arrow in \eqref{eqn:cartesian-square-ind}.
We leave the details to the reader.
\xpf

\rema
\thlabel{ind.cons.equivalence.example}
The functor \eqref{eqn:ind-lisse-external-functor} is not essentially surjective in general.
To see this, note that the functor $\D_{\indcons}(X^\Weil,\Lambda) \rightarrow \D_{\indcons}(X_{\bbF},\Lambda)$ admits a left adjoint $F$ that adds a free partial Frobenius action. 
Explicitly, for an object $M\in \D_{\indcons}(X_{\bbF},\Lambda)$ the object $F(M)$ has underlying sheaf $F(M)_{\bbF}$ given by a countable direct sum of copies of $M$. 
If $M$ was not originally in the image of the external tensor product (for example, $M$ as in \thref{AS.exam.intro}), then $F(M)$ will not be either. 
This is, however, the only obstacle for essential surjectivity: as noted in the proof of \thref{split.objects.essential.image}, the diagram \eqref{eqn:cartesian-square-ind} is Cartesian. 
\xrema 

\subsection{Cohomology of shtuka spaces}
\label{sect--shtuka.cohomology}
Finally, let us mention a key application of \thref{split.objects.essential.image}. 
Let $X$ be a smooth projective geometrically connected curve over $\bbF_q$. 
Let $N\subset X$ be a finite subscheme, and denote its complement by $Y=X\backslash N$.
Let $E\supset \bbQ_\ell$, $\ell\neq p$ be an algebraic field extension containing a fixed square root of $q$. 
Let $\calO_E$ be its ring of integers and denote by $k_E$ the residue field. 
Let $\Lambda$ be any of the topological rings $E,\calO_E,k_E$.
Let $G$ be a split (for simplicity) reductive group over $\bbF_q$.
We denote by $\widehat{G}$ the Langlands dual group of $G$ considered as a split reductive group over $\Lambda$.

In the seminal works \cite{Drinfeld:Langlands, LafforgueL:Chtoucas} ($G=\GL_n$) and \cite{Lafforgue:Chtoucas, LafforgueZhu:Decomposition} (general reductive $G$) on the Langlands correspondence over global function fields, the construction of the $\Weil(Y)$-action on automorphic forms of level $N$ is realized using the cohomology sheaves of moduli stacks of shtukas, defined in \cite{Varshavsky:Moduli} and \cite[Section 2]{Lafforgue:Chtoucas}. 
As explained in \cite{LafforgueZhu:Decomposition, GKRV.Toy.Model, zhu2021coherent}, the output of the geometric construction of Lafforgue can be encoded as a natural transformation
\eqn\label{eqn:drinfeld-lafforgue-xue}
\H_{N,I} \colon \Rep^{\fp}_{\Lambda}\big({\widehat{G}}^I\big) \to \Rep^{\cts}_{\Lambda}\big(\Weil(Y)^I\big), \quad  I\in \FinSet
\xeqn
of functors $\FinSet \to \Cat$ from the category of finite sets to the category of 1-categories.
Here the functor $\Rep^{\mathrm{\fp}}_{\Lambda}\big({\widehat{G}}^{\bullet}\big)$ assigns to a finite set $I$ the category of algebraic representations of ${\widehat{G}}^I$ on finite free $\Lambda$-modules, and $\Rep^{\cts}_{\Lambda}\big(\Weil(Y)^\bullet\big)$ the category of continuous representations of $\Weil(Y)^I$ in $\Lambda$-modules. 
In both cases, the transition maps are given by restriction of representations.

Let us recall some elements of its construction. 
For a finite set $I$,
\cite{Varshavsky:Moduli} and \cite[Section 2]{Lafforgue:Chtoucas} define the ind-algebraic stack $\Cht_{N,I}$ classifying $I$-legged $G$-shtukas on $X$ with full level-$N$-structure.
The morphism sending a $G$-shtuka to its legs
\eqn\label{eqn:shtukas.to.legs}
\mathfrak{p}_{N,I} \colon \Cht_{N,I} \rightarrow Y^I,
\xeqn
is locally of finite presentation. 
For every $W\in \Rep^{\fp}_{\Lambda}\big({\widehat{G}}^{I}\big)$, there is the normalized Satake sheaf $\calF_{N,I,W}$ on $\Cht_{N,I} $, see \cite[Définition 2.14]{Lafforgue:Chtoucas}. 
Base changing to $\bbF$ and taking compactly supported cohomology, we obtain the object
\[
\calH_{N,I}(W) \defined (\mathfrak{p}_{N,I,\bbF})_!(\calF_{N,I,W,\bbF}) \in \D_{\indcons}\big(Y^{I}_{\bbF},\Lambda\big),
\]
see \cite[D\'{e}finition 4.7]{Lafforgue:Chtoucas} and \cite[Definition 2.5.1]{Xue:Cuspidal}. 
Under the normalization of the Satake sheaves, the degree $0$ cohomology sheaf
\[
\H_{N,I} (W)\defined \H^0(\calH_{I}(W)) \in \D_{\indcons}\big(Y^{I}_{\bbF},\Lambda\big)^{\heartsuit}
\]
corresponds to the middle degree compactly supported intersection cohomology of $\Cht_{N,I}$.
Using the symmetries of the moduli stacks of shtukas, the sheaf $\H_{N,I}(W)$ is endowed with a partial Frobenius equivariant structure \cite[\textsection 6]{LafforgueL:Chtoucas}.
So we obtain objects
\eqn\label{eqn-laffogue-cohomology-sheaf}
\H_{N,I} (W) \in \D_{\indcons}\big((Y^\Weil)^{I},\Lambda\big)^{\heartsuit}.
\xeqn
Next, using the finiteness \cite{Xue:Finiteness} and smoothness \cite[Theorem 4.2.3]{Xue:Smoothness} results, the classical Drinfeld's lemma (\thref{Drinfelds_lemma}) applies to give objects $\H_{N,I}(W) \in \Rep^{\cts}_{\Lambda}\big(\Weil(Y)^I\big)$. 
The construction of the natural transformation \eqref{eqn:drinfeld-lafforgue-xue} encodes the functoriality and fusion satisfied by the objects $\{\H_{N,I}(W)\}$ for varying $I$ and $W$.

However, in order to analyze construction \eqref{eqn:drinfeld-lafforgue-xue} further,
it is desirable to upgrade the natural transformation of functors \eqref{eqn:drinfeld-lafforgue-xue} to the derived level. Namely, to have construction for the complexes $\{\calH_{I}(W)\}_{I,W}$ and not just for their cohomology sheaves, compare with \cite{zhu2021coherent}. 
Such an upgrade is possible using the derived version of Drinfeld's lemma, as given in the following proposition. 
A further study of this construction will appear in future work of the first named author (T.\ H.). 

\prop \thlabel{prop:derived.lafforgue.cohomologies}
For $\Lambda\in \{E,\calO_E,k_E\}$ and any $W\in \Rep_\Lambda(\widehat{G}^I)$, the shtuka cohomology \eqref{eqn-laffogue-cohomology-sheaf} lies in the essential image of the fully faithful functor
\eqn
\label{last.equation.finish}
\D_{\indlis}(Y^{\Weil},\Lambda)^{\otimes I} \r \D_{\indcons}\big((Y^\Weil)^I,\Lambda\big).
\xeqn
\xprop
\pf
By \cite[Theorem 4.2.3]{Xue:Smoothness}, the ind-constructible sheaf $\H_{N,I}(W)$ is ind-lisse.   
By \cite[Proposition 3.2.15]{Xue:Finiteness}, the action of $\FWeil(Y^I)$ on $\H_{N,I}(W)$ factors through the product $\Weil(Y)^I$. 
In particular, the action of $\pi_1(X_{\bbF}^{I})$ on $\H_{N,I}(W)$ factors through the product $\pi_1(X_{\bbF})^I$. 
So it is ind-(split lisse) in the sense of \thref{split.weil.definition}, and we are done by \thref{split.objects.essential.image}.
\xpf

\rema
One can upgrade the above construction in a homotopy coherent way to show that the whole complex $\calH_{N,I}(W)$ lies in $\D_{\indcons}\big((Y^\Weil)^{I},\Lambda\big)$.
If $N\not =\varnothing$ so that $\calH_{N,I}(W)$ is known to be bounded, then \thref{prop:derived.lafforgue.cohomologies} implies that $\calH_{N,I}(W)$ lies in the essential image of \eqref{last.equation.finish}. 
\xrema

\bibliographystyle{alphaurl}
\bibliography{bib}

\end{document}

%% file: constructible_citations.tex
\def\myconsref#1{\ifthenelse{\equal{#1}{adic.comparison.prop}}{Proposition~7.6}{}\ifthenelse{\equal{#1}{adjoints.closed.open.immersion.eq}}{3.3}{}\ifthenelse{\equal{#1}{all.colimit.commutation.lemm}}{Lemma~8.5}{}\ifthenelse{\equal{#1}{Artin.vanishing}}{Lemma~8.6}{}\ifthenelse{\equal{#1}{bounded.complexes.lemm}}{Lemma~4.10}{}\ifthenelse{\equal{#1}{bounded.complexes.nota}}{Definition~4.9}{}\ifthenelse{\equal{#1}{bounded.lisse.objects.lemm}}{Lemma~4.12}{}\ifthenelse{\equal{#1}{classical.comparison.colimit.lemm}}{Lemma~7.9}{}\ifthenelse{\equal{#1}{classical.comparison.lemm}}{Theorem~7.7}{}\ifthenelse{\equal{#1}{coherent.ring.Mod.fp.abelian}}{Lemma~6.5}{}\ifthenelse{\equal{#1}{cohomological.dimension.eq}}{8.1}{}\ifthenelse{\equal{#1}{cohomological.dimension.warn}}{Remark~5.4}{}\ifthenelse{\equal{#1}{cohomology.sheaves.colimit.eq}}{5.4}{}\ifthenelse{\equal{#1}{colimit.coefficients.lemm}}{Proposition~5.2}{}\ifthenelse{\equal{#1}{colimit.commutation.lemm}}{Lemma~5.3}{}\ifthenelse{\equal{#1}{compact.objects.intro}}{Proposition~1.8}{}\ifthenelse{\equal{#1}{compact.objects}}{Proposition~8.2}{}\ifthenelse{\equal{#1}{connected.component.map.eq}}{4.1}{}\ifthenelse{\equal{#1}{connected.component.pullback}}{Proposition~4.3}{}\ifthenelse{\equal{#1}{constant.comparison.eq}}{7.1}{}\ifthenelse{\equal{#1}{constant.objects.functor}}{3.1}{}\ifthenelse{\equal{#1}{constructible.bounded.coro}}{Corollary~4.11}{}\ifthenelse{\equal{#1}{ctf.comparison.rema}}{Remark~7.2}{}\ifthenelse{\equal{#1}{definition.lisse.constructible.introduction}}{Definition~1.1}{}\ifthenelse{\equal{#1}{descent.ind.cons}}{Corollary~8.7}{}\ifthenelse{\equal{#1}{descent.proof.of.condensed.shape}}{A.3}{}\ifthenelse{\equal{#1}{discrete.coefficients.lemm}}{Lemma~7.3}{}\ifthenelse{\equal{#1}{discrete.comparison.prop.proof}}{7.2}{}\ifthenelse{\equal{#1}{discrete.comparison.prop}}{Proposition~7.1}{}\ifthenelse{\equal{#1}{dualizable.objects.contractible.lem.intro}}{Lemma~1.2}{}\ifthenelse{\equal{#1}{dualizable.objects.contractible.lem}}{Lemma~4.1}{}\ifthenelse{\equal{#1}{eqn--Cat.idem.monoidal}}{2.2}{}\ifthenelse{\equal{#1}{eqn--CatExIdem.etc}}{2.1}{}\ifthenelse{\equal{#1}{etale.vs.proetale.fundamental.group}}{Remark~5.7}{}\ifthenelse{\equal{#1}{evaluation.formula.eq}}{4.2}{}\ifthenelse{\equal{#1}{evaluation.formula.eq2}}{4.3}{}\ifthenelse{\equal{#1}{examples.t.admissible.rings}}{Corollary~6.11}{}\ifthenelse{\equal{#1}{filtered.colimit.Hausdorff.rings}}{Lemma~5.6}{}\ifthenelse{\equal{#1}{filtered.colimits.limits.lemm}}{Lemma~2.3}{}\ifthenelse{\equal{#1}{formulas.closed.open.immersion.eq}}{3.4}{}\ifthenelse{\equal{#1}{fun.example}}{Example~3.6}{}\ifthenelse{\equal{#1}{Gaitsgory.Lurie.comparison}}{Corollary~8.3}{}\ifthenelse{\equal{#1}{global.sections.dualizable}}{3.2}{}\ifthenelse{\equal{#1}{ind.cons.t.structure}}{Remark~8.4}{}\ifthenelse{\equal{#1}{ind.lis.cons.defi.intro}}{Definition~1.7}{}\ifthenelse{\equal{#1}{ind.lis.cons.defi}}{Definition~8.1}{}\ifthenelse{\equal{#1}{insane.lattice.coro}}{Corollary~6.13}{}\ifthenelse{\equal{#1}{item--all.colimit.commutation.lemm.1}}{1}{}\ifthenelse{\equal{#1}{item--all.colimit.commutation.lemm.2}}{2}{}\ifthenelse{\equal{#1}{item--Artin.vanishing.1}}{1}{}\ifthenelse{\equal{#1}{item--Artin.vanishing.2}}{2}{}\ifthenelse{\equal{#1}{item--bounded.lisse.objects.lemm.1}}{1}{}\ifthenelse{\equal{#1}{item--bounded.lisse.objects.lemm.2}}{2}{}\ifthenelse{\equal{#1}{item--bounded.lisse.objects.lemm.3}}{3}{}\ifthenelse{\equal{#1}{item--colimit.commutation.lemm.1}}{1}{}\ifthenelse{\equal{#1}{item--colimit.commutation.lemm.2}}{2}{}\ifthenelse{\equal{#1}{item--connected.component.pullback.1}}{1}{}\ifthenelse{\equal{#1}{item--connected.component.pullback.2}}{2}{}\ifthenelse{\equal{#1}{item--connected.component.pullback.3}}{3}{}\ifthenelse{\equal{#1}{item--dualizable.objects.contractible.lem.1}}{1}{}\ifthenelse{\equal{#1}{item--dualizable.objects.contractible.lem.3}}{2}{}\ifthenelse{\equal{#1}{item--examples.t.admissible.rings.1}}{1}{}\ifthenelse{\equal{#1}{item--examples.t.admissible.rings.2}}{2}{}\ifthenelse{\equal{#1}{item--neighborhood.extension.coro.1}}{1}{}\ifthenelse{\equal{#1}{item--neighborhood.extension.coro.2}}{2}{}\ifthenelse{\equal{#1}{item--perfect.complexes.colimit.lemm.1}}{1}{}\ifthenelse{\equal{#1}{item--perfect.complexes.colimit.lemm.3}}{2}{}\ifthenelse{\equal{#1}{item--restriction.t.structures.1}}{1}{}\ifthenelse{\equal{#1}{item--restriction.t.structures.2}}{2}{}\ifthenelse{\equal{#1}{item--restriction.t.structures.3}}{3}{}\ifthenelse{\equal{#1}{item--restriction.t.structures.4}}{4}{}\ifthenelse{\equal{#1}{item--semi.hereditary.sections.combined.1}}{1}{}\ifthenelse{\equal{#1}{item--semi.hereditary.sections.combined.2}}{2}{}\ifthenelse{\equal{#1}{item--semi.hereditary.sections.combined.3}}{3}{}\ifthenelse{\equal{#1}{item--summary.sheaves.1}}{1}{}\ifthenelse{\equal{#1}{item--summary.sheaves.2}}{2}{}\ifthenelse{\equal{#1}{item--summary.sheaves.3}}{3}{}\ifthenelse{\equal{#1}{item--summary.sheaves.4}}{4}{}\ifthenelse{\equal{#1}{item--summary.sheaves.5}}{5}{}\ifthenelse{\equal{#1}{item--t.admissible.properties.1}}{2}{}\ifthenelse{\equal{#1}{item--t.admissible.properties.4}}{1}{}\ifthenelse{\equal{#1}{item--t.admissible.properties.5}}{3}{}\ifthenelse{\equal{#1}{item--t.structure.condensed.1}}{2}{}\ifthenelse{\equal{#1}{item--t.structure.condensed.2}}{3}{}\ifthenelse{\equal{#1}{limit.coefficients.lemm}}{Proposition~5.1}{}\ifthenelse{\equal{#1}{lisse.constructible.defi}}{Definition~3.3}{}\ifthenelse{\equal{#1}{lisse.constructible.hyperdescent.coro}}{Corollary~4.7}{}\ifthenelse{\equal{#1}{lisse.constructible.local.property.lemm}}{Lemma~4.5}{}\ifthenelse{\equal{#1}{lisse.evaluation.contractible.lemm}}{Corollary~4.8}{}\ifthenelse{\equal{#1}{lisse.sheaves.via.condensed.shape}}{Proposition~A.1}{}\ifthenelse{\equal{#1}{local.constancy.diagram.eq}}{4.4}{}\ifthenelse{\equal{#1}{localization.lisse.cons.lemm}}{Proposition~5.5}{}\ifthenelse{\equal{#1}{locally.constant.cor}}{Corollary~6.3}{}\ifthenelse{\equal{#1}{locally.constant.lemm}}{Corollary~7.4}{}\ifthenelse{\equal{#1}{locally.constant.prop}}{Theorem~4.13}{}\ifthenelse{\equal{#1}{neighborhood.extension.coro}}{Corollary~4.4}{}\ifthenelse{\equal{#1}{not.t.exact.example}}{Example~4.2}{}\ifthenelse{\equal{#1}{not.t.structure.rema}}{Remark~6.4}{}\ifthenelse{\equal{#1}{perfect.complexes.colimit.lemm}}{Lemma~2.2}{}\ifthenelse{\equal{#1}{perfect.lc}}{Definition~3.5}{}\ifthenelse{\equal{#1}{preservation.constructibility}}{Corollary~4.6}{}\ifthenelse{\equal{#1}{profinite.galois.equivalence}}{A.2}{}\ifthenelse{\equal{#1}{reduction.eq.2}}{5.2}{}\ifthenelse{\equal{#1}{reduction.eq}}{5.1}{}\ifthenelse{\equal{#1}{reduction.I}}{Lemma~7.5}{}\ifthenelse{\equal{#1}{restriction.t.structures.eq}}{6.2}{}\ifthenelse{\equal{#1}{restriction.t.structures}}{Proposition~6.6}{}\ifthenelse{\equal{#1}{sect--adic.coefficients}}{7.1.1}{}\ifthenelse{\equal{#1}{sect--coefficients.colimits.sect}}{5.2}{}\ifthenelse{\equal{#1}{sect--coefficients.limits.sect}}{5.1}{}\ifthenelse{\equal{#1}{sect--coefficients.localizations.sect}}{5.3}{}\ifthenelse{\equal{#1}{sect--coefficients.sect}}{5}{}\ifthenelse{\equal{#1}{sect--comparison.results}}{7}{}\ifthenelse{\equal{#1}{sect--condensed.shape}}{A}{}\ifthenelse{\equal{#1}{sect--discrete.coefficients}}{7.1}{}\ifthenelse{\equal{#1}{sect--examples.comparison.l.adic}}{7.2}{}\ifthenelse{\equal{#1}{sect--generalities.proetale}}{3}{}\ifthenelse{\equal{#1}{sect--hyperdescent}}{4.2}{}\ifthenelse{\equal{#1}{sect--ind.cons.sheaves}}{8}{}\ifthenelse{\equal{#1}{sect--limits.filtered.colimits}}{2.2}{}\ifthenelse{\equal{#1}{sect--lisse.locally.connected.sect}}{4.4}{}\ifthenelse{\equal{#1}{sect--monoidal.aspects}}{2.1}{}\ifthenelse{\equal{#1}{sect--prelude}}{2}{}\ifthenelse{\equal{#1}{sect--t.structures.sect}}{6}{}\ifthenelse{\equal{#1}{semi.hereditary.sections.combined}}{Theorem~6.12}{}\ifthenelse{\equal{#1}{semi.hereditary.t.admissible}}{Lemma~6.10}{}\ifthenelse{\equal{#1}{shape.w-contractible.cover}}{A.1}{}\ifthenelse{\equal{#1}{six.functor.rema}}{Remark~7.10}{}\ifthenelse{\equal{#1}{split.sequence.eq}}{6.3}{}\ifthenelse{\equal{#1}{summary.sheaves}}{Theorem~1.6}{}\ifthenelse{\equal{#1}{t-admissible}}{Definition~1.4}{}\ifthenelse{\equal{#1}{t.admissible.defi}}{Definition~6.1}{}\ifthenelse{\equal{#1}{t.admissible.properties}}{Lemma~6.7}{}\ifthenelse{\equal{#1}{t.exactness.lemm}}{Lemma~3.4}{}\ifthenelse{\equal{#1}{t.structure.condensed}}{Theorem~6.2}{}\ifthenelse{\equal{#1}{t.structure.discrete}}{Proposition~6.8}{}\ifthenelse{\equal{#1}{t.structures.eq}}{6.1}{}\ifthenelse{\equal{#1}{tensor.extension.by.zero.eq}}{5.3}{}\ifthenelse{\equal{#1}{Zl.no.t-structure}}{Remark~6.9}{}}

%% file: kuenneth.bbl
\begin{thebibliography}{GKRV22}

\bibitem[BL19]{BhattLurie:Riemann}
Bhargav Bhatt and Jacob Lurie.
\newblock A {R}iemann-{H}ilbert correspondence in positive characteristic.
\newblock {\em Camb. J. Math.}, 7(1-2):71--217, 2019.
\newblock \href {https://doi.org/10.4310/CJM.2019.v7.n1.a3}
  {\path{doi:10.4310/CJM.2019.v7.n1.a3}}.

\bibitem[BM21]{BhattMathew:Arc}
Bhargav Bhatt and Akhil Mathew.
\newblock The arc-topology.
\newblock {\em Duke Math. J.}, 170(9):1899--1988, 2021.
\newblock \href {https://doi.org/10.1215/00127094-2020-0088}
  {\path{doi:10.1215/00127094-2020-0088}}.

\bibitem[Bou12]{Bourbaki:Algebre8}
N.~Bourbaki.
\newblock {\em \'{E}l\'{e}ments de math\'{e}matique. {A}lg\`ebre. {C}hapitre 8.
  {M}odules et anneaux semi-simples}.
\newblock Springer, Berlin, 2012.
\newblock Second revised edition of the 1958 edition [MR0098114].
\newblock \href {https://doi.org/10.1007/978-3-540-35316-4}
  {\path{doi:10.1007/978-3-540-35316-4}}.

\bibitem[BS15]{BhattScholze:ProEtale}
Bhargav Bhatt and Peter Scholze.
\newblock The pro-\'{e}tale topology for schemes.
\newblock {\em Ast\'{e}risque}, (369):99--201, 2015.

\bibitem[BZFN10]{BenZviFrancisNadler:Integral}
David Ben-Zvi, John Francis, and David Nadler.
\newblock Integral transforms and {D}rinfeld centers in derived algebraic
  geometry.
\newblock {\em J. Amer. Math. Soc.}, 23(4):909--966, 2010.
\newblock \href {https://doi.org/10.1090/S0894-0347-10-00669-7}
  {\path{doi:10.1090/S0894-0347-10-00669-7}}.

\bibitem[CR06]{CurtisReiner:RepresentationsOfFiniteGroups}
Charles~W. Curtis and Irving Reiner.
\newblock {\em Representation theory of finite groups and associative
  algebras}.
\newblock AMS Chelsea Publishing, Providence, RI, 2006.
\newblock Reprint of the 1962 original.
\newblock \href {https://doi.org/10.1090/chel/356}
  {\path{doi:10.1090/chel/356}}.

\bibitem[Del80]{Deligne:Weil2}
Pierre Deligne.
\newblock La conjecture de {W}eil. {II}.
\newblock {\em Inst. Hautes \'Etudes Sci. Publ. Math.}, 52:137--252, 1980.

\bibitem[Dri80]{Drinfeld:Langlands}
V.~G. Drinfeld.
\newblock Langlands' conjecture for {${\rm GL}(2)$} over functional fields.
\newblock In {\em Proceedings of the {I}nternational {C}ongress of
  {M}athematicians ({H}elsinki, 1978)}, pages 565--574. Acad. Sci. Fennica,
  Helsinki, 1980.

\bibitem[Dri87]{Drinfeld:FSheaves}
V.~G. Drinfeld.
\newblock Moduli varieties of {$F$}-sheaves.
\newblock {\em Funktsional. Anal. i Prilozhen.}, 21(2):23--41, 1987.

\bibitem[Gei04]{Geisser:Weil}
Thomas Geisser.
\newblock Weil-\'etale cohomology over finite fields.
\newblock {\em Math. Ann.}, 330(4):665--692, 2004.
\newblock \href {https://doi.org/10.1007/s00208-004-0564-8}
  {\path{doi:10.1007/s00208-004-0564-8}}.

\bibitem[GKRV22]{GKRV.Toy.Model}
Dennis Gaitsgory, David Kazhdan, Nick Rozenblyum, and Yakov Varshavsky.
\newblock A toy model for the {D}rinfeld--{L}afforgue shtuka construction.
\newblock {\em Indagationes Mathematicae}, 33(1):39--189, 2022.
\newblock \href {http://arxiv.org/abs/1908.05420} {\path{arXiv:1908.05420}}.

\bibitem[GL19]{GaitsgoryLurie:Weil}
Dennis Gaitsgory and Jacob Lurie.
\newblock {\em Weil's conjecture for function fields. {V}ol. 1}, volume 199 of
  {\em Annals of Mathematics Studies}.
\newblock Princeton University Press, Princeton, NJ, 2019.

\bibitem[GR17]{GaitsgoryRozenblyum:StudyI}
Dennis Gaitsgory and Nick Rozenblyum.
\newblock {\em A study in derived algebraic geometry. {V}ol. {I}.
  {C}orrespondences and duality}, volume 221 of {\em Mathematical Surveys and
  Monographs}.
\newblock American Mathematical Society, Providence, RI, 2017.

\bibitem[Hei18]{Heinloth:Langlands}
Jochen Heinloth.
\newblock Langlands parameterization over function fields following {V}.
  {L}afforgue.
\newblock {\em Acta Math. Vietnam.}, 43(1):45--66, 2018.

\bibitem[HRS23]{HemoRicharzScholbach:Constructible}
Tamir Hemo, Timo Richarz, and Jakob Scholbach.
\newblock {Constructible sheaves on schemes}, 2023.
\newblock \href {http://arxiv.org/abs/2305.18131} {\path{arXiv:2305.18131}}.

\bibitem[HS21]{HansenScholze:RelativePerversity}
David Hansen and Peter Scholze.
\newblock Relative perversity, 2021.
\newblock \href {http://arxiv.org/abs/2109.06766} {\path{arXiv:2109.06766}}.

\bibitem[Kah03]{Kahn:Equivalences}
Bruno Kahn.
\newblock \'{E}quivalences rationnelle et num\'{e}rique sur certaines
  vari\'{e}t\'{e}s de type ab\'{e}lien sur un corps fini.
\newblock {\em Ann. Sci. \'{E}cole Norm. Sup. (4)}, 36(6):977--1002 (2004),
  2003.
\newblock \href {https://doi.org/10.1016/j.ansens.2003.02.002}
  {\path{doi:10.1016/j.ansens.2003.02.002}}.

\bibitem[Ked19]{Kedlaya:Shtukas}
Kiran~S. Kedlaya.
\newblock {\em Sheaves, Stacks, and Shtukas}, volume 242 of {\em in Perfectoid
  spaces: Mathematical Surveys and Monographs}.
\newblock American Mathematical Society, Providence, RI, 2019.
\newblock Lectures from the 2017 Arizona Winter School, held in Tucson, AZ,
  March 11--17, Edited and with a preface by Bryden Cais, With an introduction
  by Peter Scholze.
\newblock \href {https://doi.org/10.1090/surv/242}
  {\path{doi:10.1090/surv/242}}.

\bibitem[Laf97]{LaurentLafforgue:Chtoucas}
Laurent Lafforgue.
\newblock Chtoucas de {D}rinfeld et conjecture de {R}amanujan-{P}etersson.
\newblock {\em Ast\'{e}risque}, (243):ii+329, 1997.

\bibitem[Laf02]{LafforgueL:Chtoucas}
Laurent Lafforgue.
\newblock Chtoucas de {D}rinfeld et correspondance de {L}anglands.
\newblock {\em Invent. Math.}, 147(1):1--241, 2002.

\bibitem[Laf18]{Lafforgue:Chtoucas}
Vincent Lafforgue.
\newblock Chtoucas pour les groupes r\'{e}ductifs et param\'{e}trisation de
  {L}anglands globale.
\newblock {\em J. Amer. Math. Soc.}, 31(3):719--891, 2018.

\bibitem[Lau04]{Lau:Shtukas}
Eike~S\"{o}ren Lau.
\newblock {\em On generalised {$D$}-shtukas}, volume 369 of {\em Bonner
  Mathematische Schriften [Bonn Mathematical Publications]}.
\newblock Universit\"{a}t Bonn, Mathematisches Institut, Bonn, 2004.
\newblock Dissertation, Rheinische Friedrich-Wilhelms-Universit\"{a}t Bonn,
  Bonn, 2004.

\bibitem[Lic05]{Lichtenbaum:Weil}
S.~Lichtenbaum.
\newblock The {W}eil-\'etale topology on schemes over finite fields.
\newblock {\em Compos. Math.}, 141(3):689--702, 2005.
\newblock \href {https://doi.org/10.1112/S0010437X04001150}
  {\path{doi:10.1112/S0010437X04001150}}.

\bibitem[Lur09]{Lurie:Higher}
Jacob Lurie.
\newblock {\em Higher topos theory}, volume 170 of {\em Annals of Mathematics
  Studies}.
\newblock Princeton University Press, Princeton, NJ, 2009.

\bibitem[Lur17]{Lurie:HA}
Jacob Lurie.
\newblock {Higher Algebra}, 2017.
\newblock URL: \url{http://www.math.harvard.edu/~lurie/}.

\bibitem[Lur18]{Lurie:SAG}
Jacob Lurie.
\newblock {Spectral Algebraic Geometry}, 2018.
\newblock URL: \url{http://www.math.harvard.edu/~lurie/}.

\bibitem[LZ19]{LafforgueZhu:Decomposition}
Vincent Lafforgue and Xinwen Zhu.
\newblock D\'ecomposition au-dessus des param\`etres de langlands elliptiques,
  2019.
\newblock \href {http://arxiv.org/abs/1811.07976} {\path{arXiv:1811.07976}}.

\bibitem[Ser64]{Serre:GroupesAnal}
Jean-Pierre Serre.
\newblock Groupes analytiques {$p$}-adiques.
\newblock In {\em S\'{e}minaire {B}ourbaki, 1963/64, {F}asc. 2, {E}xpos\'{e}
  270}, page~10. Secr\'{e}tariat math\'{e}matique, Paris, 1964.

\bibitem[SGA03]{SGA1}
{\em Rev\^etements \'etales et groupe fondamental ({SGA} 1)}.
\newblock Documents Math\'ematiques (Paris) [Mathematical Documents (Paris)],
  3. Soci\'et\'e Math\'ematique de France, Paris, 2003.
\newblock S\'eminaire de g\'eom\'etrie alg\'ebrique du Bois Marie 1960--61.
  [Algebraic Geometry Seminar of Bois Marie 1960-61], Directed by A.
  Grothendieck, With two papers by M. Raynaud, Updated and annotated reprint of
  the 1971 original [Lecture Notes in Math., 224, Springer, Berlin; MR0354651
  (50 \#7129)].

\bibitem[{Sta}17]{StacksProject}
{Stacks Project Authors}.
\newblock {Stacks Project}.
\newblock \url{http://stacks.math.columbia.edu}, 2017.

\bibitem[SW20]{ScholzeWeinstein:Berkeley}
Peter Scholze and Jared Weinstein.
\newblock {\em Berkeley lectures on {{\(p\)}}-adic geometry}, volume 207 of
  {\em Ann. Math. Stud.}
\newblock Princeton, NJ: Princeton University Press, 2020.
\newblock \href {https://doi.org/10.1515/9780691202150}
  {\path{doi:10.1515/9780691202150}}.

\bibitem[Tor22]{Torii:Perfect}
Takeshi Torii.
\newblock A perfect pairing for monoidal adjunctions, 2022.
\newblock \href {http://arxiv.org/abs/arXiv:2202.02493}
  {\path{arXiv:arXiv:2202.02493}}.

\bibitem[{Var}04]{Varshavsky:Moduli}
Yakov {Varshavsky}.
\newblock {Moduli spaces of principal $F$-bundles.}
\newblock {\em {Sel. Math., New Ser.}}, 10(1):131--166, 2004.

\bibitem[Xue20a]{Xue:Cuspidal}
Cong Xue.
\newblock Cuspidal cohomology of stacks of shtukas.
\newblock {\em Compos. Math.}, 156(6):1079--1151, 2020.
\newblock \href {https://doi.org/10.1112/s0010437x20007058}
  {\path{doi:10.1112/s0010437x20007058}}.

\bibitem[Xue20b]{Xue:Finiteness}
Cong Xue.
\newblock Finiteness of cohomology groups of stacks of shtukas as modules over
  {H}ecke algebras, and applications.
\newblock {\em \'{E}pijournal G\'{e}om. Alg\'{e}brique}, 4:Art. 6, 42 pp.--41,
  2020.

\bibitem[Xue20c]{Xue:Smoothness}
Cong Xue.
\newblock Smoothness of cohomology sheaves of stacks of shtukas, 2020.
\newblock \href {http://arxiv.org/abs/2012.12833} {\path{arXiv:2012.12833}}.

\bibitem[Zhu21]{zhu2021coherent}
Xinwen Zhu.
\newblock Coherent sheaves on the stack of {Langlands} parameters, 2021.
\newblock \href {http://arxiv.org/abs/2008.02998} {\path{arXiv:2008.02998}}.

\end{thebibliography}
